\UseRawInputEncoding
\documentclass[12pt, reqno]{amsart}
\usepackage[margin=1in]{geometry}
\usepackage{amssymb,latexsym,amsmath,amscd,amsfonts}
\usepackage{latexsym}
\usepackage[mathscr]{eucal}
\usepackage{bm}
\usepackage{tabularx}
\usepackage{mathptmx}
\usepackage{amssymb}
\usepackage{amsthm}
\usepackage{dcolumn}
\usepackage[all]{xy}
\usepackage{enumitem}
\usepackage[utf8]{inputenc}
\usepackage{mathtools}

\def \qed {\hfill \vrule height6pt width 6pt depth 0pt}
\def\textmatrix#1&#2\\#3&#4\\{\bigl({#1 \atop #3}\ {#2 \atop #4}\bigr)}
\def\dispmatrix#1&#2\\#3&#4\\{\left({#1 \atop #3}\ {#2 \atop #4}\right)}
\newcommand{\beg}{\begin{equation}}
	\newcommand{\eeg}{\end{equation}}
\newcommand{\ben}{\begin{eqnarray*}}
	\newcommand{\een}{\end{eqnarray*}}

\newcommand{\Q}{\mathbb{H}}
\newcommand{\CQ}{\overline{\mathbb{H}}}
\newcommand{\C}{\mathbb{C}}

\newcommand{\D}{\mathbb{D}}
\newcommand{\DC}{\overline{\mathbb{D}}}
\newcommand{\T}{\mathbb{T}}

\newcommand{\Gg}{\mathbb{G}_2}

\newcommand{\E}{\mathbb E}
\newcommand{\Ebar}{\overline{\mathbb E}}

\newcommand{\lm}{\lambda}

\newcommand{\al}{\alpha}
\newcommand{\be}{\beta}

\newcommand{\HB}{\mathbb{H}_\mu}
\newcommand{\CHB}{\overline{\mathbb{H}}_\mu}
\newcommand{\QN}{\mathbb{H}_{N}}
\newcommand{\CQN}{\overline{\mathbb{H}}_{N}}
\newcommand{\Pe}{\mathbb{P}}
\newcommand{\Pbar}{\overline{\mathbb{P}}}
\newcommand{\Pz}{\psi_{z_1, z_2}}

\newtheorem{thm}{Theorem}[section]
\newtheorem{cor}[thm]{Corollary}
\newtheorem{lem}[thm]{Lemma}

\newtheorem{prop}[thm]{Proposition}
\numberwithin{equation}{section} \theoremstyle{definition}
\newtheorem{defn}[thm]{Definition}

\newtheorem{eg}[thm]{Example}

\def\textmatrix#1&#2\\#3&#4\\{\bigl({#1 \atop #3}\ {#2 \atop #4}\bigr)}
\def\dispmatrix#1&#2\\#3&#4\\{\left({#1 \atop #3}\ {#2 \atop #4}\right)}
\begin{document}
		\title{The Hexablock: a domain associated with the $\mu$-synthesis in $M_2(\mathbb C)$}

	\author[I. Biswas]{Indranil Biswas}
	
\address{Department of Mathematics,
Shiv Nadar University, NH91, Tehsil Dadri, Greater Noida,
Uttar Pradesh 201314, India.}

\email{indranil.biswas@snu.edu.in, indranil29@gmail.com}
	
\author[S. Pal]{Sourav Pal}

\address{Mathematics Department, Indian Institute of Technology Bombay,
Powai, Mumbai - 400076, India.}

\email{sourav@math.iitb.ac.in, souravmaths@gmail.com}
	
\author[N. Tomar]{Nitin Tomar}

\address{Mathematics Department, Indian Institute of Technology Bombay, Powai, Mumbai-400076, India.}

\email{tomarnitin414@gmail.com, tnitin@math.iitb.ac.in}		
	
	
	
\begin{abstract}
We introduce a domain named \textit{hexablock} in $\mathbb C^4$ and show that its origin is a 
special case of $\mu$-synthesis in $M_2(\mathbb C)$, more precisely the $\mu_E$-unit ball with respect to the 
linear subspace $E$ consisting of $2 \times 2$ upper triangular matrices. The hexablock is denoted by $\mathbb 
H$ and is defined by
\[
\mathbb{H}=\left\{(a, x_1, x_2, x_3) \,\in\, \mathbb{C} \times \mathbb{E}\,\,\big\vert\,\,
\sup_{z_1,\, z_2 \,\in\,
\D}\left|\frac{a\sqrt{(1-|z_1|^2)(1-|z_2|^2)}}{1-x_1z_1-x_2z_2+x_3z_1z_2}\right|
<1\right\},
\]
where $\mathbb{E}$ is the \textit{tetrablock}, another domain in $\mathbb C^3$ associated with a
different case of $\mu$-synthesis, and is given by 
\[
\mathbb{E}=\{(x_1, x_2, x_3) \in \mathbb{C}^3 :
1-x_1z_1-x_2z_2+x_3z_1z_2 \ne 0 \ \text{for all } \, z_1, z_2 \in \overline{\mathbb D}\}.
\]
We show that two other objects in $\mathbb C^4$ namely, the $\mu$-hexablock $\mathbb H_{\mu}$
and the normed hexablock $\mathbb H_N$ naturally arise in the $\mu_E$-unit ball and the norm
unit ball of $M_2(\mathbb C)$, respectively and pave the way to reach the domain
$\mathbb H$. A set of independent characterizations for the points in $\mathbb H_{\mu}, \mathbb H_N$
and $\mathbb H$ are obtained. Geometric and function theoretic aspects of $\mathbb H$ are
studied and its connections with the popular domains such as symmetrized bidisc
$\mathbb G_2$, tetrablock $\mathbb E$ and pentablock $\mathbb P$ are explored.
\end{abstract}	
	
	\maketitle

	\section{Introduction}\label{ch_intro} 
	
	\noindent Throughout the paper, the sets of complex numbers, real numbers, rational numbers, integers and natural numbers are denoted by $\mathbb{C}, \mathbb{R}, \mathbb{Q}, \mathbb{Z}$ and $\mathbb{N}$, respectively. The symbols $\mathbb{D}, \mathbb{T}$ stand for the open unit disc and the unit circle in the complex plane with their centres at the origin, respectively. We write $M_n(\mathbb{C})$ for the space of all $n \times n$ complex matrices. The adjoint, trace, determinant, spectral radius and operator norm of a matrix $A \in M_n(\mathbb{C})$ are denoted by $A^*, \text{tr}(A)$, $\det(A)$, $r(A)$ and $\|A\|$, respectively. For a set $\Omega$ in a topological space, we use $int(\Omega)$ and $\partial \Omega$ to refer to its interior and topological boundary, respectively.
	
	\smallskip
	
	In the theory of control engineering (e.g. see \cite{DoyleII, Francis}), the structured singular value of an $n\times n$ matrix $A$ is a cost function that generalizes the usual operator norm of $A$ and encodes structural information about the perturbations of $A$ that are being studied. The structured singular value is always computed with respect to a specified linear subspace of $ M_n(\mathbb{C})$. The linear subspace is usually referred to as the structure. Given a linear subspace $E$ of $M_n(\mathbb{C})$, the function
	\begin{equation} \label{eqn:NEW-01}
		\mu_E(A)=\frac{1}{\inf\{\|X\| \ : \ X \in E, \ \text{det}(I-AX)=0 \}}   \quad (A \in M_n(\mathbb{C})),
	\end{equation} 
	is called the \textit{structured singular value}. Here, $\|.\|$ denotes the operator norm of a matrix relative to the Euclidean norm on $\C^n$. If ${\rm det}\, (I-AX) \ne 0$ for all $X \in E$, then define $\mu_E(A)=0$. The notion of structured singular value of a matrix was introduced by Doyle \cite{Doyle} and was motivated by the problem of efficient stabilization of systems subjected to structured
	uncertainties. Two special cases of the structured singular value are the operator norm $\|.\|$ and the spectral radius $r$. If $E=M_n(\mathbb{C})$, then $\mu_E(A)=\|A\|$, and if $E$ is the space
	of scalar matrices, (i.e., scalar-times identity matrix), then $\mu_E(A)=r(A)$. These two instances of the structured singular value are extreme in the sense that $r(A) \leq \mu_E(A) \leq \|A\|$ for any linear subspace $E$ of $M_n(\mathbb{C})$ containing the identity matrix $I$. An interested reader is referred to \cite{Doyle} for control theoretic motivations behind $\mu_E$ and to \cite{Francis} for further details on $H^{\infty}$ control theory. For a linear subspace $E \subseteq M_n(\C)$, the aim of $\mu$-synthesis is to find a necessary and sufficient condition for the following interpolation problem: 
	
	\smallskip 
	
	\noindent ``Given distinct points $\lambda_1, \dotsc, \lambda_m$ in $\D$ and $n \times n$ complex matrices $A_1, \dotsc, A_m \in M_n(\C)$, does there exist an analytic $n \times n$ matrix-valued function $F$ on $\D$ such that
	\begin{align}\label{eqn_101}
		F(\lambda_j)=A_j \quad (j=1, \dotsc, m)
	\end{align}
	and 
	\begin{align}\label{eqn_102}
		\mu_E(F(\lambda)) \leq 1 \quad (\lambda \in \D)?\, "
	\end{align}
	Such an interpolation problem is referred to as the \textit{$\mu$-synthesis problem}. When $E=M_n(\C)$, the function $\mu_E$ coincides with the operator norm and the $\mu$-synthesis problem reduces to the classical \textit{Nevanlinna-Pick interpolation problem} in $M_n(\C)$, that is obtained by replacing the condition \eqref{eqn_102} by $\|F(\lambda)\| \leq 1$ for $\lambda \in \D$. This classical Nevanlinna-Pick interpolation problem has a complete solution (e.g., see Chapter X in \cite{Foias_Frazho} and Chapter XVIII in \cite{Ball}), where the existence of an interpolating function $F$ is equivalent to the semi-positive definiteness of the matrix given by 
	\begin{equation}\label{eqn_103}
		\begin{bmatrix}
			(1-\overline{\lambda}_i\lambda_j)^{-1}(I-A_i^*A_j)
		\end{bmatrix}_{i, j=1}^m.
	\end{equation}
The matrix in \eqref{eqn_103} is known as the \textit{Nevanlinna-Pick matrix} associated with the 
interpolation data. Furthermore, if $E\,=\,\{\alpha I\,\big\vert\,\, \alpha \in \C\}\subseteq M_n(\C)$, then $\mu_E$ is equal to 
the spectral radius and the $\mu$-synthesis problem becomes the \textit{spectral Nevanlinna-Pick} problem in 
which the condition \eqref{eqn_102} reduces to $r(F(\lambda)) \leq 1$ for $\lambda \in \D$.
	
	\smallskip
	
	We now show the readers by an example how the $\mu$-synthesis problem with respect to a linear subspace $E$ of $M_n(\C)$ gives rise to a domain. Suppose $E$ is the linear subspace of $M_2(\C)$ consisting of diagonal matrices, i.e., $E=\left\{\begin{pmatrix}
		z_1 & 0 \\
		0 & z_2
	\end{pmatrix}\,\,\big\vert\,\, z_1,\, z_2 \in \mathbb{C} \right\}$. It was shown in \cite{Abouhajar} that the $\mu$-synthesis with respect to this subspace induces the domain \textit{tetrablock}, which is defined by
	\begin{equation}\label{eqn:NEW-02}
		\mathbb{E}=\{(x_1, x_2, x_3) \in \mathbb{C}^3\,\big\vert\,\,
		1-x_1z_1-x_2z_2+x_3z_1z_2 \ne 0 \ \text{for all } \, z_1, z_2 \in \overline{\mathbb D}\}.
	\end{equation}
	Indeed, for $X=\begin{pmatrix}
		z_1 & 0\\
		0& z_2
	\end{pmatrix} \in E$ and $A=(a_{ij})_{i,j=1}^2 \in M_2(\mathbb{C})$, we have that
	$\|X\|=\max\{|z_1|, |z_2|\}$ and $\det(I-AX)=1-a_{11}z_1-a_{22}z_2+\det(A)z_1z_2$ and it follows from the definition of $\mu_E$ (see Equation-(\ref{eqn:NEW-01})) that
	\begin{align*}
		\mu_E(A)<1
		& \iff \|X\|>1 \ \text{for all} \ X \in E \ \text{with} \ \det(I-AX)=0\\
		& \iff \det(I-AX)\ne 0 \ \text{for all} \ X \in E \ \text{with} \ \|X\| \leq 1\\
		& \iff 1-a_{11}z_1-a_{22}z_2+\det(A)z_1z_2\ne 0 \ \text{for all} \ z_1, z_2 \in \overline{\mathbb{D}}\\
		& \iff (a_{11}, a_{22}, \det(A)) \in \mathbb{E}.
	\end{align*}
	Thus, $A=(a_{ij})_{i,j=1}^2$ is in the $\mu_E$-unit ball if and only if the point $(a_{11}, a_{22}, \det(A))$ belongs to the tetrablock $\mathbb E$. It turns out that an interpolation problem between $\D$ and the $\mu_E$-unit ball is equivalent with a similar interpolation problem between $\D$ and the tetrablock $\mathbb E$,  e.g., see the last Section of \cite{Abouhajar}. Evidently, it is easier to deal with a bounded domain like $\mathbb E \subset \C^3$ than the norm-unbounded object, the $\mu_E$ unit ball in $\C^4$. This is in fact the main reason behind studying such domains arising in the $\mu$-synthesis problem. The tetrablock has been studied extensively in recent past, see \cite{AbouhajarI, Abouhajar, Tirtha, Kosinski, zwonek1, Zwonek} and the references therein.
	
	\smallskip
	
	In a similar fashion, the subspace $E=\{\alpha I: \alpha \in \C\}\subseteq M_n(\C)$, for which $\mu_E$ is spectral radius, induces the \textit{symmetrized polydisc} (e.g. see \cite{Costara2005_II, EdigarianII}), a family of domains $\{\mathbb G_n: n\geq 2 \}$ defined by
	\[
	\mathbb{G}_n:=\left\{ \left(\sum_{1 \leq i \leq n}z_i, \sum_{1 \leq i<j \leq n}z_iz_j, \dotsc, \prod_{i=1}^{n} z_i\right) : z_1, \dotsc, z_n \in \D \right\} \subset \C^n.
	\]
	The relationship between $\mathbb{G}_n$ and the spectral unit ball $\mathbb B^{sp}_n$ of $M_n(\C)$ can be explained in the following way: a matrix $A \in M_n(\C)$ is in the spectral unit ball, i.e., $r(A)<1$ if and only if $\pi_n(\lambda_1, \dotsc, \lambda_n) \in \mathbb{G}_n$ (see \cite{Costara2005_II}), where $\lambda_1, \dotsc, \lambda_n$ belong to the spectrum $\sigma(A)$ of $A$ and $\pi_n$ is the symmetrization map on $\C^n$ defined by
	\[
	\pi_n(z_1, \dotsc, z_n):=\left(\sum_{1 \leq i \leq n}z_i, \sum_{1 \leq i<j \leq n}z_iz_j, \dotsc, \prod_{i=1}^{n} z_i\right).
	\]
	Thus, the map 
	\begin{align}\label{eqn_Pi_n}
		\Pi_n: M_n(\C) \to \C^n \quad \text{defined by} \quad  \Pi_n(A)=\pi_n(\sigma(A))
	\end{align}
	shows that $\mathbb{G}_n=\Pi_n(\mathbb{B}_n^{sp})$, see \cite{Costara2004, Costara2005_II}) for further details. The symmetrized polydisc has attracted considerable attentions in past two decades from perspectives of function theory, e.g., \cite{zwonek1, PZ, zwonek3, Ogle}, complex geometry \cite{AglerIII, bharali, EdigarianII, NN, Nikolov, zwonek4} and operator theory \cite{Agler2003, AglerYoung, Tirtha_Pal, tirtha-sourav1}.
	
	\smallskip
	
	Let us go back to the context of $\mu$-synthesis problem based on the linear subspaces of $M_2(\C)$ and the corresponding induced domains. To this end, we denote the norm unit ball and the $\mu_E$-unit ball with respect to a linear subspace $E$ of $M_2(\C)$ by
	\[
	\mathbb{B}_{\|.\|}=\{A \in M_2(\C) : \|A\| <1  \}\quad \text{and} \quad \mathbb{B}_{\mu_E}=\{A \in M_2(\mathbb{C}) \ : \ \mu_E(A) <1\}
	\] 
	respectively. As a special case of $\mu$-synthesis in $M_2(\C)$, we have already mentioned about the tetrablock, which arises when the linear subspace $E$ is chosen to be the space of $2 \times 2$ diagonal matrices. For this particular case, if we denote the cost function $\mu_E$ by $\mu_{\text{tetra}}$ and the corresponding $\mu_E$-unit ball $\mathbb{B}_{\mu_E}$ by $\mathbb{B}_{\mu_{\text{tetra}}}$, then it was proved in \cite{Abouhajar} that
	\begin{equation*}
		\begin{split}
			\E=\left\{(a_{11}, a_{22}, \det(A)) : A=(a_{ij})_{i, j=1}^2 \in \mathbb{B}_{\mu_{\text{tetra}}}\right\}
			=\left\{(a_{11}, a_{22}, \det(A)) : A=(a_{ij})_{i, j=1}^2 \in \mathbb{B}_{\|.\|}\right\}.
		\end{split}
	\end{equation*}
	Consequently, the images of $\mathbb{B}_{\mu_{\text{tetra}}}$ and $\mathbb{B}_{\|.\|}$ under the map 
	\begin{align}\label{eqn_pi_E}
		\pi_{\E}: M_2(\C) \to\C^3, \quad  A=(a_{ij})_{i, j=1}^2 \mapsto (a_{11}, a_{22}, \det(A))
	\end{align}
	coincide. However, the symmetrized bidisc $\mathbb G_2$ (which is associated with the $\mu_E$-unit ball with $E=\{\alpha I\,:\, \alpha \in \C  \} \subset M_2(\C)$) was the first born of $\mu$-synthesis in $M_2(\C)$ and appeared much earlier than the tetrablock, e.g., see \cite{AglerYoung}. Moreover, $\mathbb G_2$ is a member of the family of symmetrized polydisc $\{\mathbb G_n \, : \, n \geq 2 \}$ and like tetrablock the image of the $\mu_E$-unit ball and the norm unit ball under the $\Pi_2$ map coincides here too, i.e.,
	\[
	\Gg=\{(\text{tr}(A), \det(A)) : A \in M_2(\C), r(A)<1 \}=\{(\text{tr}(A), \det(A)) : A \in \mathbb{B}_{\|.\|}\}.
	\]
	A next step in this direction was taken by Agler et al. \cite{AglerIV}, where they considered the $\mu$-synthesis problem associated with the subspace
	\[
	E=\left\{ \begin{pmatrix}
		z & w \\
		0 & z
		
	\end{pmatrix} \,:\, z, w \in \C \right\} \subset M_2(\C)
	\]
	and this induced the domain pentablock $\Pe$. Denoting the map $\mu_E$ in this case by $\mu_{\text{penta}}$ and
	the corresponding $\mu_E$-unit ball by $\mathbb{B}_{\mu_{\text{penta}}}$, we have from \cite{AglerIV} that
	\[
	\Pe = \{(a_{21}, \text{tr}(A), \det(A)) : A=(a_{ij})_{i, j=1}^2 \in \mathbb{B}_{\mu_{\text{penta}}}\}= \left\{(a_{21}, \text{tr}(A), \det(A)) : A \in \mathbb{B}_{\|.\|} \right\}.
	\] 
	Thus, in this case also $\mathbb{B}_{\mu_E}$ and $\mathbb{B}_{\|.\|}$ have same images under the map
	\begin{align}\label{eqn_pi_Pe}
		\pi_{\Pe}: M_2(\C) \to \C^3, \quad A=(a_{ij})_{i, j=1}^2 \mapsto (a_{21}, \text{tr}(A), \det(A)).
	\end{align}
	An interested reader is referred to \cite{AglerIV, KosinskiII, Guicong} for a comprehensive reading on complex geometric aspects of the pentablock.
	
	\smallskip
	
	In this article, we move one step ahead and consider the linear subspace $E$ consisting of upper-triangular matrices, i.e.,
	\[
	E=\bigg\{\begin{pmatrix}
		z_1 & w\\
		0 & z_2
	\end{pmatrix} \ : \ z_1, z_2, w \in \mathbb{C} \bigg\} \subset M_2(\C).
	\]
	Denoting by $\mu_{\text{hexa}}$ the respective cost function $\mu_E$ and considering the map 
	\begin{align}\label{eqn_pi}
		\pi: M_2(\C) \to \C^4, \quad A=(a_{ij})_{i, j=1}^2 \mapsto (a_{21}, a_{11}, a_{22}, \text{det}(A)),
	\end{align}
	we obtain the following two sets in $\C^4$:
	\begin{equation*}
		\HB:=\{\pi(A)  :  A \in M_2(\C), \  \mu_{\text{hexa}}(A)<1 \}, \quad  \quad \QN:=\{\pi(A) : A \in M_2(\C), \|A\|<1\},
	\end{equation*}
	which we refer to as the \textit{$\mu$-hexablock} and the \textit{normed hexablock}, respectively. Like the relation between $\mu_{tetra}$-unit ball and tetrablock (discussed above), one can easily show that $\mu_{hexa}(A)<1$ if and only if $(a_{21},a_{11},a_{22}, \det(A)) \in \Q_{\mu}$, see Theorem \ref{thm2.4} for a proof to this. Evidently, this particular choice of $E$ generalizes 
	all other previous linear subspaces of $M_2(\C)$ (that correspond to $\Gg, \mathbb{E}$ and $\mathbb{P}$) and we have the following natural inequality:
	\begin{align}\label{eqn_mu_rel}
		r(A) \leq \mu_{\text{tetra}}(A), \mu_{\text{penta}}(A) \leq \mu_{\text{hexa}}(A) \leq \|A\|
	\end{align}
	for every $A \in M_2(\C)$. Also, the choice of the map $\pi$ is not arbitrary. It originates from the map $\pi_\E$ associated with $\E$ in 
	a similar manner that the map $\pi_{\Pe}$ associated with $\Pe$ arises from the map $\Pi_2$ associated with $\Gg$. The first notable facts about $\HB$ and $\QN$ are that neither they are domains in $\C^4$ nor they are equal (unlike the previous cases of $\Gg, \E$ and $\Pe$). Indeed, it is proved in Section \ref{norm_hexa} that $\HB \setminus \QN$ contains a non-empty open set in $\C^4$. However, both $\HB$ and $\QN$ are proved to be connected sets. In Sections \ref{mu_hexa} and \ref{norm_hexa}, we find some geometric properties of $\HB$ and $\QN$, 
	respectively and establish several characterizations of these two sets. In this direction, we mention another deviation of $\HB$ 
	and $\QN$ from the preceding three domains $\Gg, \E$ and $\Pe$. Note that if $(a, s, p) \in \Pe$, then $(s, p) \in \Gg$. Moreover, if 
	$(s, p) \in \Gg$, then there exist $z_1, z_2 \in \D$ such that $(s, p)=(z_1+z_2, z_1z_2)$. For a $2 \times 2$ 
	diagonal matrix $A$ with $z_1$ and $z_2$ on its diagonal, we have that $\|A\|<1$ and so, $\pi_{\Pe}(A)=(0, s, p) 
	\in \Pe$. Consequently, $\{0\} \times \Gg \subset \Pe$. We prove in Sections \ref{mu_hexa} and \ref{norm_hexa} 
	that $(x_1, x_2, x_3) \in \E$ whenever $(a, x_1, x_2, x_3)$ is in either $\HB$ or $\QN$. However, we also show that 
	$\{0\} \times \E$ is not contained in either of the sets $\HB$ or $\QN$. These deviations put forth a challenge in 
	defining an appropriate domain in $\C^4$ that contains $\{0\} \times \E$ providing an embedding $(x_1, x_2, x_3) 
	\mapsto (0, x_1, x_2, x_3)$ of $\E$ into this domain. To overcome this, we recall from \cite{AglerIV} a 
	characterization of the pentablock given by
	\begin{align}\label{eqn_109}
		\Pe=\left\{(a, s, p) \in \C \times \Gg : \sup_{z \in \D}|\psi_z(a, s, p)|<1\right\},   \quad \text{where} \quad 
		\psi_z(a, s, p)=\frac{a(1-|z|^2)}{1-sz+pz^2}.
	\end{align}
	The class of fractional linear maps $\{\psi_z: z \in \D\}$ play an important role in the theory of pentablock. In a similar manner, we consider in Section \ref{psi_map} a family of fractional linear maps defined by
	\[
	\Pz(a, x_1, x_2, x_3):=\frac{a\sqrt{(1-|z_1|^2)(1-|z_2|^2)}}{1-x_1z_1-x_2z_2+x_3z_1z_2},
	\]
	where $z_1, z_2 \in \overline{\D}$ and $(x_1, x_2, x_3) \in \E$. We discuss these maps and the motivation behind 
	them in Section \ref{psi_map}. Next, similar to the description of the pentablock in \eqref{eqn_109}, we define 
	\[
	\Q:=\left\{(a, x_1, x_2, x_3) \in \C \times \E: \sup_{z_1, z_2 \in \D}|\Pz(a, x_1, x_2, x_3)|<1\right\}. 
	\]
	The set $\Q$ is referred to as the \textit{hexablock}. This name arises from the fact that $\Q \cap \mathbb{R}^4$ is a bounded set, whose topological boundary consists of four extreme sets and two hypersurfaces. Specifically, the four extreme sets in the topological boundary of $\Q \cap \mathbb{R}^4$ correspond to the four faces of the convex set $\E \cap \mathbb{R}^3$, while the two hypersurfaces are induced by the two curved surfaces in the topological boundary of the real pentablock  $\Pe \cap \mathbb{R}^3$. In Section \ref{hexa}, we study the geometric properties of $\Q$, and prove that it is a domain in $\C^4$, unlike $\HB$ and $\QN$. Also, we have a natural embedding $(x_1, x_2, x_3) \mapsto (0, x_1, x_2, x_3)$ from $\E$ to $\Q$. By an \textit{embedding}, we mean an injective continuous map that is a homeomorphism onto its image. Consequently, one can employ the theory of $\E$ to explore the properties of $\Q$, in much the same way the theory of $\Pe$ is built upon the theory of $\Gg$.

	\smallskip 
	
	Although the sets $\HB$ and $\QN$ are not open in $\C^4$, one can still obtain the open sets out of them by considering the interiors of $\HB, \QN, \CHB$ and $\CQN$, that are denoted by $int(\HB), int(\QN), int(\CHB)$ and $int(\CQN)$ respectively. In fact, we prove in Section \ref{hexa} that each of these four sets is a domain in $\C^4$, and that $int(\CQN)=int(\QN) \subsetneq int(\HB) \subsetneq int(\CHB)$. However, we shall see in Sections \ref{mu_hexa}, \ref{norm_hexa} and \ref{hexa} that $\{0\} \times \E$ is contained in $int(\CHB)$, but not in $int(\HB), int(\QN)$, or $int(\CQN)$. As mentioned earlier, our aim is to obtain a domain in $\C^4$ that is closely related to $\HB$ and $\QN$, with $\{0\} \times \E$ being a part of it. For this reason, in addition to the hexablock $\Q$, we also focus on the domain $int(\CHB)$. Interestingly, it turns out that
	\[
	\Q= int(\CHB)=\HB \cup(\{0\} \times \E)= int(\HB) \cup (\{0\} \times \E), 
	\]
	which we prove in Section \ref{hexa}. This gives an alternative description of $\Q$ in terms of $\HB$. As discussed above, the term `hexablock' is justified by the fact that $\Q \cap \mathbb{R}^4$ contains six different blocks. We establish this geometric phenomenon in Section \ref{realslice} by showing that the topological boundary of $\Q \cap \mathbb{R}^4$ decomposes into precisely four extreme sets and two hypersurfaces.

\section{A class of linear fractional maps}\label{psi_map}

\noindent Attempts to solve different cases of $\mu$-synthesis problem led to the study of domains such as symmetrized bidisc, tetrablock and pentablock, as mentioned in the previous section. These domains exhibit many interesting properties from both complex function theoretic and operator theoretic points of view. A useful characterization for each of these domains is obtained via certain classes of linear fractional maps. For the symmetrized bidisc $\Gg$, the purpose is served by the following family of rational maps: 
\begin{align}\label{eqn_301}
	\Phi_z(s, p)=\frac{2zp-s}{2-zs}, \quad z \in \overline{\D}.
\end{align}
By Theorem 2.1 in \cite{AglerIII}, the class $\Phi_z$ determines a point in the symmetrized bidisc in the sense that $(s, p) \in \Gg$ if and only if $|s|<2$ and $\displaystyle \sup_{z \in \DC}|\Phi_z(s, p)|<1$. The tetrablock $\mathbb E$ has a similar characterization via the family of maps given by
\begin{align}\label{eqn_302}
	\Psi(z, x_1, x_2, x_3)  =\frac{x_3 z-x_1}{x_2 z-1}, \quad z \in \DC,
\end{align}
where $\Psi(., x_1, x_2, x_3)$ is the constant function equal to $x_1$ when $x_1x_2= x_3$. Also, $\Psi(z, x_1, x_2, x_3)$ is often denoted simply by $\Psi(z, x)$. An interested reader may refer to \cite{Abouhajar} for further details. Abouhajar et al. \cite{Abouhajar} proved that $(x_1, x_2, x_3) \in \E$ if and only if  $\displaystyle \sup_{z \in \mathbb{D}}|\Psi(z, x_1, x_2, x_3)|<1$ and $|x_2|<1$. It is not difficult to see that $\Phi_z(s, p)=-\Psi(z, s\slash 2, s\slash 2, p)$ and so, $|\Phi_z(s, p)|=|\Psi(z, s\slash 2, s\slash 2, p)|$. Consequently, the map $\Psi(z, .)$ as in \eqref{eqn_302} is a generalization of the map $\Phi_z$ as in \eqref{eqn_301}. In the same spirit, the points in the pentablock $\Pe$ can be described via the following family of linear fractional maps: 
\begin{align}\label{eqn_303}
	\psi_{z}\left(a, s, p\right)=\frac{a(1-|z|^2)}{1-sz+pz^2}, \quad z \in \DC.
\end{align}
Indeed, it follows from Theorem 5.2 in \cite{AglerIV} that  a point $(a, s, p) \in \Pe$ if and only if $(s, p) \in \Gg$ and $\underset{z \in \D}{\sup}|\psi_{z}(a, s, p)|<1$. Here we introduce another family of linear fractional maps $\Pz$ defined by
\begin{align}\label{eqn_304}
	\Pz(a, x_1, x_2, x_3):=\frac{a\sqrt{(1-|z_1|^2)(1-|z_2|^2)}}{1-x_1z_1-x_2z_2+x_3z_1z_2}, \quad z_1, z_2 \in \overline{\D}.
\end{align}
If $(x_1, x_2, x_3) \in \E$, it follows from the definition of $\E$ that $1-x_1z_1-x_2z_2+x_3z_1z_2 \ne 0$ for all $z_1, z_2 \in \overline{\D}$. Thus, $\Pz$ is a well-defined continuous function when $(a, x_1, x_2, x_3)  \in \C \times \E$. The map $\Pz$ generalizes the map $\psi_{z}$
as in \eqref{eqn_303} in the sense that
\[
\psi_{z, z}(a, x_1, x_2, x_3)=\frac{a\sqrt{(1-|z|^2)(1-|z|^2)}}{1-x_1z-x_2z+x_3z^2}=\frac{a(1-|z|^2)}{1-(x_1+x_2)z+x_3z^2}=\psi_{z}(a, x_1+x_2, x_3)
\]
for $z$ in $\DC$. Consequently, $\psi_z(a, s, p)=\psi_{z, z}(a, s\slash 2, s\slash 2, p)$ showing that the choice of the map $\Pz$ is not arbitrary. Indeed, the above discussion shows that the maps $\Pz$ are generalization of the maps $\psi_z$, much like how the maps $\Psi(z, .)$ generalize $\Phi_z$. In this section, we study the class of fractional linear maps 
$
\{\Pz: z_1, z_2 \in \DC\}
$
and identify the points in $\DC^2$ at which $\displaystyle \sup_{z_1, z_2 \in \DC}|\Pz(a, x_1, x_2, x_3)|$ is attained, where $(a, x_1, x_2, x_3) \in \C \times \E$. For this purpose, it suffices to find the supremum of the function given by
\begin{align}\label{eqn_kappa}
	\kappa(., x): \DC^2 \to \C, \quad  \kappa(z_1, z_2, x)=\frac{\sqrt{(1-|z_1|^2)(1-|z_2|^2)}}{1-x_1z_1-x_2z_2+x_3z_1z_2},
\end{align}
where $x=(x_1, x_2, x_3) \in \E$. Clearly, $\displaystyle \sup_{z_1, z_2 \in \DC}|\kappa(z_1, z_2, x)|$ exists since $\kappa(., x)$ is a continuous function in $(z_1, z_2)$ over the compact set $\overline{\mathbb{D}}^2$. Also, $\kappa(z_1, z_2, x)=0$ for $(z_1, z_2) \in \partial\D^2$ and so, $\sup_{z_1, z_2 \in \DC}|\kappa(z_1, z_2, x)|$ is attained at some point of $\D^2$. We explicitly determine these points below, which is the main result of this section. 

\begin{prop}\label{prop2.1}
	Let $x=(x_1, x_2, x_3) \in \mathbb{E}$. Then the supremum of $\bigg|\displaystyle\frac{\sqrt{(1-|z_1|^2)(1-|z_2|^2)}}{1-x_1z_1-x_2z_2+x_3z_1z_2}\bigg|$ over $(z_1, z_2)$ in $\DC^2$ is attained uniquely at a point in $\D^2$ given by
	\begin{align*}
	z_1(x)&=\frac{2\overline{\be}_1}{1+|\be_1|^2-|\be_2|^2+\sqrt{(1+|\be_1|^2-|\be_2|^2)^2-4|\be_1|^2}},\\	z_2(x)&=\frac{2\overline{\be}_2}{1+|\be_2|^2-|\be_1|^2+\sqrt{(1+|\be_2|^2-|\be_1|^2)^2-4|\be_2|^2}},
	\end{align*}
	where  $\displaystyle 
	\be_1=\frac{x_1-\overline{x}_2x_3}{1-|x_3|^2}$ and $\displaystyle \be_2=\frac{x_2-\overline{x}_1x_3}{1-|x_3|^2}$.
		
\end{prop}

\begin{proof} The proof follows a similar approach to that of Proposition 4.2 in \cite{AglerIV}.
The map given by 
$
h(z_1, z_2)=1-x_1z_1-x_2z_2+x_3z_1z_2=u(z_1, z_2)+iv(z_1, z_2)
$
is analytic on $\C^2$, where $u$ and $v$ are real and imaginary parts of $h$, respectively. Thus, $h(z_1, z_2)$ satisfies Cauchy- Riemann equations separately in both $z_1$ and $z_2$.  Since $(x_1, x_2, x_3) \in \E$, it follows from the definition of $\E$ that $|h(z_1, z_2)|>0$  for every $z_1, z_2 \in \mathbb{D}$. Let us define 
$\displaystyle 
g(z_1, z_2)=\frac{\sqrt{(1-|z_1|^2)(1-|z_2|^2)}}{|h(z_1, z_2)|}$ for $z_1, z_2 \in \D$. 
For any $z_1, z_2 \in \overline{\D}$, $h(z_1, z_2)$ is a non-zero complex number as $(x_1, x_2, x_3) \in \E$. Hence, $g$ extends to a continuous function on $\overline{\D}^2$ which takes value $0$ on the topological boundary of $\D^2$. Let $z_1=w_1+iy_1$ and $z_2=w_2+iy_2$. Then 
\begin{equation*}
	\begin{split}
		& \frac{\partial}{\partial w_1}|h(z_1, z_2)|= \frac{\partial}{\partial w_1}(u^2+v^2)^{\frac{1}{2}}=\frac{uu_{w_1}+vv_{w_1}}{|h(z_1, z_2)|},\\
		& \frac{\partial}{\partial y_1}|h(z_1, z_2)|= \frac{\partial}{\partial y_1}(u^2+v^2)^{\frac{1}{2}}=\frac{uu_{y_1}+vv_{y_1}}{|h(z_1, z_2)|}=\frac{vu_{w_1}-uv_{w_1}}{|h(z_1, z_2)|},\\
		& \frac{\partial}{\partial w_2}|h(z_1, z_2)|= \frac{\partial}{\partial w_2}(u^2+v^2)^{\frac{1}{2}}=\frac{uu_{w_2}+vv_{w_2}}{|h(z_1, z_2)|},\\
		& \frac{\partial}{\partial y_2}|h(z_1, z_2)|= \frac{\partial}{\partial y_2}(u^2+v^2)^{\frac{1}{2}}=\frac{uu_{y_2}+vv_{y_2}}{|h(z_1, z_2)|}=\frac{vu_{w_2}-uv_{w_2}}{|h(z_1, z_2)|}.
	\end{split}
\end{equation*}
Consequently, we have
\begin{equation*}
	\begin{split}
		&
		\frac{\partial}{\partial w_1}g(z_1, z_2)
		=\frac{\sqrt{1-w_2^2-y_2^2}}{|h(z_1, z_2)|^2}\ \ \left[\frac{-w_1|h(z_1, z_2)|}{\sqrt{1-w_1^2-y_1^2}}-\bigg(\frac{uu_{w_1}+vv_{w_1}}{|h(z_1, z_2)|}\bigg)\sqrt{(1-w_1^2-y_1^2)}\right], \\
		&
		\frac{\partial}{\partial y_1}g(z_1, z_2)=\frac{\sqrt{1-w_2^2-y_2^2}}{|h(z_1, z_2)|^2}\ \ \left[\frac{-y_1|h(z_1, z_2)|}{\sqrt{1-w_1^2-y_1^2}}-\bigg(\frac{vu_{w_1}-uv_{w_1}}{|h(z_1, z_2)|}\bigg)\sqrt{(1-w_1^2-y_1^2)}\right], \\
		&
		\frac{\partial}{\partial w_2}g(z_1, z_2)=\frac{\sqrt{1-w_1^2-y_1^2}}{|h(z_1, z_2)|^2}\ \ \left[\frac{-w_2|h(z_1, z_2)|}{\sqrt{1-w_2^2-y_2^2}}-\bigg(\frac{uu_{w_2}+vv_{w_2}}{|h(z_1, z_2)|}\bigg)\sqrt{(1-w_2^2-y_2^2)}\right],
		\\
		&					\frac{\partial}{\partial y_2}g(z_1, z_2)=\frac{\sqrt{1-w_1^2-y_1^2}}{|h(z_1, z_2)|^2}\ \ \left[\frac{-y_2|h(z_1, z_2)|}{\sqrt{1-w_2^2-y_2^2}}-\bigg(\frac{vu_{w_2}-uv_{w_2}}{|h(z_1, z_2)|}\bigg)\sqrt{(1-w_2^2-y_2^2)}\right].
	\end{split}
\end{equation*}
At critical points of $g$ in $\D^2$, we have that
\begin{equation*}
	\begin{split}
		&|h(z_1, z_2)|^2w_1=-(1-|z_1|^2)(uu_{w_1}+vv_{w_1}),\quad
		|h(z_1, z_2)|^2y_1=-(1-|z_1|^2)(vu_{w_1}-uv_{w_1}),\\
		& |h(z_1, z_2)|^2w_2=-(1-|z_2|^2)(uu_{w_2}+vv_{w_2}), \quad 
		|h(z_1, z_2)|^2y_2=-(1-|z_2|^2)(vu_{w_2}-uv_{w_2}).\\
	\end{split}
\end{equation*}
We solve these equations to obtain
\begin{equation*}
	\begin{split}
		& u_{w_1}=\frac{1}{|z_1|^2-1}(w_1u+y_1v), \qquad v_{w_1}=\frac{1}{|z_1|^2-1}(w_1v-y_1u),\\
		& u_{w_2}=\frac{1}{|z_2|^2-1}(w_2u+y_2v), \qquad v_{w_2}=\frac{1}{|z_2|^2-1}(w_2v-y_2u).\\	
	\end{split}
\end{equation*}
Hence, we have 
\begin{equation*}
	\begin{split}
		\frac{\partial}{\partial z_1}h(z_1, z_2)=u_{w_1}+iv_{w_1}
		=\frac{\overline{z}_1}{|z_1|^2-1}h(z_1, z_2), \quad 		\frac{\partial}{\partial z_2}h(z_1, z_2)=u_{w_2}+iv_{w_2}
		=\frac{\overline{z}_2}{|z_2|^2-1}h(z_1, z_2).
	\end{split}
\end{equation*}
Thus, the critical points of $g$ in $\D^2$ are the points $(z_1, z_2) \in \D^2$ such that 
\[
(-x_1+x_3z_2)(|z_1|^2-1)=\overline{z}_1(1-x_1z_1-x_2z_2+x_3z_1z_2) 
\]
and 
\[
(-x_2+x_3z_1)(|z_2|^2-1)=\overline{z}_2(1-x_1z_1-x_2z_2+x_3z_1z_2). 
\]
Consequently, it follows that the critical points $(z_1, z_2)$ of $g$ in $\D^2$ satisfy the following:
\begin{equation}\label{eqn3.1}
	\overline{z}_1=\frac{x_1-x_3z_2}{1-x_2z_2} \qquad \text{and} \qquad \overline{z}_2=\frac{x_2-x_3z_1}{1-x_1z_1}.
\end{equation}
Some laborious but routine calculations give that 
\begin{equation}\label{eqn3.2}
	\overline{z}_1^2(\overline{x}_1-x_2\overline{x}_3)-\overline{z}_1(1+|x_1|^2-|x_2|^2-|x_3|^2)+(x_1-\overline{x}_2x_3)=0.
\end{equation}
We compute the values of $z_1$ and $z_2$ depending on the choices of $\be_1$ and $\be_2$.

\medskip

\noindent Case $1$. Let $\be_1=0$ or $\be_2=0$. If $\be_1=0$, then $x_1-\overline{x}_2x_3=0$. By \eqref{eqn3.2}, $-\overline{z}_1(1+|x_1|^2-|x_2|^2-|x_3|^2)=0$. We have by part-$(4)$ of Theorem \ref{tetrablock} that $1+|x_1|^2-|x_2|^2-|x_3|^2>0$ since $(x_1, x_2, x_3) \in \E$. Thus, $z_1=0$ and so, $z_2=\overline{x}_2$ by \eqref{eqn3.1}. Similarly, one can show that $z_2=0$ and $z_1=\overline{x}_1$ when $\be_2=0$.

\medskip		 

\noindent Case $2$.	 Let $\be_1=\be_2=0$. Then $x_1=\overline{x}_2x_3$ and $x_2=\overline{x}_1x_3$. Thus, $x_1(1-|x_3|^2)=0$ and $x_2(1-|x_3|^2)=0$. Note that $|x_3|<1$ as $(x_1, x_2, x_3) \in \E$. Therefore, $x_1=x_2=0$. We have by \eqref{eqn3.1} and \eqref{eqn3.2} that $z_1=z_2=0$.

\medskip

\noindent Case $3$. Let $\be_1, \be_2 \ne 0$. Dividing both sides in  \eqref{eqn3.2} by $(1-|x_3|^2)$, we have that 
\[
\overline{z}_1^2\overline{\beta}_1-\overline{z}_1(1+|\be_1|^2-|\be_2|^2)+\be_1=0.
\]
The discriminant of the above quadratic equation (in the variable $z_1$) is given by 
\[
(1+|\be_1|^2-|\be_2|^2)^2-4|\be_1|^2 =(1-|\be_1|-|\be_2|)(1-|\be_1|+|\be_2|)(1+|\be_1|-|\be_2|)(1+|\be_1|+|\be_2|)
\]
which is strictly positive as $|\be_1|+|\be_2|<1$. Hence, the solutions to \eqref{eqn3.2} are 
\[
\overline{z}_1=\frac{(1+|\be_1|^2-|\be_2|^2) \pm \sqrt{(1+|\be_1|^2-|\be_2|^2)^2-4|\be_1|^2}}{2\overline{\be_1}}.
\]
Note that $(1+|\be_1|^2-|\be_2|^2)>\pm \sqrt{(1+|\be_1|^2-|\be_2|^2)^2-4|\be_1|^2}$. The above possible solutions of \eqref{eqn3.2} can be re-written as 
\[
\overline{z}_1=\frac{2\be_1}{(1+|\be_1|^2-|\be_2|^2)\pm \sqrt{(1+|\be_1|^2-|\be_2|^2)^2-4|\be_1|^2}}.
\]
Let if possible, 
\[
\overline{z}_1=\frac{2\be_1}{(1+|\be_1|^2-|\be_2|^2)- \sqrt{(1+|\be_1|^2-|\be_2|^2)^2-4|\be_1|^2}}.
\]
Since $|z_1|<1$, we have that
\begin{equation*}
	\begin{split}
		& \qquad \quad 2|\be_1|<(1+|\be_1|^2-|\be_2|^2) - \sqrt{(1+|\be_1|^2-|\be_2|^2)^2-4|\be_1|^2}\\
		& \implies \sqrt{(1+|\be_1|^2-|\be_2|^2)^2-4|\be_1|^2}<(1+|\be_1|^2-|\be_2|^2)-2|\be_1|\\
		& \implies 
		(1+|\be_1|^2-|\be_2|^2)^2-4|\be_1|^2<(1+|\be_1|^2-|\be_2|^2)^2+4|\be_1|^2-4|\be_1|(1+|\be_1|^2-|\be_2|^2)\\
		& \implies 	-|\be_1|^2< |\be_1|^2-|\be_1|(1+|\be_1|^2-|\be_2|^2)\\
		& \implies 	-|\be_1|< |\be_1|-(1+|\be_1|^2-|\be_2|^2) \qquad \qquad [\text{since} \ \be_1 \ne 0] \\
		& \implies 	1+|\be_1|^2-|\be_2|^2-2|\be_1|< 0 \\
		& \implies 	(1-|\be_1|+|\be_2|)(1-|\be_1|-|\be_2|)< 0, \\
	\end{split}
\end{equation*}
which is a contradiction as $|\be_1|+|\be_2|<1$.  Thus, we have 
\begin{equation}\label{eqn_z1}
	z_1=\frac{2\overline{\be}_1}{(1+|\be_1|^2-|\be_2|^2)+ \sqrt{(1+|\be_1|^2-|\be_2|^2)^2-4|\be_1|^2}}.
\end{equation}
Note that $|z_1|<1$, because
\begin{equation*}
	\begin{split}
		&(1+|\be_1|^2-|\be_2|^2)+ \sqrt{(1+|\be_1|^2-|\be_2|^2)^2-4|\be_1|^2}-2|\be_1|\\
		&=(1-|\be_1|)^2-|\be_2|^2)+ \sqrt{(1+|\be_1|^2-|\be_2|^2)^2-4|\be_1|^2}\\
		&=(1-|\be_1|-|\be_2|)(1-|\be_1|+|\be_2|)+ \sqrt{(1+|\be_1|^2-|\be_2|^2)^2-4|\be_1|^2}\\
		& >0 \qquad \qquad  [\text{as } |\be_1|+|\be_2|<1].
	\end{split}
\end{equation*}
Employing similar arguments and computations as above, we find that 
\begin{equation}\label{eqn_z2}
	z_2=\frac{2\overline{\be}_2}{(1+|\be_2|^2-|\be_1|^2)+\sqrt{(1+|\be_2|^2-|\be_1|^2)^2-4|\be_2|^2}}.
\end{equation}
In either case, we have only one choice for the critical points $(z_1, z_2)$ of $g$  in $\D^2$. Also, it is easy to see that the solutions for $z_1$ and $z_2$ in \eqref{eqn_z1} \& \eqref{eqn_z2} also include the solutions of $z_1$ and $z_2$ from Cases $1$ and $2$. Since $g$ is a continuous function over $\DC^2$ and $g|_{\partial \D^2}=0$, it follows that $g$ attains its maximum over $\D^2$. Since $(z_1, z_2)$ is a unique critical point of $g$ in $\D^2$, it must be a point of global maximum for $g$ in $\D^2$. The proof is now complete.
\end{proof} 

Using the same techniques as in the proof of Proposition \ref{prop2.1}, we obtain the following result.

\begin{cor}\label{cor2.2}
	Let $x=(x_1, x_2, x_3) \in \mathbb{E} \cap \mathbb{R}^3$. Then the supremum of $\displaystyle\frac{\sqrt{(1-|z_1|^2)(1-|z_2|^2)}}{1-x_1z_1-x_2z_2+x_3z_1z_2}$ over $(z_1, z_2)$ in $[-1, 1] \times [-1, 1]$ is attained uniquely at a point in $(-1, 1) \times (-1, 1)$ given by
	\begin{align*}
	z_1(x)&=\frac{2\be_1}{1+\be_1^2-\be_2^2+\sqrt{(1+\be_1^2-\be_2^2)^2-4\be_1^2}},\\
	z_2(x)&=\frac{2\be_2}{1+\be_2^2-\be_1^2+\sqrt{(1+\be_2^2-\be_1^2)^2-4\be_2^2}},
	\end{align*}
	where $\displaystyle \be_1=\frac{x_1-x_2x_3}{1-x_3^2}$ and $\displaystyle \be_2=\frac{x_2-x_1x_3}{1-x_3^2}$.
\end{cor}

It is clear from Proposition \ref{prop2.1} and Corollary \ref{cor2.2} that if $x=(x_1, x_2, x_3) \in \E \cap \mathbb{R}^3$, then the points $z_1(x), z_2(x)$ as in the above corollary belong to  the interval $(-1, 1)$ and so, we have the following result. 

\begin{cor}\label{cor_sup}
	Let $x=(x_1, x_2, x_3) \in \mathbb{E} \cap \mathbb{R}^3$. Then 
	\begin{align*}
		\sup_{(z_1, z_2)\in \DC^2}\left|\displaystyle\frac{\sqrt{(1-|z_1|^2)(1-|z_2|^2)}}{1-x_1z_1-x_2z_2+x_3z_1z_2}\right|
		&=	\sup_{z_1, z_2\in (-1, 1)}\left|\displaystyle\frac{\sqrt{(1-|z_1|^2)(1-|z_2|^2)}}{1-x_1z_1-x_2z_2+x_3z_1z_2}\right|\\
		&=\left|\displaystyle\frac{\sqrt{(1-|z_1(x)|^2)(1-|z_2(x)|^2)}}{1-x_1z_1(x)-x_2z_2(x)+x_3z_1(x)z_2(x)}\right|,
	\end{align*}
	where $z_1(x)$ and $z_2(x)$ are as in Corollary \ref{cor2.2}. 
\end{cor}

For $x=(x_1, x_2, x_3) \in \E$, let $(z_1(x), z_2(x)) \in \D^2$ be the unique point where $\kappa(., x)$ from \eqref{eqn_kappa} attains its maximum modulus. The existence of such a point is always guaranteed by Proposition \ref{prop2.1}. Furthermore, the function given by
\begin{align}\label{eqn_K_*}
	K_*: \E \to [1, \infty), \quad  K_*(x_1, x_2, x_3)=|\kappa(z_1(x), z_2(x), x)|
\end{align}
is a well-defined continuous function. The continuity of $K_*$ follows from the facts that for any $x \in \E$, the maps $z_1(x)$ and $z_2(x)$ are continuous functions of $x$ and $x \mapsto \kappa(z_1(x), z_2(x), x)$ is also continuous on $\E$. Also, $K_*(x_1, x_2, x_3) \geq |\kappa(0, 0, x)| =1$ and it follows from Proposition \ref{prop2.1} that
\[
\sup_{z_1, z_2 \in \DC}|\psi_{z_1, z_2}(a, x_1, x_2, x_3)|=\sup_{z_1, z_2 \in \D}|\psi_{z_1, z_2}(a, x_1, x_2, x_3)|= |a|K_*(x_1, x_2, x_3)
\]
for $(a, x_1, x_2, x_3) \in \C \times \E$. Putting everything together, we arrive at the following corollary.

\begin{cor}\label{cor2.3}
	Let $(x_1, x_2, x_3) \in \E$ and let $a \in \C$. Then the following are equivalent:
	\begin{enumerate}
		\item $\displaystyle	\sup_{z_1, z_2 \in \D}|\psi_{z_1, z_2}(a, x_1, x_2, x_3)|<1$;
		\item $|\psi_{z_1, z_2}(a, x_1, x_2, x_3)|<1$ for all $(z_1, z_2) \in \D^2$;
		\item $|a|K_*(x_1, x_2, x_3)<1$.
	\end{enumerate}
\end{cor}

Let us discuss some useful consequences of Proposition \ref{prop2.1} depending on $(x_1, x_2, x_3) \in \E$. We need to describe briefly a few characterizations of the symmetrized bidisc $\mathbb G_2$ (and its closure $\Gamma_2$) along with the same for the tetrablock $\E$ and its closure $\overline{\E}$. Let us recall from the `Introduction' the definitions of $\mathbb{G}_2, \E$ and their closures:
\begin{align*}
\Gg&=\{(z_1+z_2, z_1z_2) \in \C^2 : z_1, z_2 \in \D \}, \quad \Gamma_2=\{(z_1+z_2, z_1z_2) \in \C^2 : z_1, z_2 \in \overline{\D} \}; \\
\E&=\{(x_1, x_2, x_3) \in \C^3 \ : \ 1-x_1z_1-x_2z_2+x_3z_1z_2 \ne 0 \ \text{whenever} \ |z_1| \leq 1, |z_2| \leq 1 \}; \\
\overline{\E}&=\{(x_1, x_2, x_3) \in \C^3 \ : \ 1-x_1z_1-x_2z_2+x_3z_1z_2 \ne 0 \ \text{whenever} \ |z_1| < 1, |z_2| < 1 \}.
\end{align*}

\begin{thm}[\cite{AglerIII}, Theorem 2.1 \& \cite{Agler2003}, Theorem 1.1]\label{thmG_2}
	Let $(s,p)\in\C^2$. The following are equivalent.
	
	\begin{enumerate}
		\item $(s,p)\in \Gg$   $($or, $(s,p)\in \Gamma_2$$)$;
		
		\item $|s-\overline{s}p|<1-|p|^2$ 
		$($or, $|s|\le 2$ and
		$|s-\overline{s}p|\le 1-|p|^2$$)$;
		
		\item
		$
		2|s-\overline{s}p|+|s^2-4p|
		<4-|s|^2
		$ 
		$($or,
		$
		2|s-\overline{s}p|+|s^2-4p|
		\le 4-|s|^2$$)$;

		\item $|p|<1$ and there exists $\beta\in\D$ such that
		$s=\beta+\overline{\beta}p$ 
		$($or, $|p|\le 1$ and there exists
		$\beta\in\overline{\D}$ such that
		$s=\beta+\overline{\beta}p$$)$;
		
		\item there exists $A\in M_2(\C)$ with $\|A\|<1$
		such that
		$(s,p)=(\text{tr}(A),\det(A))$
		$($or, there exists $A\in M_2(\C)$ with
		$\|A\|\le 1$ such that
		$(s,p)=(\text{tr}(A),\det(A))$$)$.
	\end{enumerate}
\end{thm}

\begin{thm}[\cite{Abouhajar}, Theorem 2.2]\label{tetrablock}
	For $(x_1, x_2, x_3) \in \C^3$, the following are equivalent:
	\begin{enumerate}
		\item $(x_1, x_2, x_3) \in \E$ ;
		\vspace{0.1cm}
		\item $\sup_{z \in \mathbb{D}}|\Psi(z, x_1, x_2, x_3)|<1$ and if $x_1x_2=x_3$ then, in addition, $|x_2|<1;$
		
		\item $|x_1|^2+|x_2-\overline{x}_1x_3|+|x_1x_2-x_3|<1$;
		
		\item $1+|x_1|^2-|x_2|^2-|x_3|^2-2|x_1-\overline{x}_2x_3|>0$ and $|x_1|<1$;
		\vspace{0.1cm}
		
		\item $|x_1-\overline{x}_2x_3|+|x_2-\overline{x}_1x_3| < 1-|x_3|^2$;
		
		\item there is a $2 \times 2$ matrix $A=(a_{ij})$ such that $\|A\|<1$ and $x=(a_{11}, a_{22}, \det(A))$;
		\vspace{0.1cm}
		\item $|x_3|<1$ and there exist $\beta_1, \beta_2 \in \mathbb{C}$ such that $|\beta_1|+|\beta_2| <1$ and 
		\[
		x_1=\beta_1+\overline{\beta}_2x_3, \quad x_2=\beta_2+\overline{\beta}_1x_3;
		\]
		
		\item $1-|x_1|^2-|x_2|^2+|x_3|^2>2|x_1x_2-x_3|$ and if $x_1x_2=x_3$ then, in addition, $|x_1|+|x_2|<2$.
	\end{enumerate}	
\end{thm}

\begin{thm}[\cite{Abouhajar}, Theorem 2.4]\label{closedtetrablock}
	For $(x_1, x_2, x_3) \in \C^3$, the following are equivalent:
	\begin{enumerate}
		\item $(x_1, x_2, x_3) \in \overline{\E}$ ;
		\vspace{0.1cm} 
		\item $1-x_1z_1-x_2z_2+x_3z_1z_2 \ne 0$ for all $z_1, z_2 \in \D$;
		\vspace{0.1cm}
		\item $\sup_{z \in \mathbb{D}}|\Psi(z, x_1, x_2, x_3)| \leq 1$ and if $x_1x_2=x_3$ then, in addition, $|x_2|\leq 1;$
		\vspace{0.1cm}
		\item $1+|x_1|^2-|x_2|^2-|x_3|^2-2|x_1-\overline{x}_2x_3| \geq 0$ and $|x_1|\leq 1$;
		\vspace{0.1cm}
		\item there is a $2 \times 2$ matrix $A=(a_{ij})$ such that $\|A\| \leq 1$ and $x=(a_{11}, a_{22}, \det(A))$;
		\vspace{0.1cm}
		\item $|x_3| \leq 1$ and there exist $\beta_1, \beta_2 \in \mathbb{C}$ such that $|\beta_1|+|\beta_2| \leq 1$ and 
		\[
		x_1=\beta_1+\overline{\beta}_2x_3, \quad x_2=\beta_2+\overline{\beta}_1x_3.
		\]
	\end{enumerate}	
\end{thm}

We begin with the following corollary that follows directly from Proposition \ref{prop2.1}. 

\begin{cor}\label{cor_302}
	Let $(x_1, x_2, x_3) \in \mathbb{E}$ with $x_1=x_2$. Then $\bigg|\displaystyle\frac{\sqrt{(1-|z_1|^2)(1-|z_2|^2)}}{1-x_1z_1-x_2z_2+x_3z_1z_2}\bigg|$ attains its supremum over $\DC^2$ at a unique point $(z_1, z_2)$ given by
	\[
	z_1=z_2=\frac{2\overline{\be_1}}{1+\sqrt{1-4|\be_1|^2}}, \quad \text{where} \quad
	\be_1=\frac{x_1-\overline{x}_1x_3}{1-|x_3|^2}.
	\]	
\end{cor}

The next result, appearing as Proposition 4.2 in \cite{AglerIV}, follows simply from the above corollary.

\begin{cor}
	Let $(s, p) \in \Gg$. Then the supremum of $\displaystyle \frac{1-|z|^2}{|1-sz+pz^2|}$ over $z \in \DC$ is attained uniquely at a point in $\D$ given by
	\[
	\lambda=\frac{\overline{\be}}{1+\sqrt{1-|\be|^2}}, \quad \text{where} \quad \beta=\frac{s-\overline{s}p}{1-|p|^2}.
	\]
\end{cor}

\begin{proof} 
We have by part-$(2)$ of Theorem \ref{thmG_2} and part-$(4)$ of Theorem \ref{tetrablock} that $(s\slash 2, s\slash 2, p)\in \E$. By Corollary \ref{cor_302}, it follows that
\[
\bigg|\displaystyle\frac{\sqrt{(1-|z_1|^2)(1-|z_2|^2)}}{1-(s\slash 2)z_1-(s\slash 2)z_2+pz_1z_2}\bigg|
\]
attains its supremum over $\DC^2$ at $(\lm, \lm) \in\D^2$ with
$
\lm=\overline{\be}(1+\sqrt{1-|\be|^2})^{-1}$. Then
\begin{equation*}
	\begin{split}
		\frac{1-|\lm|^2}{|1-s\lm+p\lm^2|} \leq \sup_{z \in \DC} \frac{1-|z|^2}{|1-sz+pz^2|}
		\leq \sup_{z_1, z_2 \in \DC}\bigg|\frac{\sqrt{(1-|z_1|^2)(1-|z_2|^2)}}{1-(s\slash 2)z_1-(s\slash 2)z_2+pz_1z_2}\bigg|
		=\frac{1-|\lm|^2}{|1-s\lm+p\lm^2|}
	\end{split}
\end{equation*}
and so, the desired conclusion follows. 
\end{proof} 

We now consider the case when $x=(x_1, x_2, x_3) \in \E$ with $x_1x_2=x_3$ and find the supremum of $|\kappa(z_1, z_2, x)|$ over $(z_1, z_2)$ in $\DC^2$. Note that a point $(x_1, x_2, x_3) \in \C^3$ for which $x_1x_2=x_3$ is called a \textit{triangular point}.

\begin{lem}\label{lem3.4}
	Let $(a, x_1, x_2, x_3) \in \C \times \E$ and let $x_1x_2=x_3$. Then 
	\[
	\sup_{z_1, z_2 \in \D} \bigg|\frac{\sqrt{(1-|z_1|^2)(1-|z_2|^2)}}{1-x_1z_1-x_2z_2+x_3z_1z_2}\bigg|=\frac{1}{\sqrt{(1-|x_1|^2)(1-|x_2|^2)}}.
	\]
	In particular, $\sup_{z_1, z_2 \in \D}|\psi_{z_1, z_2}(a, x_1, x_2, x_3)|<1$
	if and only if $|a| < (1-|x_1|^2)^{1\slash 2}(1-|x_2|^2)^{1\slash 2}$.
\end{lem}

\begin{proof}
Since $x_3=x_1x_2$, we have that
\begin{equation*}
	\begin{split}
		\Pz(a, x_1, x_2, x_3)&=\frac{a\sqrt{(1-|z_1|^2)(1-|z_2|^2)}}{1-x_1z_1-x_2z_2+x_3z_1z_2}=a\left(\frac{\sqrt{1-|z_1|^2}}{1-x_1z_1}\right)\left(\frac{\sqrt{1-|z_2|^2}}{1-x_2z_2}\right)
	\end{split}
\end{equation*}
and so,
\begin{equation*}
	\begin{split}
		\sup_{z_1, z_2 \in \D} |\Pz(a, x_1, x_2, x_3)|=|a|\ \ \sup_{z_1 \in \D}\left|\frac{\sqrt{1-|z_1|^2}}{1-x_1z_1}\right| \ \ \sup_{z_2 \in \D}\left|\frac{\sqrt{1-|z_2|^2}}{1-x_2z_2}\right|.
	\end{split}
\end{equation*}
For $x,z \in \D$, we have that $|1-xz|^2-(1-|z|^2)(1-|x|^2)=|x-\overline{z}|^2 \geq 0$. Then
\begin{equation*}
	\begin{split}
		\sup_{z \in \D}\frac{\sqrt{1-|z|^2}}{|1-xz|} \leq \frac{1}{\sqrt{1-|x|^2}}
		\quad \text{and so,} \quad	\sup_{z_1, z_2 \in \D} |\Pz(a, x_1, x_2, x_3)| \leq \frac{|a|}{\sqrt{(1-|x_1|^2)(1-|x_2|^2)}}.
	\end{split}
\end{equation*}
The equality holds when $(z_1, z_2)=(\overline{x}_1, \overline{x}_2)$ and so, the desired conclusion follows.
\end{proof}

So far, we have studied the map $\kappa(z_1, z_2, x)$ from \eqref{eqn_kappa} for $x=(x_1, x_2, x_3) \in \E$ and $(z_1, z_2) \in \DC^2$. The underlying reason is that $1-x_1z_1-x_2z_2+x_3z_1z_2 \ne 0$ for such a choice of $x$ and $(z_1, z_2)$. However, for $x \in \Ebar$ the map $\kappa(z_1, z_2, x)$ need not be defined for $(z_1, z_2) \in \DC^2$. For example, if $x=(x_1, x_2, x_3)=(1, 0, 0) \in \Ebar$, then 
\begin{align}\label{eqn_k2}
	\kappa(z_1, z_2, x)=\frac{\sqrt{(1-|z_1|^2)(1-|z_2|^2)}}{1-x_1z_1-x_2z_2+x_3z_1z_2}=\frac{\sqrt{(1-|z_1|^2)(1-|z_2|^2)}}{1-z_1}
\end{align}
and so, $\kappa(1, 0, x)$ is not well-defined. Nevertheless, there is a natural way to study the map $\kappa(z_1, z_2, x)$ when $x \in \Ebar$. To do so, we must consider $(z_1, z_2)$ in $\D^2$, since part-$(2)$ of Theorem \ref{closedtetrablock} guarantees that $1-x_1z_1-x_2z_2+x_3z_1z_2 \ne 0$ for $(x_1, x_2, x_3) \in \Ebar$ and $(z_1, z_2) \in \D^2$. As a result, $\kappa(z_1, z_2, x)$ and hence $\Pz(a, x_1, x_2, x_3)$ is well-defined for $(a, x) \in \C \times \Ebar$ and $(z_1, z_2) \in \D^2$. In line with the results of this section, one may ask whether $\displaystyle \sup_{z_1, z_2 \in \D}|\kappa(z_1, z_2, x)|$ exists when $x \in \Ebar$. For $x=(1, 0, 0)$, it follows from \eqref{eqn_k2} that 
\[
\sup_{z_1, z_2 \in \D}|\kappa(z_1, z_2, x)| \geq \sup_{z_1 \in \D}|\kappa(z_1, 0, x)| \geq \sup_{z_1 \in \D}\frac{\sqrt{1-|z_1|^2}}{|1-z_1|} \geq \sup_{\alpha \in (-1, 1)}\frac{\sqrt{1-\alpha^2}}{1-\alpha} =\sup_{\alpha \in (-1, 1)}\frac{\sqrt{1+\alpha}}{\sqrt{1-\alpha}}
\]
and so,  $\displaystyle \sup_{z_1, z_2 \in \D}|\kappa(z_1, z_2, x)|$ does not exist. However, we show that  $\sup_{z_1, z_2 \in \D}|\kappa(z_1, z_2, x)|$ exists for a class of points $x$ belonging to the distinguished boundary $b\E$ of $\Ebar$. We recall from \cite{Abouhajar} a description of the distinguished boundary of the tetrablock.

\begin{thm}[\cite{Abouhajar}, Theorem 7.1]\label{thm6.1}
	For $x=(x_1, x_2, x_3) \in \C^3$, the following are equivalent:
	\begin{enumerate}
		\item $x_1=\overline{x}_2x_3, |x_3|=1$ and $|x_2|\leq 1$;
		\item $x \in b\E$;
		\item $x \in \overline{\E}$ and $|x_3|=1$. 
	\end{enumerate}	
\end{thm}  

Now, we compute $\displaystyle \sup_{z_1, z_2 \in \D}|\kappa(z_1, z_2, x)|$ when $x=(x_1, x_2, x_3) \in b\E$ with $|x_1|<1$.

\begin{thm}\label{thm_sup_bE}
	Let $x=(x_1, x_2, x_3) \in b\E$ and let $|x_1|< 1$. Then 
	\[
	\sup_{z_1, z_2 \in \D} \bigg|\frac{\sqrt{(1-|z_1|^2)(1-|z_2|^2)}}{1-x_1z_1-x_2z_2+x_3z_1z_2}\bigg|=\frac{1}{\sqrt{1-|x_1|^2}}.
	\]
\end{thm}

\begin{proof} Let $(z_1, z_2) \in \D^2$ and let $x_r=(rx_1, rx_2, r^2x_3)$ for $r \in (0,1)$. Since $x \in \Ebar$, we have by part-$(6)$ of Theorem \ref{tetrablock} and part-$(5)$ of Theorem \ref{closedtetrablock} that $x_r \in \E$. It follows from Proposition \ref{prop2.1} that 
\begin{align}\label{eqn_312}
	\bigg|\frac{\sqrt{(1-|z_1|^2)(1-|z_2|^2)}}{1-rx_1z_1-rx_2z_2+r^2x_3z_1z_2}\bigg| \leq \bigg|\frac{\sqrt{(1-|z_1(x_r)|^2)(1-|z_2(x_r)|^2)}}{1-rx_1z_1(x_r)-rx_2z_2(x_r)+r^2x_3z_1(x_r)z_2(x_r)}\bigg|,
\end{align} 
where $(z_1(x_r), z_2(x_r))$ is the unique point in $\D^2$ as in Proposition \ref{prop2.1} that corresponds to $x_r \in \E$. Let us define 
\[
\be_{1r}:=\frac{rx_1-r^3\overline{x}_2x_3}{1-r^4|x_3|^2} \quad \text{and} \quad  \beta_{2r}:=\frac{rx_2-r^3\overline{x}_1x_3}{1-r^4|x_3|^2}.
\]
It follows from Theorem \ref{thm6.1} that $x_1=\overline{x}_2x_3, x_2=\overline{x}_1x_3$ and $|x_3|=1$ as $x \in b\E$. Then
\begin{align}\label{eqn_313}
	\beta_{1r}=\frac{rx_1}{1+r^2} \quad \text{and} \quad \beta_{2r}=\frac{rx_2}{1+r^2}.
\end{align}
Since $|x_1|=|x_2|$, we have that $|\beta_{1r}|=|\beta_{2r}|$. Again by Proposition \ref{prop2.1}, we have that
\begin{align}\label{eqn_314}
	z_1(x_r)=\frac{2\overline{\be}_{1r}}{1+|\be_{1r}|^2-|\be_{2r}|^2+\sqrt{(1+|\be_{1r}|^2-|\be_{2r}|^2)^2-4|\be_{1r}|^2}}=\frac{2\overline{\be}_{1r}}{1+\sqrt{1-4|\be_{1r}|^2}}, \notag \\
	z_2(x_r)=\frac{2\overline{\be}_{2r}}{1+|\be_{2r}|^2-|\be_{1r}|^2+\sqrt{(1+|\be_{2r}|^2-|\be_{1r}|^2)^2-4|\be_{2r}|^2}}=\frac{2\overline{\be}_{2r}}{1+\sqrt{1-4|\be_{1r}|^2}}.
\end{align}
So, $|z_1(x_r)|=|z_2(x_r)|$, because $|\be_{1r}|=|\be_{2r}|$. A routine calculation gives that
\begin{align}\label{eqn_315}
	& \quad 1-rx_1z_1(x_r)-rx_2z_2(x_r)+r^2x_3z_1(x_r)z_2(x_r) \notag \\ 
	&=1-\frac{2rx_1\overline{\be}_{1r}}{1+\sqrt{1-4|\be_{1r}|^2}}-\frac{2rx_2\overline{\be}_{2r}}{1+\sqrt{1-4|\be_{1r}|^2}}+\frac{4r^2x_3\overline{\be}_{1r}\overline{\be}_{2r}}{(1+\sqrt{1-4|\be_{1r}|^2})^2} \notag \\
	&=1-\frac{4rx_1\overline{\be}_{1r}}{1+\sqrt{1-4|\be_{1r}|^2}}+\frac{4r^2x_3\overline{\be}_{1r}\overline{\be}_{2r}}{(1+\sqrt{1-4|\be_{1r}|^2})^2} \qquad \qquad [\text{as } x_1\overline{\be}_{1r}=x_2\overline{\be}_{2r} ] \notag \\
	&=1-\frac{4r^2|x_1|^2}{(1+r^2)(1+\sqrt{1-4|\be_{1r}|^2})}+\frac{4r^4|x_1|^2}{(1+r^2)^2(1+\sqrt{1-4|\be_{1r}|^2})^2},
\end{align}
where the last equality follows from \eqref{eqn_313} and the fact that $x_1=\overline{x}_2x_3$. It is clear from \eqref{eqn_313} that $\lim_{r \to 1}\be_{1r}=x_1\slash 2$ and $\lim_{r \to 1}\be_{2r}=x_2\slash 2$. Thus, we have from \eqref{eqn_314} that
\begin{align*}
	\lim_{r \to 1}(1-|z_j(x_r)|^2)
	=\lim_{r\to1}\left[1-\frac{4|\be_{1r}|^2}{(1+\sqrt{1-4|\be_{1r}|^2})^2}\right]
	=1-\frac{|x_1|^2}{(1+\sqrt{1-|x_1|^2})^2}=\frac{2\sqrt{1-|x_1|^2}}{1+\sqrt{1-|x_1|^2}}
\end{align*}
for $j=1, 2$ since $|z_1(x_r)|=|z_2(x_r)|$. Also, we have by \eqref{eqn_315} that 
\begin{align*}
	&\quad \lim_{r\to 1}(1-rx_1z_1(x_r)-rx_2z_2(x_r)+r^2x_3z_1(x_r)z_2(x_r))\\
	&=\lim_{r\to 1}\left[1-\frac{4r^2|x_1|^2}{(1+r^2)(1+\sqrt{1-4|\be_{1r}|^2})}+\frac{4r^4|x_1|^2}{(1+r^2)^2(1+\sqrt{1-4|\be_{1r}|^2})^2}\right] \quad [\text{by} \ \eqref{eqn_315}]\\
	&=1-\frac{2|x_1|^2}{1+\sqrt{1-|x_1|^2}}+\frac{|x_1|^2}{(1+\sqrt{1-|x_1|^2})^2}\\
	&=\frac{2(1-|x_1|^2)}{1+\sqrt{1-|x_1|^2}}.
\end{align*}
Taking limits as $r \to 1$ in \eqref{eqn_312}, we have 
\begin{align*}
	\bigg|\frac{\sqrt{(1-|z_1|^2)(1-|z_2|^2)}}{1-x_1z_1-x_2z_2+x_3z_1z_2}\bigg|
	& \leq \lim_{r \to 1} \bigg|\frac{\sqrt{(1-|z_1(x_r)|^2)(1-|z_2(x_r)|^2)}}{1-rx_1z_1(x_r)-rx_2z_2(x_r)+r^2x_3z_1(x_r)z_2(x_r)}\bigg|\\
	& = \lim_{r \to 1} \bigg|\frac{1-|z_1(x_r)|^2}{1-rx_1z_1(x_r)-rx_2z_2(x_r)+r^2x_3z_1(x_r)z_2(x_r)}\bigg|\\
	&=\frac{1}{\sqrt{1-|x_1|^2}}.
\end{align*} 
Consequently, it follows that
\[
\sup_{z_1, z_2 \in \D} \bigg|\frac{\sqrt{(1-|z_1|^2)(1-|z_2|^2)}}{1-x_1z_1-x_2z_2+x_3z_1z_2}\bigg| \leq \frac{1}{\sqrt{1-|x_1|^2}},
\]
and the supremum is attained when $z_1=\overline{x}_1$ and $z_2=0$. The proof is now complete. 
\end{proof} 

\section{The $\mu$-hexablock}\label{mu_hexa}

\noindent The motivation for this work is to identify a new instance of the inherently challenging $\mu$-synthesis problem that is suitable for a more tractable analysis. This section serves as an initial step in this endeavor, wherein we introduce the $\mu$-hexablock, a subset of $\mathbb{C}^4$ closely related to the $\mu$-synthesis problem in $M_2(\C)$. Recall from Section \ref{ch_intro} that $\mu$ denotes a cost function on $M_n(\C)$ given by
\[
\mu_E(A)=\frac{1}{\inf\{\|X\| \ : \ X \in E, \ \text{det}(I-AX)=0 \}},
\] 
with the understanding that $\mu_E(A)=0$ if $\det(I-AX) \ne 0$ for all $X \in E$. Here $E$ is a linear subspace of $M_n(\C)$. Previous attempts to study a few instances of $\mu$-synthesis problem led to the domains such as symmetrized bidisc $\Gg$, tetrablock $\E$ and pentablock $\Pe$. As discussed in Section \ref{ch_intro}, these domains can be written as $\Gg=\{\Pi_2(A): A \in M_2(\C), r(A)<1 \}, \E=\{\pi_\E(A) : A \in M_2(\C), \mu_{\text{tetra}}(A)<1 \}$ and $\Pe=\{\pi_{\Pe}(A) : A \in M_2(\C), \mu_{\text{penta}}(A)<1\}$, where $\Pi_2, \pi_\E$ and $\pi_\Pe$ are the maps as in \eqref{eqn_Pi_n}, \eqref{eqn_pi_E} and \eqref{eqn_pi_Pe} respectively. Also, $\mu_{\text{tetra}}$ and $\mu_{\text{penta}}$ represent the cost functions in $\mu$-synthesis associated with $\E$ and $\Pe$ respectively, as discussed in Section \ref{ch_intro}.

\smallskip 

In this section, we consider the linear subspace $E$ consisting of all $2 \times 2$ upper triangular complex matrices and we denote $\mu_E$ in this case by $\mu_{\text{hexa}}$. Let us define the map 
\[
\pi: M_2(\C) \to \C^4, \quad \ A=(a_{ij})_{i, j=1}^2 \mapsto (a_{21}, a_{11}, a_{22}, \det(A)).
\] 
Motivated by the aforementioned characterizations of the domains $\Gg, \E$ and $\Pe$, we define 
\begin{equation}\label{eqn1.1}
	\HB=\{\pi(A) : A \in M_2(\C), \mu_{hexa}(A) <1\}.
\end{equation}
We refer to the set $\HB$ as the \textit{$\mu$-hexablock}. In general, computing the value of a cost function $\mu$ is difficult. So, it is natural to seek an easier criteria for a point to belong to $\HB$. Most of this section is devoted to providing a few alternative descriptions of $\HB$, and later we discuss some of its geometric properties. We begin with a basic lemma.

\begin{lem}\label{lem2.3}
	Let $X=\begin{pmatrix}
		z_1 & w\\
		0 & z_2
	\end{pmatrix} \in M_2(\C)$. Then the following hold:
	\begin{enumerate}
		\item $\|X\| \leq 1$ if and only if $|z_1|, |z_2| \leq 1$ and $|w| \leq \sqrt{(1-|z_1|^2)(1-|z_2|^2)}$.
		\item $\|X\|<1$ if and only if $|z_1|, |z_2| < 1$ and $|w| < \sqrt{(1-|z_1|^2)(1-|z_2|^2)}$. 		
	\end{enumerate}		
\end{lem}

\begin{proof}
A Hilbert space operator $T$ has norm at most $1$ if and only if $I-T^*T \geq 0$. Therefore, $\|X\| \leq 1$ if and only if the matrix given by
\[
I-X^*X=\begin{pmatrix}
	1-|z_1|^2  & -\overline{z}_1w \\
	-z_1\overline{w} & 1-|z_2|^2-|w|^2
\end{pmatrix}
\]
is positive, i.e., $I-X^*X \geq 0$, which is possible if and only if $|z_1| \leq 1$ and $\det(I-X^*X) \geq 0$. Note that $\det(I-X^*X)=(1-|z_1|^2)(1-|z_2|^2)-|w|^2$ and so, the part $(1)$ follows. Since $X$ is an operator on $\C^2$, $\|X\|$ is attained at some point of the unit ball in $\C^2$ and so, $\|X\|<1$ if and only if $\|Xh\|<\|h\|$ for every non-zero $h$ in $\C^2$. Combining these things, we have that $\|X\|<1$ if and only if $I-X^*X>0$. Now, the part $(2)$ follows directly from the above computations. 
\end{proof} 

Before coming to the main results of this section, we need to state a few characterizations of the points in the pentablock and its closure. 

\begin{thm}[\cite{AglerIV}, Theorem 5.2]\label{pentablock}
	Let $a \in \C$ and let $(s, p)=(\lm_1+\lm_2, \lm_1\lm_2) \in \Gg$. Then the following are equivalent:	
	\begin{enumerate}
		\item $(a, s, p) \in \mathbb{P}$;
		\item there exists $A=(a_{ij})_{i, j=1}^2 \in M_2(\C)$ with $\mu_{\text{penta}}(A)<1$ and $(a, s, p)=(a_{21}, \text{tr}(A), \det(A))$;
		\item $
		\displaystyle \sup_{z \in \D} \bigg|\displaystyle\frac{a(1-|z|^2)}{1-sz+pz^2}\bigg|<1$;
		\item $
		|a|< \frac{1}{2}|1-\overline{\lambda}_2\lambda_1|+\frac{1}{2}\sqrt{(1-|\lambda_1|^2)(1-|\lambda_2|^2)}$.
	\end{enumerate}
\end{thm}

\begin{thm}[\cite{AglerIV}, Theorem 5.3]\label{pentablock_c}
	Let $a \in \C$ and let $(s, p)=(\lm_1+\lm_2, \lm_1\lm_2) \in \Gamma_2$. Then the following are equivalent:	
	\begin{enumerate}
		\item $(a, s, p) \in \Pbar $;
		\item there exists $A=(a_{ij})_{i, j=1}^2 \in M_2(\C)$ such that $\|A\|\leq 1$ and $(a, s, p)=(a_{21}, \text{tr}(A), \det(A))$;
		\item there exists $A=(a_{ij})_{i, j=1}^2 \in M_2(\C)$ with $\mu_{\text{penta}}(A)\leq 1$ and $(a, s, p)=(a_{21}, \text{tr}(A), \det(A))$;
		\item $
		\displaystyle \sup_{z \in \D} \bigg|\displaystyle\frac{a(1-|z|^2)}{1-sz+pz^2}\bigg|\leq 1$;
		\item $
		|a|\leq  \frac{1}{2}|1-\overline{\lambda}_2\lambda_1|+\frac{1}{2}\sqrt{(1-|\lambda_1|^2)(1-|\lambda_2|^2)}$.
	\end{enumerate}
\end{thm}

As discussed in the beginning of this section, the holomorphic map $\pi_{\E}$ is central to the theory of tetrablock $\E$. Evidently, the last three components of $\pi(A)=(a_{21}, a_{11}, a_{22}, \det(A))$ constitute the element $\pi_\E(A)=(a_{11}, a_{22}, \det(A))$, where $A=(a_{ij})_{i, j=1}^2 \in M_2(\C)$. So, if $\mu_\text{hexa}(A)<1$, then by \eqref{eqn_mu_rel}, we have that $\mu_{\text{tetra}}(A) \leq \mu_{\text{hexa}}(A)<1$. Consequently, $\pi_\E(A) \in \E$ and so, the last three components of an element in $\HB$ form an element in $\E$. Hence, one can capitalize the tetrablock theory to study the $\mu$-hexablock $\HB$. In fact, we find below several characterizations for the elements in $\HB$ using the theory of tetrablock, along with the fractional maps $\Pz$ for $z_1, z_2 \in \D$, as discussed in Section \ref{psi_map}.

\begin{thm}\label{thm2.4}
	For $A=(a_{ij})_{i, j=1}^2 \in M_2(\C)$, the following are equivalent:
	\begin{enumerate}
		\item $\mu_{\text{hexa}}(A) <1$;
		
		\smallskip
		
		\item $
		1-a_{11}z_1-a_{22}z_2+\det(A)z_1z_2 \ne a_{21}w \ $ for $z_1, z_2, w \in \overline{\D}$ and $|w| \leq \sqrt{(1-|z_1|^2)(1-|z_2|^2)}$ ;
		
		\smallskip 
		
		\item $(a_{11}, a_{22}, \det(A)) \in \E$ and $\sup_{z_1, z_2 \in \mathbb{D}}\left|\Pz(a_{21}, a_{11}, a_{22}, \det(A))\right| <1$;
		
		\smallskip			
		
		\item $(a_{11}, a_{22}, \det(A)) \in \E$ and $\left|\Pz(a_{21}, a_{11}, a_{22}, \det(A) )\right| <1$ for every $z_1, z_2 \in \D$;
		
		\smallskip		
		
		\item $(a_{11}, a_{22}, \det(A)) \in \E$ and	$|a_{21}|K_*(a_{11}, a_{22}, \det(A))<1$,
	\end{enumerate}
	where $K_*$ is the map as in \eqref{eqn_K_*}. Therefore,	
	\[
	\HB=	\left \{(a_{21}, a_{11}, a_{22}, \det (A)) \in \C \times \E \ : \ A=(a_{ij})_{i,j=1}^2 \in M_2(\C), \ |a_{21}|K_*(a_{11}, a_{22}, \det(A))<1\right\}.
	\]
\end{thm}

\begin{proof} The equivalences $(3) \iff (4)$ and $(4) \iff (5)$ follow directly from Corollary \ref{cor2.3}. So, it suffices to prove that $(1) \implies (2) \implies (3) \implies (1)$.

\medskip

\noindent $(1) \implies (2).$ If $\mu_{\text{hexa}}(A)<1$, then $\inf\{\|X\| :  X \in E, \ \det(I-AX)=0  \} >1$. Thus, we have $\det(I-AX) \neq 0$ for every $X \in E$ with $\|X\| \leq 1$. Let $z_1, z_2 \in \overline{\D}$ and let $w \in \C$ be such that $|w| \leq \sqrt{(1-|z_1|^2)(1-|z_2|^2)}$.
Define $X=\begin{pmatrix}
	z_1 & w \\
	0 & z_2
\end{pmatrix}$. Then $X \in E$ and we have by Lemma \ref{lem2.3} that $\|X\|\leq 1$. Consequently, $\det(I-AX)=1-a_{11}z_1-a_{22}z_2+\det(A)z_1z_2-a_{21}w \ne 0$.

\medskip

\noindent $(2) \implies (3).$ Assume that $1-a_{11}z_1-a_{22}z_2+\det(A)z_1z_2 \ne a_{21}w$ for every $z_1, z_2, w \in \overline{\D}$ such that $|w| \leq \sqrt{(1-|z_1|^2)(1-|z_2|^2)}$. In particular, if we choose $w=0$ and $z_1, z_2 \in \overline{\D}$, we have that
\begin{equation}\label{eqn2.2}
	1-a_{11}z_1-a_{22}z_2+\det(A)z_1z_2 \neq 0.
\end{equation}
By definition of $\E$, it follows that $(a_{11}, a_{22}, \det(A)) \in \E$. For every $z_1, z_2 \in \overline{\D}$, we show that 
\begin{equation}\label{eqn2.3}
	|1-a_{11}z_1-a_{22}z_2+\det(A)z_1z_2|> |a_{21}|\sqrt{(1-|z_1|^2)(1-|z_2|^2)}
\end{equation}
holds. If $a_{21}=0$, then  \eqref{eqn2.3} follows from \eqref{eqn2.2}. Let $a_{21} \neq 0$ and let $z_1, z_2 \in \overline{\D}$ be such that 

\begin{equation}\label{eqn2.4}
	|1-a_{11}z_1-a_{22}z_2+\det(A)z_1z_2| \leq  |a_{21}|\sqrt{(1-|z_1|^2)(1-|z_2|^2)}. 
\end{equation}
Define $w=a_{21}^{-1}(1-a_{11}z_1-a_{22}z_2+\det(A)z_1z_2)$ and so, $|w| \leq \sqrt{(1-|z_1|^2)(1-|z_2|^2)}$ by \eqref{eqn2.4}. Moreover, we have that $1-a_{11}z_1-a_{22}z_2+\det(A)z_1z_2-a_{21}w=0$ which is a contradiction and so, \eqref{eqn2.3} holds. Consequently, we have
\[
|\Pz(a_{21}, a_{11}, a_{22}, \det(A))|=		\bigg|\frac{a_{21}\sqrt{(1-|z_1|^2)(1-|z_2|^2)}}{1-a_{11}z_1-a_{22}z_2+\det(A)z_1z_2}\bigg| <1
\] 
for all $z_1, z_2 \in \overline{\D}$. The map $(z_1, z_2) \mapsto |\Pz(a_{21}, a_{11}, a_{22}, \det(A))|$ attains its supremum over $\mathbb{D}^2$ as it is a continuous map on $\overline{\D}^2$ and $|\Pz(a_{21}, a_{11}, a_{22}, \det(A))|=0$ for $(z_1, z_2) \in \partial \D^2$. Thus, the desired conclusion is achieved.
\medskip

\noindent $(3) \implies (1).$ Assume that $(a_{11}, a_{22},\det \ A) \in \E$ and 
$\displaystyle \underset{z_1, z_2 \in \mathbb{D}}{\sup}\bigg|\frac{a_{21}\sqrt{(1-|z_1|^2)(1-|z_2|^2)}}{1-a_{11}z_1-a_{22}z_2+\det(A)z_1z_2}\bigg| <1$. Take $X=\begin{pmatrix}
	z_1 & w\\
	0 & z_2\\
\end{pmatrix} \in E$ with $\det(I-AX)=0$. Let if possible, $\|X\| \leq 1$. We have by Lemma \ref{lem2.3} that $|z_1|, |z_2| \leq 1$ and $|w| \leq \sqrt{(1-|z_1|^2)(1-|z_2|^2)}$. Note that $\det(I-AX)=1-a_{11}z_1-a_{22}z_2+\det(A)z_1z_2-a_{21}w$. If $z_1, z_2 \in \partial \D^2$,  then $w=0$ and so, $1-a_{11}z_1-a_{22}z_2+\det(A)z_1z_2=0$ which is not possible as $(a_{11}, a_{22},\det(A)) \in \E$. Thus, $z_1, z_2 \in \D$. Since $\det(I-AX)=0$, we have that
\[
|1-a_{11}z_1-a_{22}z_2+\det(A)z_1z_2|=|a_{21}w| \leq |a_{21}|\sqrt{(1-|z_1|^2)(1-|z_2|^2)}
\]
and thus,
\[
\bigg|\frac{a_{21}\sqrt{(1-|z_1|^2)(1-|z_2|^2)}}{1-a_{11}z_1-a_{22}z_2+\det(A)z_1z_2}\bigg| \geq 1,
\]
contrary to the given hypothesis. Consequently, $\|X\|>1$ whenever $X \in E$ and $\det(I-AX)=0$. Note that the set $\{\|X\| : X \in E, \det(I-AX)=0\}$ is closed and bounded below and so, it attains its infimum. Therefore, $\inf\{\|X\| : X \in E, \det(I-AX)=0\}>1$ which implies that $\mu_{\text{hexa}}(A)<1$.
\end{proof}

Using the above theorem, we obtain a few characterizations for $A \in M_2(\C)$ with $\mu_{\text{hexa}}(A) \leq 1$.

\begin{prop}\label{prop2.5}
	For $A=(a_{ij})_{i, j=1}^2 \in M_2(\C)$, the following are equivalent:
	\begin{enumerate}
		\item $\mu_{\text{hexa}}(A) \leq 1$;
		
		\item $1-a_{11}z_1-a_{22}z_2+\det(A)z_1z_2 \ne a_{21}w \ $ for $z_1, z_2, w \in \D$ with $|w| < \sqrt{(1-|z_1|^2)(1-|z_2|^2)}$;		
		
		\item $(a_{11}, a_{22},\det(A)) \in \overline{\E}$ and $\underset{z_1, z_2 \in \mathbb{D}}{\sup}\left|\Pz(a_{21}, a_{11}, a_{22}, \det(A))\right| \leq 1$.
	\end{enumerate}
\end{prop}

\begin{proof}
We prove that $(1) \implies (2) \implies (3) \implies (1)$.

\noindent $(1) \implies (2).$ Let $\mu_{\text{hexa}}(A) \leq 1$. Then $\inf\{\|X\| :  X \in E, \det(I-AX)=0  \} \geq 1$ and so, it follows that $\det(I-AX) \neq 0$ for every $X \in E$ with $\|X\|<1$. Let $z_1, z_2, w \in \D$ and $|w| < \sqrt{(1-|z_1|^2)(1-|z_2|^2)}$. Define $X=\begin{pmatrix}
	z_1 & w \\
	0 & z_2
\end{pmatrix}$. Then $X \in E$ and Lemma \ref{lem2.3} gives that $\|X\| < 1$. Therefore, $
\det(I-AX)=1-a_{11}z_1-a_{22}z_2+\det(A)z_1z_2-a_{21}w \neq 0$.

\medskip

\noindent $(2) \implies (3)$. Assume that $1-a_{11}z_1-a_{22}z_2+\det(A)z_1z_2 \ne a_{21}w$	for every $z_1, z_2, w$ in $\D$ satisfying $|w| < \sqrt{(1-|z_1|^2)(1-|z_2|^2)}$. Consequently, $1-a_{11}z_1-a_{22}z_2+\det(A)z_1z_2 \neq 0$ for $z_1, z_2 \in \D$ and so, by part-$(2)$ of Theorem \ref{closedtetrablock}, $(a_{11}, a_{22},\det(A)) \in \overline{\E}$. For every $z_1, z_2 \in \D$, we show that 
\begin{equation}\label{eqn2.9}
	|1-a_{11}z_1-a_{22}z_2+\det(A)z_1z_2|\geq |a_{21}|\sqrt{(1-|z_1|^2)(1-|z_2|^2)}
\end{equation}
holds. Evidently, it holds trivially if $a_{21}=0$. Let $a_{21} \neq 0$ and let $z_1, z_2 \in \D$ be such that 
\begin{equation}\label{eqn2.10}
	|1-a_{11}z_1-a_{22}z_2+\det(A)z_1z_2| <  |a_{21}|\sqrt{(1-|z_1|^2)(1-|z_2|^2)}. 
\end{equation}
Define $w=a_{21}^{-1}(1-a_{11}z_1-a_{22}z_2+\det(A)z_1z_2)$ and so, $|w| < \sqrt{(1-|z_1|^2)(1-|z_2|^2)}$ by \eqref{eqn2.10}. Furthermore, $1-a_{11}z_1-a_{22}z_2+\det(A)z_1z_2-a_{21}w=0$ which is a contradiction to the given assumption. Hence, \eqref{eqn2.9} holds which gives the desired conclusion. 

\vspace{0.2cm}

\noindent $(3) \implies (1)$. Let $(a_{11}, a_{22},\det (A)) \in \overline{\E}$ and let 
\begin{equation}\label{eqn2.12}
	\underset{z_1, z_2 \in \mathbb{D}}{\sup}\left|\Pz(a_{21}, a_{11}, a_{22}, \det(A))\right|= 	\underset{z_1, z_2 \in \mathbb{D}}{\sup}\bigg|\frac{a_{21}\sqrt{(1-|z_1|^2)(1-|z_2|^2)}}{1-a_{11}z_1-a_{22}z_2+\det(A)z_1z_2}\bigg| \leq 1.
\end{equation}
It follows from part-$(6)$ of Theorem \ref{tetrablock} and part-$(5)$ of Theorem \ref{closedtetrablock} that $(ra_{11}, ra_{22}, r^2 \det (A)) \in \E$ for every $r \in (0,1)$. Let $z_1, z_2 \in \D$ and let $r \in (0, 1)$. Then
\begin{equation*}
	\begin{split}
		\bigg|\frac{ra_{21}\sqrt{(1-|z_1|^2)(1-|z_2|^2)}}{1-ra_{11}z_1-ra_{22}z_2+r^2\det(A)z_1z_2}\bigg| & 
		\leq 		\bigg|\frac{ra_{21}\sqrt{(1-|rz_1|^2)(1-|rz_2|^2)}}{1-ra_{11}z_1-ra_{22}z_2+r^2\det(A)z_1z_2}\bigg|\\
		& <	\bigg|\frac{a_{21}\sqrt{(1-|rz_1|^2)(1-|rz_2|^2)}}{1-a_{11}(rz_1)-a_{22}(rz_2)+\det(A)(rz_1)(rz_2)}\bigg|\\
		& \leq 1,
	\end{split}
\end{equation*}
where	the last inequality follows from  \eqref{eqn2.12}. We have by Theorem \ref{thm2.4}  that $\mu_{\text{hexa}}(rA) <1$ for every $r \in (0,1)$. It is not difficult to see that $\mu_{hexa}(rA)=r\mu_{hexa}(A)$ and so, $
\mu_{\text{hexa}}(A) < 1\slash r$. Taking limits as $r \to 1$, it follows that $\mu_{\text{hexa}}(A) \leq 1$. The proof is now complete. 
\end{proof} 

It is clear from the above results that if $(a, x_1, x_2, x_3) \in \HB$ (respectively, $\CHB$), then $(x_1, x_2, x_3) \in \E$ (respectively, $\Ebar$). A similar conclusion is also true for $\Pe$ and $\Gg$. Indeed, if $(a, s, p) \in \Pe$ (respectively, $\Pbar$), then $(s, p) \in \Gg$ (respectively, $\Gamma_2$). Furthermore, we have that $\{0\} \times \Gg \subset \Pe$ and $\{0\} \times \Gamma_2 \subset \Pbar$ (see Theorems \ref{pentablock} \& \ref{pentablock_c}). However, this is no longer true for $\HB$ and $\E$ as shown in the following example.  

\begin{eg}\label{eg_1}
	Let $\alpha \in \D\setminus \{0\}$. We have by part-$(3)$ of Theorem \ref{tetrablock} that $(0, 0, \alpha) \in \E$. Let if possible, $(0, 0, 0, \alpha) \in \HB$. By definition, there exists $A=\begin{pmatrix}
		a_{11} & a_{12} \\
		a_{21} & a_{22}
	\end{pmatrix}$ in $M_2(\C)$ such that $\mu_{hexa}(A)<1$ and $\pi(A)=(0, 0, 0, \alpha)$. Thus, $a_{21}=a_{11}=a_{22}=0$ and $\det(A)=\alpha$. This is possible if and only if $\alpha=0$ which is a contradiction. Hence, $\{0\} \times \E$ is not a subset of $\HB$. \qed
\end{eg}

Motivated by the above example, we describe the sets $\HB \cap (\{0\} \times \E)$ and $\HB \setminus (\{0\} \times \E)$. To do so, we need the following lemma. 

\begin{lem}\label{rem3.6}
	For $(a, x_1, x_2, x_3) \in \C^4$ with $a \ne 0$, there exists a unique matrix $A \in M_2(\C)$ with $\pi(A)=(a, x_1, x_2, x_3)$.
\end{lem}

\begin{proof} Let $A=\begin{pmatrix}
	x_1 & a^{-1}(x_1x_2-x_3) \\
	a & x_2
\end{pmatrix}$. Then  $\pi(A)=(a, x_1, x_2, x_3)$. Let $B=(b_{ij})_{i, j=1}^2 \in M_2(\C)$ be such that $\pi(B)=(a, x_1, x_2, x_3)$. Then $
(b_{21}, b_{11}, b_{22}, \det(B))=(a, x_1, x_2, x_3)$ which is possible if and only if $b_{12}=a^{-1}(x_1x_2-x_3)$. Therefore, $A=B$ and we reach the conclusion.
\end{proof} 

As a corollary to Lemma \ref{rem3.6}, we obtain the following description of $\HB \cap (\{0\} \times \E)$. 
\begin{cor}\label{cor2.6}  
	Let $(a, x_1, x_2, x_3) \in \C^4$. Then the following are equivalent.
	\begin{enumerate}
		\item $(a, x_1, x_2, x_3) \in \HB$;
		\item $(x_1, x_2, x_3) \in \E$ and $|a|K_*(x_1, x_2, x_3)<1$. In addition, $x_1x_2=x_3$ if $a=0$;
		
		\item  $(x_1, x_2, x_3) \in \E$ and $\underset{z_1, z_2 \in \D}{\sup}|\psi_{z_1, z_2}(a, x_1, x_2, x_3)|<1$. In addition, $x_1x_2=x_3$ if $a=0$.
	\end{enumerate}
	Therefore, $\HB \cap (\{0\} \times \E)=\{(0, x_1, x_2, x_3) \in \HB : x_1x_2=x_3\}$
\end{cor}

\begin{proof} The equivalence of $(2)$ and $(3)$ follows from Corollary \ref{cor2.3}.
\medskip 

\noindent $(1) \implies (2)$. Let $(a, x_1, x_2, x_3) \in \HB$. By Theorem \ref{thm2.4}, $(x_1, x_2, x_3) \in \E$ and $|a|K_*(x_1, x_2, x_3)<1$. By definition, there is a $2 \times 2$ matrix $A$ with $\mu_{\text{hexa}}(A) <1$ such that $\pi(A)=(a, x_1, x_2, x_3)$. If $a=0$, then $x_3=\det(A)=x_1x_2$.

\medskip	

\noindent $(2) \implies (1)$. Let $(x_1, x_2, x_3) \in \E$ and $|a|K_*(x_1, x_2, x_3)<1$. It follows from Corollary \ref{cor2.3} that $\underset{z_1, z_2 \in \D}{\sup}|\psi_{z_1, z_2}(a, x_1, x_2, x_3)|<1$. If $a \ne 0$, then it follows from Lemma \ref{rem3.6} that there exists a unique $A \in M_2(\C)$ such that $\pi(A)=(a, x_1, x_2, x_3)$. For $a=0$ and $x_1x_2=x_3$, let us define $A=\begin{pmatrix}
	x_1 & 0 \\
	0 & x_2
\end{pmatrix}$ so that $\pi(A)=(a, x_1, x_2, x_3)$. In either case, it follows from Theorem \ref{thm2.4} that $\mu_{hexa}(A)<1$ and so, $\pi(A)=(a, x_1, x_2, x_3) \in \HB$. The proof is now complete.
\end{proof} 

So, a description of $\HB \cap (\{0\} \times \E)$ is obtained. Now we proceed ahead and identify the set $\HB \setminus (\{0\} \times \E)$. In fact, we shall prove in Proposition \ref{prop2.9} that $int(\HB)=\HB \setminus (\{0\} \times \E)$, the interior of $\HB$. This immediately shows that $\HB$ is not an open subset of $\C^4$. To prove these results, we first recall the following properties of $\mu_E$ relative to a linear subspace $E$ of $M_n(\C)$:
\begin{enumerate}
	\item $\mu_E$ is upper semicontinuous,
	\item $\mu_E(\lambda A)=|\lambda|\mu_E(A)$ for any $\lambda \in \C, A\in M_n(\C)$,
	\item $\mathbb{B}_{\mu_E}=\{A \in M_n(\C) : \mu_E(A) <1\}$ is a domain.
\end{enumerate}
Recall that in the present context, $E$ is the linear space of all upper triangular matrices in $M_2(\C)$. An interested reader can refer to Remark 2.2.1 in \cite{Jarnicki} (also see the discussion followed by Remark 3.1 in \cite{Zap}) for further details. 

\begin{prop}\label{prop2.7}
	$\HB$ is connected but not open in $\C^4$. Also, if $(x_1, x_2, x_1x_2) \in \E$, then $(0, x_1, x_2, x_1x_2)$ is not an interior point of $\HB$.
\end{prop}

\begin{proof}
$\HB$ is connected since it is the image of the domain $\mathbb B_{\mu_E}$ under the continuous map $\pi$. Let $(x_1, x_2, x_1x_2) \in \E$. It follows from Corollary \ref{cor2.6} that $(0, x_1, x_2, x_1x_2) \in \HB$. Let if possible, $(0, x_1, x_2, x_1x_2)$ be an interior point of $\HB$. Then there exists some $\epsilon>0$ such that 
\[
(0, x_1, x_2, x_1x_2) \in B(0, \epsilon) \times B(x_1, \epsilon) \times B(x_2, \epsilon) \times B(x_1x_2, \epsilon) \subset \HB,
\]
where $B(z_0, r)=\{z: |z-z_0|<r\}$ for $z_0 \in \C$ and $r>0$. Thus, $(0, x_1, x_2, x_1x_2+\epsilon \slash 2) \in \HB$ which is a contradiction to Corollary \ref{cor2.6}. Consequently, $\HB$ is not open.
\end{proof}

It follows from Corollary \ref{cor2.6} and Proposition \ref{prop2.7} that no point of $\HB \cap (\{0\} \times \E)$ is an interior point of $\HB$. Indeed, we show that all other remaining points in $\HB$ are its interior points.

\begin{prop}\label{prop2.9}
	$int(\HB)=\left\{(a, x_1, x_2, x_3) \in \HB \ : \ a \ne 0\right\}=\HB \setminus (\{0\} \times \E)$.
\end{prop}

\begin{proof}
Consider the function $
\phi: (\C\setminus \{0\}) \times \E \to \mathbb{R}$ defined as $\phi(a, x_1, x_2, x_3)=|a|K_*(x_1, x_2, x_3)$. Then $\phi^{-1}(-\infty, 1)=\{(a, x_1, x_2, x_3) \in (\C\setminus \{0\}) \times \E : |a|K_*(x_1, x_2, x_3) <1\}$. We have by Corollary \ref{cor2.6} that 
$
\phi^{-1}(-\infty, 1)=\{(a, x_1, x_2, x_3) \in \HB : a \ne 0\}.
$
Since $\phi$ is a continuous function and $(\C\setminus\{0\}) \times \E$ is open in $\C^4$, it follows that $\phi^{-1}(-\infty, 1)$ is an open set in $\C^4$ contained in $\HB$. Consequently, $\phi^{-1}(-\infty, 1) \subseteq int(\HB) \subseteq \HB$ and $\HB \setminus \phi^{-1}(-\infty, 1)=\HB \cap (\{0\} \times \E)$. It follows from Proposition \ref{prop2.7} that no point in $\HB \setminus \phi^{-1}(-\infty, 1)$ belongs to $int(\HB)$ and so, $int(\HB)=\phi^{-1}(-\infty, 1)=\left\{(a, x_1, x_2, x_3) \in \HB \ : \ a \ne 0\right\}$. It is only left to show that \[
\left\{(a, x_1, x_2, x_3) \in \HB \ : \ a \ne 0\right\} =\HB\setminus (\{0\} \times \E).
\]
Evidently, if $(a, x_1, x_2, x_3) \in \HB$ with $a \ne 0$, then $(a, x_1, x_2, x_3) \in \HB\setminus (\{0\} \times \E)$. Conversely, let $(a, x_1, x_2, x_3) \in \HB\setminus (\{0\} \times \E)$. Since $(a, x_1, x_2, x_3) \notin \{0\} \times \E$, it is obvious that either $a \ne 0$ or $(x_1, x_2, x_3) \notin \E$. As $(a, x_1, x_2, x_3) \in \HB$, Corollary \ref{cor2.6} shows that $(x_1, x_2, x_3) \in \E$ and so, $a \ne 0$. The proof is now complete. 
\end{proof} 

Thus, Corollary \ref{cor2.6} reveals the crucial fact that $\{0\} \times \E$ is not contained in $\HB$. However, it turns out that $\{0\} \times \Ebar \subset \CHB$ as the following result explains. 

\begin{thm}\label{CHB}
	The closure of $\HB$ is given by	
	\[
	\CHB=\left\{(a, x_1, x_2, x_3) \in \C \times \overline{\E} \ : \ \bigg|\frac{a\sqrt{(1-|z_1|^2)(1-|z_2|^2)}}{1-x_1z_1-x_2z_2+x_3z_1z_2}\bigg| \leq 1 \ \text{for every} \ z_1, z_2 \in \D \right\}.
	\]	
	Moreover, $\overline{int(\HB)}=\CHB$ and $\left\{\pi(A) : \ A \in M_2(\C), \mu_{\text{hexa}}(A) \leq 1 \right\} \subseteq  \CHB$. 
\end{thm}

\begin{proof}
Let $\HB^{(0)}=\left\{\pi(A) : \ A \in M_2(\C), \mu_{\text{hexa}}(A) \leq 1 \right\}$ and let
\[
\HB^{(1)}=\left\{(a, x_1, x_2, x_3) \in \C \times \overline{\E} \ : \ \bigg|\frac{a\sqrt{(1-|z_1|^2)(1-|z_2|^2)}}{1-x_1z_1-x_2z_2+x_3z_1z_2}\bigg| \leq 1 \ \text{for every} \ z_1, z_2 \in \D \right\}.
\]
We prove that $\HB^{(0)} \subseteq \CHB \subseteq \HB^{(1)} \subseteq \CHB$.	Let $A=(a_{ij})_{i, j=1}^2$ with $\mu_{\text{hexa}}(A) \leq 1$. For $n \in \mathbb{N}$, let us define $A_n=\frac{n-1}{n}A$. Then $\mu_{\text{hexa}}(A_n) \leq 1-1\slash n <1$ and so, $\pi(A_n) \in \HB$. It is easy to see that $\pi(A_n) \to \pi (A)$ as $n \to \infty$. Therefore, $\pi(A) \in \CHB$  which implies that $\HB^{(0)} \subset \CHB$. 

\smallskip 

Let $(a, x_1, x_2, x_3) \in \CHB$. Then there is a sequence $\{(a^{(n)}, x_1^{(n)}, x_2^{(n)}, x_3^{(n)})\}$ in $\HB$  that converges to  $(a, x_1, x_2, x_3)$ as $n \to \infty$. It follows straight-forward from Theorem \ref{thm2.4} that $(x_1^{(n)}, x_2^{(n)}, x_3^{(n)}) \in \E$ for every $n \in \mathbb{N}$ and so, $(x_1, x_2, x_3) \in \overline{\E}$. For every $n \in \mathbb{N}$ and $z_1, z_2 \in \D$, we have by Theorem \ref{thm2.4} that 
\[
\bigg|\frac{a^{(n)}\sqrt{(1-|z_1|^2)(1-|z_2|^2)}}{1-x_1^{(n)}z_1-x_2^{(n)}z_2+x_3^{(n)}z_1z_2}\bigg| < 1 
\quad \text{and so,} \quad 
\bigg|\frac{a\sqrt{(1-|z_1|^2)(1-|z_2|^2)}}{1-x_{1}z_1-x_{2}z_2+x_{3}z_1z_2}\bigg| \leq 1.
\]
Thus, $(a, x_1, x_2, x_3) \in \HB^{(1)}$ and so, $\CHB \subseteq \HB^{(1)}$.

\smallskip

Let $(a, x_1, x_2, x_3) \in \HB^{(1)}$. Then $(x_1, x_2, x_3) \in \Ebar$. If $a \ne 0$, then Lemma \ref{rem3.6} ensures the existence of $A \in M_2(\C)$ such that $\pi(A)=(a, x_1, x_2, x_3)$. It follows from Proposition \ref{prop2.5} that $\mu_{\text{hexa}}(A) \leq 1$ and so, $(a, x_1, x_2, x_3) \in \HB^{(0)}$. Therefore, $(a, x_1, x_2, x_3) \in \CHB$ when $a \ne 0$. Assume that $a=0$. We proceed from here with two cases, depending on whether $(x_1, x_2, x_3)$ is in $\E$ or in $\partial \E$. Assume that $(x_1, x_2, x_3) \in \E$. The discussion above Corollary \ref{cor2.3} implies that $K_*(x_1, x_2, x_3)$ is a finite positive number. Choose $N \in \mathbb{N}$ such that $K_*(x_1, x_2, x_3)<N$. It follows from Proposition \ref{prop2.9} that $\left(1 \slash n, x_1, x_2, x_3 \right) \in int(\HB)$ for every $n \geq N$ and thus, $(0, x_1, x_2, x_3) \in \overline{int(\HB)}$. This also shows that $\{0\} \times \E \subseteq \CHB$ and so, $\{0\} \times \Ebar \subseteq \CHB$. Consequently, $(0, x_1, x_2, x_3) \in \CHB$ when $(x_1, x_2, x_3) \in \partial \E$. Combining everything together, we have proved that $\HB^{(1)} = \CHB$. It only remains to prove that $\CHB=\overline{int(\HB)}$. It is trivial to show that $\overline{int(\HB)} \subseteq \CHB$. Let $(a, x_1, x_2, x_3) \in \CHB$. By above discussion, $(x_1, x_2, x_3) \in \Ebar$ and $\{0\} \times \overline{\E} \subset \overline{int(\HB)}$. Thus, we have that $(0, x_1, x_2, x_3) \in \overline{int(\HB)}$. Suppose that $a \ne 0$. Choose a sequence $\{(a^{(n)}, x_1^{(n)}, x_2^{(n)}, x_3^{(n)})\}$ in $\HB$ with $a^{(n)} \ne 0$ for every $n$ such that $\lim_{n \to \infty}
(a^{(n)}, x_1^{(n)}, x_2^{(n)}, x_3^{(n)})= (a, x_1, x_2, x_3)$. It follows from Proposition \ref{prop2.9} that $(a^{(n)}, x_1^{(n)}, x_2^{(n)}, x_3^{(n)}) \in int(\HB)$ and so, $(a, x_1, x_2, x_3) \in \overline{int(\HB)}$. The proof is now complete.
\end{proof}

So, Theorem \ref{CHB} concludes that $\HB^{(0)}=\left\{\pi(A) : \ A \in M_2(\C), \mu_{\text{hexa}}(A) \leq 1 \right\} \subseteq  \CHB$. Needless to mention that $\HB^{(0)} \ne \CHB$. Similar arguments from Example \ref{eg_1} show that if $\alpha \in \overline{\D} \setminus \{0\}$, then $(0, 0, \alpha) \in \Ebar$ and there is no $A \in M_2(\C)$ with $\pi(A)=(0, 0, 0, \alpha)$. Thus, $(0, 0, 0, \alpha) \notin \HB^{(0)}$ but $(0, 0, 0, \alpha) \in \CHB$ which follows from Theorem \ref{CHB}. We now prove some interesting geometric properties of $\HB$ for which we need to recall a few definitions and facts from the literature. 

\begin{defn}
	For a bounded set $\Omega$ in $\C^n$, we say that
	\begin{enumerate}
		\item $\Omega$ is \textit{quasi-balanced} if there exist relatively prime natural numbers $m_1, \dotsc, m_n$ such that $(\lambda^{m_1}z_1, \dotsc, \lambda^{m_n}z_n) \in \Omega$ for every $\lambda \in \overline{\D}$ and $(z_1, \dotsc, z_n) \in \Omega$. Sometimes $\Omega$ is also referred to as \textit{$(m_1, \dotsc, m_n)$-quasi-balanced};
		
		\item $\Omega$ is \textit{starlike about a point $z_0 \in \Omega$} if whenever a point belongs to $\Omega$, then the straight line segment joining $z_0$ to the point is contained in $\Omega$;
		
		\item $\Omega$ is \textit{circled} if $cz \in \Omega$ for every $z \in \Omega$ and $c \in \T$.
	\end{enumerate}
\end{defn}

\begin{defn}
	For a compact set $K \subset \C^n$, the \textit{polynomial convex hull} $\widehat{K}$ of $K$ is defined as
	\[
	\widehat{K}=\left\{z \in \C^n \ : \ |p(z)| \leq \|p\|_{\infty, K}  \ \ \text{for all} \ p \in \C[z_1, \dotsc, z_n] \right\},
	\]
	where $\|p\|_{\infty, K}=\sup\{|p(z)|: z \in K\}$. 	The set $K$ is said to be \textit{polynomially convex} if $\widehat{K}=K$ or equivalently, if for every $z \notin K$, there exists $p \in \C[z_1, \dotsc, z_n]$ such that $\|p\|_{\infty, K} \leq 1<|p(z)|$. 
\end{defn}

The subsequent result was established separately for $\Gg$, $\E$, and $\Pe$ in \cite{Abouhajar}, \cite{AglerIII} and \cite{AglerIV} respectively.

\begin{thm}\label{thm_gp}
The sets $\Gamma_2, \Ebar$, $\Pbar$ are polynomially convex and starlike about the origin but not circled.  
\end{thm}

\begin{prop}\label{prop2.8}
	$\HB$ is $(1,1,1,2)$-quasi-balanced but is neither starlike about $(0,0,0,0)$ nor circled. Also, $\HB$ and $\HB \cap \mathbb{R}^4$ are not convex.
\end{prop}

\begin{proof} The quasi-balanced property follows from the fact that $\pi(zA)=(za, zx_1, zx_2, z^2x_3)$ whenever $\pi(A)=(a, x_1, x_2, x_3)$ for $A \in M_2(\C)$ and $z \in \C$.  For any $x_1, x_2 \in \D \setminus \{0\}$, it follows from Corollary \ref{cor2.6} that $x=(0, x_1, x_2, x_1x_2) \in \HB$. For any $r \in (0,1)$, we have by Corollary \ref{cor2.6} that $rx=(0, rx_1, rx_2, rx_1x_2) \notin \HB$. Thus, $\HB$ is not starlike about $(0,0,0,0)$. By Theorem \ref{CHB}, the point $x=(0,1,1,1) \in \CHB$ but $ix=(0, i, i, i) \notin \CHB$, because $(i, i, i) \notin \Ebar$ which follows from Theorem \ref{closedtetrablock}. Therefore, $\CHB$ is not circled and so, $\HB$ is also not circled. It only remains to show the non-convexity of $\HB$ and $\HB \cap \mathbb{R}^4$. We have by Corollary \ref{cor2.6} that $x=(0, \frac{1}{2}, \frac{1}{2}, \frac{1}{4}), y=(0, \frac{-1}{2}, \frac{-1}{2}, \frac{1}{4})$ are in $\HB$ but 
$
(x+y)\slash 2=\left(0, 0, 0, \frac{1}{4}\right)
$ is not in $\HB$. 
\end{proof}  

\begin{thm}\label{thm2.9}
	$\CHB$ is a polynomially convex set.
\end{thm} 

\begin{proof}
Let $x=(a, x_1, x_2, x_3) \in \C^4 \setminus \CHB$. We must find a polynomial $f$ in 4-variables such that $\|f\|_{\infty, \CHB} \leq 1<|f(x)|$. Let $(x_1, x_2, x_3) \notin \overline{\E}$. We have by Theorem \ref{thm_gp} that $\overline{\E}$ is polynomially convex and so, there is a polynomial $g$ in three variables such that $\|g\|_{\infty, \Ebar} \leq 1 < |g(x_1, x_2, x_3)|$. For any $(b, y_1, y_2, y_3) \in \CHB$, the point $(y_1, y_2, y_3) \in \overline{\E}$ which follows from Theorem \ref{CHB}. Consequently, $f(b, y_1, y_2, y_3) =g(y_1, y_2, y_3)$ is a polynomial that separates $x$ from $\CHB$. Now suppose that $(x_1, x_2, x_3) \in \overline{\E}$. If $|a|>1$, then $f(b, y_1, y_2, y_3)=b$ has the desired properties. Here, we further assume that $|a|\leq 1$ with $(x_1, x_2, x_3) \in \overline{\E}$ and $(a, x_1, x_2, x_3) \notin \CHB$ . By Theorem \ref{CHB}, there must be some $(z_1, z_2) \in \D^2$ such that $|\psi_{z_1, z_2}(a, x_1, x_2, x_3)|>1$. For this fixed $z_1, z_2$ in $\D$, we have $|\psi_{z_1, z_2 }(b, y_1, y_2, y_3)| \leq 1$ for every $(b, y_1, y_2, y_3)$ in $\CHB$. Thus, $\psi_{z_1, z_2}$ is an analytic function on $\CHB$ that separates $x$ from $\CHB$. The idea is to approximate the function $\psi_{z_1, z_2}$ by a sequence of polynomials in the sup-norm topology. We shall see that a polynomial sufficiently close to $\psi_{z_1, z_2}$ with respect to the sup-norm metric is the desired polynomial. For $y=(b, y_1, y_2, y_3) \in \overline{\D} \times \Ebar$, we have  
\begin{equation*}
	\begin{split}
		\psi_{z_1, z_2}(y)
		&=\frac{b\sqrt{(1-|z_1|^2)(1-|z_2|^2)}}{1-y_{1}z_1-y_{2}z_2+y_{3}z_1z_2}\\
		&=\frac{b\sqrt{(1-|z_1|^2)(1-|z_2|^2)}}{\displaystyle (1-y_{2}z_2)\bigg[1-z_1\left(\frac{y_3z_2-y_1}{y_2z_2-1}\right)\bigg]}\\
		&=\frac{b\sqrt{(1-|z_1|^2)(1-|z_2|^2)}}{(1-y_{2}z_2)\bigg[1-z_1\Psi(z_2, y_1, y_2, y_3)\bigg]}\\
		&=b\sqrt{(1-|z_1|^2)(1-|z_2|^2)}\bigg(\overset{\infty}{\underset{j=0}{\sum}}y_2^jz_2^j\bigg)\bigg(\overset{\infty}{\underset{k=0}{\sum}}(\psi(z_2, y_1, y_2, y_3))^kz_1^k\bigg),
	\end{split}
\end{equation*}
where $\Psi(.,y_1,y_2,y_3)$ is the map as in \eqref{eqn_302}.	The two series in the last equality are well-defined as $|y_2z_2| \leq |z_2| <1$ and $|z_1\Psi(z_2, y_1, y_2, y_3)| \leq |z_1|<1$ for every $(y_1, y_2, y_3) \in \overline{\E}$. Then 
\begin{equation}\label{eqn2.14}
	\bigg|(1-z)^{-1}(1-w)^{-1}-\overset{N}{\underset{j=0}{\sum}}z^j\overset{N}{\underset{k=0}{\sum}}w^k \bigg| \leq \frac{|z|^{N+1}+|w|^{N+1}}{(1-|z|)(1-|w|)}
\end{equation}
for any $z,w$ in $\D$ and $N \geq 1$. We define a sequence of functions $g_N$ analytic on $\overline{\D} \times \Ebar$ as follows:
\[
g_N(b, y_1, y_2, y_3)=b\sqrt{(1-|z_1|^2)(1-|z_2|^2)}\bigg(\overset{N}{\underset{j=0}{\sum}}y_2^jz_2^j\bigg)\bigg(\overset{N}{\underset{k=0}{\sum}}(\Psi(z_2, y_1, y_2, y_3))^kz_1^k\bigg).
\]	
Let $M=\sqrt{(1-|z_1|^2)(1-|z_2|^2)}$ and let $|z|=\max\{|z_1|, |z_2|\}$. Again, for $y=(b, y_1, y_2, y_3) \in \overline{\D} \times \Ebar$, suppose $y_0=(y_1, y_2, y_3)$ and $\Psi(., y_0)=\Psi(., y_1, y_2, y_3)$ so that the following holds for each $g_N$: 
\begin{equation*}
	\begin{split}
		|(g_N-\psi_{z_1, z_2}(y)| & = M|b|\ \bigg|(1-y_2z_2)^{-1}(1-z_1\Psi(z_2, y_0))^{-1}-\bigg(\overset{N}{\underset{j=0}{\sum}}y_2^jz_2^j\bigg)\bigg(\overset{N}{\underset{k=0}{\sum}}(\Psi(z_2, y_0))^kz_1^k\bigg)\bigg|\\
		& \leq M|b|\frac{|z_2y_2|^{N+1}+|z_1\Psi(z_2, y_0)|^{N+1}}{(1-|z_2y_2|)(1-|z_1\Psi(z_2, y_0)|)} \quad [\text{by } \eqref{eqn2.14}]\\
		& \leq M|b|\frac{|z_2|^{N+1}+|z_1|^{N+1}}{(1-|z_2y_2|)(1-|z_1\Psi(z_2, y_0)|)} \quad [\text{since } |y_2|, |\Psi(z_2, y_0)| \leq 1]\\
		& \leq M|b|\frac{|z_2|^{N+1}+|z_1|^{N+1}}{(1-|z_2|)(1-|z_1|)} \\
		& \leq \frac{2|z|^{N+1}}{1-|z|^2}.
	\end{split}
\end{equation*}
Thus, 
$
\displaystyle \sup_{y \in \overline{\D} \times \Ebar} \ |g_N(y)-\psi_{z_1, z_2}(y)| \leq \frac{2|z|^{N+1}}{1-|z|^2}.
$
Let $0< \epsilon_0 < \frac{1}{3}(|\psi_{z_1, z_2}(a, x_1, x_2, x_3)|-1)$ and choose $N$ sufficiently large so that 
$
|g_N-\psi_{z_1, z_2}|< \epsilon_0 
$
at all points $(b, x_1, x_2, x_3) \in \overline{\D} \times \Ebar$. Then $|g_N| \leq  1+ \epsilon_0$ on $\CHB$ and $|g_N(a, x_1, x_2, x_3)| > 1+2\epsilon_0$. For this fixed value of $N$, we again approximate the function $g_N$ by a sequence of polynomials. Let 
$
h_m(y_1, y_2, y_3)=(y_1-y_3z_2)(1+y_2z_2+\dotsc + y_2^mz_2^m).
$
Then for $y_0=(y_1, y_2, y_3) \in \overline{\E}$, it follows from \eqref{eqn2.14} that
\begin{equation}\label{eqn2.15}
	|h_m(y_0)-\Psi(z_2, y_0)|=\bigg|(y_1-y_3z_2)\bigg(\overset{m}{\underset{j=0}{\sum}}y_2^jz_2^j-(1-y_2z_2)^{-1}\bigg)\bigg| \leq \frac{2|z_2|^{m+1}}{1-|z_2|}.
\end{equation}
Since $|\Psi(z_2, y_0)| \leq 1$ for $y_0 \in \overline{\E}$, we have 
\begin{equation}\label{eqn2.16}
	\sup_{y_0 \in \Ebar} \ |h_m(y_0)| \leq 1+ \frac{2}{1-|z_2|} \quad \text{for every} \quad m \in \mathbb{N}.
\end{equation}
We now define a sequence of polynomials in the following way:
\[
p_m(b, y_1, y_2, y_3)=b\sqrt{(1-|z_1|^2)(1-|z_2|^2)}\bigg(\overset{N}{\underset{j=0}{\sum}}y_2^jz_2^j\bigg)\bigg(\overset{N}{\underset{k=0}{\sum}}(h_m(z_2, y_1, y_2, y_3))^kz_1^k\bigg).
\]
Let $K=\overset{N}{\underset{j=0}{\sum}}\bigg(1+\frac{2}{1-|z_2|}\bigg)^j$. For any $y_0 \in \overline{\E}$ and $1 \leq k \leq N$, note that 
\begin{equation*}
	\begin{split}
		|\Psi(z_2, y_0)^k-h_m(y_0)^k| & =|\Psi(z_2, y_0)-h_m(y_0)| \ |\Psi(z_2, y_0)^{k-1}+ \Psi(z_2, y_0)^{k-2}h_m(y_0)+ \dotsc +h_m(y_0)^{k-1}| \\
		& \leq |\Psi(z_2, y_0)-h_m(y_0)| \ (1+ |h_m(y_0)|+ \dotsc +|h_m(y_0)|^{k-1})\\
		& \leq |\Psi(z_2, y_0)-h_m(y_0)| \ (1+ |h_m(y_0)|+ \dotsc +|h_m(y_0)|^{N})\\
		& \leq \frac{2|z_2|^{m+1}}{1-|z_2|} \ \bigg[1+ \bigg(1+\frac{2}{1-|z_2|}\bigg)+ \dotsc +\bigg(1+\frac{2}{1-|z_2|}\bigg)^{N}\bigg]\\
		&=2K\frac{|z_2|^{m+1}}{1-|z_2|},
	\end{split}
\end{equation*}
where the last two inequalities follows from  \eqref{eqn2.15} \& \eqref{eqn2.16}. For $y=(b, y_1, y_2, y_3) \in \overline{\D} \times \Ebar$, let $y_0=(y_1, y_2, y_3)$. Then,
\begin{equation*}
	\begin{split}
		|(p_m-g_N)(b, y_1, y_2, y_3)|
		& \leq M|b|  \bigg(\overset{N}{\underset{j=0}{\sum}}|y_2|^j|z_2|^j\bigg)\bigg(\overset{N}{\underset{k=0}{\sum}}|\Psi(z_2, y_0)^k-h_m(y_0)^k| \ |z_1|^k\bigg) \\
		&		\leq  M|b|(N+1) \bigg(\overset{N}{\underset{k=0}{\sum}}2K\frac{|z_2|^{m+1}}{1-|z_2|} \ |z_1|^k\bigg) \\
		&   \leq 2KM|b|(N+1) \frac{|z_2|^{m+1}}{(1-|z_1|)(1-|z_2|)}.\\
	\end{split}
\end{equation*}
Since $|b| \leq 1$, we have  
\[
\sup_{y \in \overline{\D} \times \Ebar} \ |p_m(y)-g_N(y)| \leq 2KM(N+1) \frac{|z_2|^{m+1}}{(1-|z_1|)(1-|z_2|)}
\]
which tends to $0$ as $m \to \infty$. Let $0< \epsilon < \frac{1}{3}(|g_N(x_1, x_2, x_3)|-1-2\epsilon_0)$ and choose $M_0$ sufficiently large so that 
$
|p_{M_0}-g_N|< \epsilon 
$
at all points $(b, x_1, x_2, x_3) \in \overline{\D} \times \Ebar$. Then $|p_{M_0}| \leq  1+ \epsilon+ \epsilon_0$ on $\CHB$ and $|p_{M_0}(a, x_1, x_2, x_3)| > 1+2 \epsilon+ 2\epsilon_0$. The polynomial 
$
f(b, y_1, y_2, y_3)=(1+\epsilon + \epsilon_0)^{-1}p_{M_0}(b, y_1, y_2, y_3)
$
has the desired properties. Thus, $\CHB$ is polynomially convex.
\end{proof}

\section{The normed hexablock}\label{norm_hexa}

\noindent As mentioned in the previous sections, the domains $\Gg$, $\E$ and $\Pe$ are the domains that arise while studying certain instances of the $\mu$-synthesis problem in $M_2(\C)$. We discussed in Section \ref{ch_intro} that these domains can be written as 
\begin{enumerate}
	\item $\Gg=\{\Pi_2(A): A \in M_2(\C), \|A\|<1 \}=\{\Pi_2(A): A \in M_2(\C), r(A)<1 \}$, 
	\item $\E=\{\pi_\E(A) : A \in M_2(\C), \|A\|<1 \}=\{\pi_\E(A) : A \in M_2(\C), \mu_{\text{tetra}}(A)<1 \}$,
	\item $\Pe=\{\pi_{\Pe}(A) : A \in M_2(\C), \|A\|<1\}=\{\pi_{\Pe}(A) : A \in M_2(\C), \mu_{\text{penta}}(A)<1\}$,
\end{enumerate}
where $\Pi_2, \pi_\E$ and $\pi_\Pe$ are the maps as in \eqref{eqn_Pi_n}, \eqref{eqn_pi_E} and \eqref{eqn_pi_Pe} respectively. Also, recall that the maps $\mu_{\text{tetra}}$ and $\mu_{\text{penta}}$ represent the cost functions in $\mu$-synthesis problems associated with $\E$ and $\Pe$, respectively, as seen in Section \ref{ch_intro}. Motivated by the above descriptions of $\Gg, \E$ and $\Pe$, we consider the set given by
\[
\QN=\{\pi(A) : A \in M_2(\C), \|A\|<1 \},
\]
referred to as the \textit{normed hexablock}, where $\pi(A)=(a_{21}, a_{11}, a_{22}, \det(A))$ for $A=(a_{ij})_{i, j=1}^2$ as in \eqref{eqn_pi}. Indeed, $\QN$ is a norm analogue of the $\mu$-hexablock $\HB=\{\pi(A) : A \in M_2(\C), \mu_{\text{hexa}}(A)<1\}$ (see Section \ref{mu_hexa}) in the sense that $\QN$ is obtained from $\HB$ with the condition $\mu_{\text{hexa}}(A)<1$ being replaced by $\|A\|<1$.

\smallskip 

In this section, our central object of study is the normed hexablock $\QN$. We discuss its characterizations, geometric properties, and also its connection with the $\mu$-hexablock $\HB$. Unlike the aforementioned characterizations of $\Gg, \E$ and $\Pe$, we prove that $\QN$ is strictly contained in $\HB$. In fact, we show that $\HB \setminus \QN$ contains a non-empty open set in $\C^4$. We begin with a lemma.

\begin{lem}\label{lem4.1}
	There exists $A=(a_{ij})_{i, j=1}^2 \in M_2(\C)$ such that $a_{21} \ne 0$ and $\mu_{\text{hexa}}(A) < \|A\|$.	
\end{lem}

\begin{proof}
For any $A \in M_2(\C)$, it follows that $\mu_{hexa}(A) \leq \|A\|$. Let if possible, $\|A\|=\mu_{\text{hexa}}(A)$ for every $A=(a_{ij})_{i, j=1}^2 \in M_2(\C)$ with $a_{21} \ne 0$. Let $A=\begin{pmatrix}
	a_{11} & a_{12}\\
	0 & a_{22}\\
\end{pmatrix}$ and define $A_n=\begin{pmatrix}
	a_{11} & a_{12}\\
	1 \slash n & a_{22}\\
\end{pmatrix}$ for $n \in \mathbb{N}$. Since $\mu_{\text{hexa}}$ is upper semi-continuous (see the discussion prior to Proposition \ref{prop2.7}) and $\|A_n\|=\mu_{\text{hexa}}(A_n)$, we have that 
\[
\|A\|=\lim_{n \to \infty}\|A_n\|=\lim_{n \to \infty}\mu_{hexa}(A_n)=		\limsup_{n} \ \mu_{\text{hexa}}(A_n) \leq \mu_{\text{hexa}}(A) \leq \|A\|.
\]
Thus, $\mu_{\text{hexa}}(A)=\|A\|$ for all $A \in M_2(\C)$. Let $A=\begin{pmatrix}
	0 & 5 \\
	0 & 0
\end{pmatrix}$. Then $\|A\|=5$. We have by Theorem \ref{thm2.4} that $\mu_{\text{hexa}}(A)<1$ which is a contradiction. The desired conclusion now follows.
\end{proof}

\begin{thm}\label{thm4.2}
	$\QN$ is a proper subset of $\HB$.
\end{thm}

\begin{proof}
Let $x \in \QN$. Then $x=\pi(A)$ for some $A \in M_2(\C)$ with $\|A\|<1$. Thus, $\mu_{\text{hexa}}(A) \leq \|A\|<1$ and so, $x \in \HB$. Therefore, $\QN \subseteq \HB$. We show that this containment is strict. By Lemma \ref{lem4.1}, there exists $A=(a_{ij})_{i, j=1}^2 \in M_2(\C)$ such that $a_{21} \ne 0$ and $\mu_{\text{hexa}}(A) < \|A\|$. One can choose $r>0$ such that 
$
\mu_{\text{hexa}}(A)<r<\|A\|
$. Then $
\mu_{\text{hexa}}(r^{-1}A)<1<\|r^{-1}A\|
$ and so, $\pi(r^{-1}A) \in \HB$. Let if possible, $\pi(r^{-1}A) \in \QN$. By the definition of $\QN$, there exists $B=(b_{ij})_{i, j=1}^2 \in M_2(\C)$ such that $\|B\|<1$ and $\pi(B)=\pi(r^{-1}A)$. We have by Lemma \ref{rem3.6} that $B=r^{-1}A$, which is not possible as $\|B\|<1<\|r^{-1}A\|$. Therefore, $\pi(r^{-1}A) \in \HB \setminus \QN$.
\end{proof}		

Building on the proof of Lemma \ref{lem4.1}, the next example explicitly provides a matrix that satisfies its conditions, thereby providing a few elements in $\HB \setminus \QN$.

\begin{eg}
	Let $m, n \in \mathbb{N}$ with $n>m+1$. Define $\displaystyle A=\begin{pmatrix}
		0 & m\\
		-1\slash n & 0
	\end{pmatrix}$ which we re-write as $(a_{ij})_{i,j=1}^2$. Evidently, $\|A\|=m$. Let $z_1, z_2, w \in \DC$ with $|w| \leq \sqrt{(1-|z_1|^2)(1-|z_2|^2)}$. Then 
	\begin{align*}
		1-a_{11}z_1-a_{22}z_2+\det(A)z_1z_2=a_{21}w 
		\implies n+mz_1z_2=-w
		\implies |n+mz_1z_2|=|w| \leq 1.
	\end{align*}		
	Since $|z_1|, |z_2| \leq 1$ and $n>m+1$, we have that $|n+mz_1z_2| \geq n-m|z_1z_2| \geq n-m>1$, which is a contradiction. Thus, $1-a_{11}z_1-a_{22}z_2+\det(A)z_1z_2 \ne a_{21}w \ $ for every $z_1, z_2, w \in \overline{\D}$ with $|w| \leq \sqrt{(1-|z_1|^2)(1-|z_2|^2)}$. By Theorem \ref{thm2.4}, $\mu_{\text{hexa}}(A)<1$ and so, $\mu_{\text{hexa}}(A)<1\leq \|A\|$. By Lemma \ref{rem3.6}, $A$ is a unique matrix in $M_2(\C)$ such that $\pi(A)=(-1\slash n, 0, 0, m\slash n)$ and so, $\pi(A) \in \HB \setminus \QN$. \qed 
\end{eg}

We find several characterizations for the points of $\QN$. To do this, we need the following lemma.

\begin{lem}\label{lem4.3}
	For $A=\begin{pmatrix}
		x_1 & w^2\slash a\\
		a & x_2
	\end{pmatrix} \in M_2(\C)$ with $a \ne 0$, let $\be=1-|x_1|^2-|x_2|^2+|x_3|^2$ and $w^2=x_1x_2-x_3$ . If $|x_3|<1$, then the following are equivalent.
	\begin{enumerate}
		\item $\|A\|<1$;
		\smallskip 
		\item $\det(I-A^*A)=|a|^{-2}\left(-|a|^4+|a|^2(1-|x_1|^2-|x_2|^2+|x_3|^2)-|x_1x_2-x_3|^2\right)>0$;
		\smallskip 
		\item $\be-\sqrt{\be^2-4|w|^4}<2|a|^2<\be+\sqrt{\be^2-4|w|^4}$. 	
	\end{enumerate}
\end{lem}

\begin{proof} 
Some routine computations show that 
\begin{equation}\label{eqn4.1}
	\det(I-A^*A)=|a|^{-2}\left(-|a|^4+|a|^2(1-|x_1|^2-|x_2|^2+|x_3|^2)-|x_1x_2-x_3|^2\right),
\end{equation}
and that $
\text{tr}(I-A^*A)-\det(I-A^*A)=(1-|x_3|^2)>0$ since $|x_3|<1$.

\medskip

\noindent $(1) \iff (2)$. Following the proof of Lemma \ref{lem2.3}, we have that $\|A\|<1$ if and only if $I-A^*A>0$. The latter is possible if and only if $\det(I-A^*A), \text{tr}(I-A^*A)>0$. It follows from the above discussion that $\|A\|<1$ if and only if $\det(I-A^*A)>0$. 

\medskip

\noindent $(1) \iff (3).$ Note that \eqref{eqn4.1} can be re-written as $\det(I-A^*A)=-|a|^{-2}\left(|a|^4-\be |a|^2+|w|^4 \right)$.	Hence, $\det(I-A^*A)>0$ if and only if $|a|^4-\be |a|^2+|w|^4<0$. Note that
\begin{equation*}
	\begin{split} 
		\be^2-4|w|^4&=(\be-2|w|^2)(\be+2|w|^2)\\
		&=\left(1-|x_1|^2-|x_2|^2+|x_3|^2-2|x_1x_2-x_3|\right)\left(1-|x_1|^2-|x_2|^2+|x_3|^2+2|x_1x_2-x_3|\right)
	\end{split}
\end{equation*}
which is strictly positive by part-$(8)$ of Theorem \ref{tetrablock}. Thus, we have 
\[
|a|^4-\be |a|^2+|w|^4=\bigg(|a|^2-\frac{\be + \sqrt{\be^2-4|w|^4}}{2}\bigg)\bigg(|a|^2-\frac{\be - \sqrt{\be^2-4|w|^4}}{2}\bigg).
\]
Thus, $\det(I-A^*A)>0$ if and only if either of the following case holds:
\[
\frac{\be - \sqrt{\be^2-4|w|^4}}{2}<|a|^2<\frac{\be + \sqrt{\be^2-4|w|^4}}{2} \ \ \text{or} \	\ \frac{\be + \sqrt{\be^2-4|w|^4}}{2} <|a|^2<\frac{\be - \sqrt{\be^2-4|w|^4}}{2}.
\]
Clearly, the second case is not possible. The desired conclusion now follows.
\end{proof}

We have similar characterizations as in Lemma \ref{lem4.3} for $A \in M_2(\C)$ with $\|A\| \leq 1$. The proof of the next result follows from the fact that $\|A\| \leq 1$ if and only if $I-A^*A \geq 0$. Indeed, one can use the same computations as in Lemma \ref{lem4.3} for this. 

\begin{cor}\label{cor4.4}
	For $A=\begin{pmatrix}
		x_1 & w^2\slash a\\
		a & x_2
	\end{pmatrix} \in M_2(\C)$ with $a \ne 0$, let $\be=1-|x_1|^2-|x_2|^2+|x_3|^2$ and $w^2=x_1x_2-x_3$ . If $|x_3|\leq 1$, then the following are equivalent.
	\begin{enumerate}
		\item $\|A\|\leq 1$;
		\smallskip 
		\item $\det(I-A^*A)=|a|^{-2}\left(-|a|^4+|a|^2(1-|x_1|^2-|x_2|^2+|x_3|^2)-|x_1x_2-x_3|^2\right)\geq 0$;
		\smallskip 
		\item $\be-\sqrt{\be^2-4|w|^4}\leq 2|a|^2\leq \be+\sqrt{\be^2-4|w|^4}$. 	
	\end{enumerate}
\end{cor}

We now provide some useful criteria for elements in $\QN$.

\begin{thm}\label{thm4.4}
	For $(a, x_1, x_2, x_3) \in \C^4$, let $\be=1-|x_1|^2-|x_2|^2+|x_3|^2$ and $w^2=x_1x_2-x_3$. Then the following are equivalent.
	\begin{enumerate}
		\item $(a, x_1, x_2, x_3) \in \QN$;
		
		\vspace{0.15cm}
		
		\item $(x_1, x_2, x_3) \in \E$ and one of the following holds. 
		\[
		(a, x_1, x_2, x_3)=(0, x_1, x_2, x_1x_2) \ \  \text{or} \ \ 	a \ne 0, 	\ \frac{\be-\sqrt{\be^2-4|w|^4}}{2} < |a|^2 < \frac{\be+\sqrt{\be^2-4|w|^4}}{2};
		\]			
		
		\vspace{0.15cm}
		
		\item $(x_1, x_2, x_3) \in \E$ and one of the following holds.
		\[
		(a, x_1, x_2, x_3)=(0, x_1, x_2, x_1x_2) \ \  \text{or} \ \  	a \ne 0, \ 	\det(I-A^*A)>0 \ \ \text{for} \ \  A=\begin{pmatrix}
			x_1 & w^2\slash a\\
			a & x_2
		\end{pmatrix}.
		\]
		
	\end{enumerate}
\end{thm}

\begin{proof}
The equivalence of $(2)$ \& $(3)$ follows directly from Lemma \ref{lem4.3}. 

\medskip 		

\noindent $(1) \implies (3).$ 	Let $(a, x_1, x_2, x_3) \in \QN$. We have by Theorems \ref{thm2.4} \& \ref{thm4.2} that $(a, x_1, x_2, x_3) \in \HB$ and $(x_1, x_2, x_3)\in \E$. By Corollary \ref{cor2.6}, $x_1x_2=x_3$ if $a=0$. Assume that $a \ne 0$. By definition of $\QN$, there is a $2 \times 2$ matrix $A$ such that $\|A\|<1$ and $\pi(A)=(a, x_1, x_2, x_3)$. Since $a \ne 0$, we have by Lemma \ref{rem3.6} that $
A=\begin{pmatrix}
	x_1 & w^2\slash a\\
	a & x_3
\end{pmatrix}$. It follows from Lemma \ref{lem4.3} that $\det(I-A^*A)>0$.

\medskip

\noindent $(3) \implies (1).$ Let $A=\begin{pmatrix}
	x_1 & w^2\slash a\\
	a & x_2	
\end{pmatrix}$ be such that the conditions in $(3)$ hold. Since $(x_1, x_2, x_3) \in \E$, we have by Theorem \ref{tetrablock} that $|x_1|, |x_2|<1$. For $a=0$ and $x_1x_2=x_3$, define $A_0=\begin{pmatrix}
	x_1 & 0 \\
	0 & x_2
\end{pmatrix}$. Then $\|A_0\|=\max\{|x_1|, |x_2|\}<1$ and so, $\pi(A_0)=(0, x_1, x_2, x_1x_2)=(a, x_1, x_2, x_3) \in \QN$. Now assume that $a \ne 0$. It follows from Lemma \ref{lem4.3} that $\|A\|<1$ and thus, $\pi(A)=(a, x_1, x_2, x_3) \in \QN$. 
\end{proof}

Having established various characterizations of $\QN$, we now turn our attention to its closure $\CQN$ and its interior $int(\QN)$. The following results proceed in this direction.

\begin{prop}\label{prop4.6}
	$\CQN=\{\pi(A)\ : \ A \in M_2(\C), \ \|A\| \leq 1\}$.
\end{prop}

\begin{proof}
Let $A \in M_2(\C)$ and let $\|A\| \leq 1$. Define $\epsilon_n=1-1\slash n$ and $A_n=\epsilon_n A$ for $n \in \mathbb{N}$. Then $\|A_n\|<1$ and $\underset{n \to \infty}{\lim} \|A_n-A\|=0$. Thus, $\pi(A_n) \in \QN$ and $\underset{n \to \infty}{\lim}\pi(A_n)=\pi(A)$. Hence, $\pi(A) \in \CQN$ and so, $\{\pi(A) : A \in M_2(\C), \|A\| \leq 1\} \subseteq \CHB$. Let $x=(a, x_1, x_2, x_3) \in \CQN$. Then there is a sequence $\{(a^{(n)}, x_1^{(n)}, x_2^{(n)}, x_3^{(n)})\}$ in $\QN$ that converges to  $(a, x_1, x_2, x_3)$ as $n \to \infty$. We have by Theorem \ref{thm4.4} that each $(x_1^{(n)}, x_2^{(n)}, x_3^{(n)}) \in \E$ and thus, $(x_1, x_2, x_3) \in \overline{\E}$. By definition of $\QN$, there exists $A_n \in M_2(\C)$ such that $\|A_n\|<1$ and
\begin{equation}\label{eqn4.2}
	A_n=\begin{pmatrix}
		x_1^{(n)} & \al^{(n)}\\
		a^{(n)} & x_2^{(n)}
	\end{pmatrix} \quad \text{with} \quad \al^{(n)}a^{(n)}=x_{1}^{(n)}x_{2}^{(n)}-x_{3}^{(n)}
\end{equation}
for all $n \in \mathbb{N}$. Since
$\|A_n\|<1$, the sequence $\{\al^{(n)}\}$ is bounded and therefore must have a convergent subsequence. We now construct a $ 2\times 2$ matrix $A$ such that $\|A\| \leq 1$ and $\pi(A)=(a, x_1, x_2, x_3)$.

\begin{enumerate}
	\item[Case 1.] Let $a=0$. Then $\underset{n \to \infty}{\lim}a^{(n)}=0$ and passing onto subsequence, if necessary, we have that 
	\[
	0=\left(\lim_{n \to \infty}\al^{(n)}\right)\left(\lim_{n \to \infty}a^{(n)}\right)=\lim_{n \to \infty}\al^{(n)}a^{(n)}=\lim_{n \to \infty}\left(x_1^{(n)}x_2^{(n)}-x_3^{(n)}\right)=x_1x_2-x_3.
	\]
For $A=\begin{pmatrix}
		x_1 & 0\\
		0 & x_2\\
	\end{pmatrix}$, it follows that $\|A\|=\max\{|x_1|, |x_2|\} \leq 1$ and $\pi(A)=(0, x_1, x_2, x_3)=x$. 
	
	\item[Case 2.] Let $a \ne 0$. Then there is a convergent subsequence, say $\{a^{(n_k)}\}$, of $\{a^{(n)}\}$ consisting of non-zero scalars. Define $A=\begin{pmatrix}
		x_1 & \al \\
		a & x_2
	\end{pmatrix}$ where $\al=a^{-1}(x_1x_2-x_3)$. Then  
	\[
	\lim_{n_k \to \infty}\al^{(n_k)}=\lim_{n_k \to \infty}\frac{x_1^{(n_k)}x_2^{(n_k)}-x_3^{(n_k)}}{a^{(n_k)}} = \frac{x_1x_2-x_3}{a}=\al.
	\]
	Consequently, $\underset{n_k \to \infty}{\lim} \|A_{n_k}\|=\|A\|$. Thus, $\|A\| \leq 1$ and $\pi(A)=(a, x_1, x_2, x_3)=x$.	
\end{enumerate}	

\smallskip Hence, $\CQN \subseteq \left\{\pi(A) :  A\in M_2(\C), \|A\| \leq 1 \right\}$. The proof is now complete.	
\end{proof}

Using the above result, we obtain the following characterizations for the elements in $\CQN$.

\begin{thm}\label{thm4.6}
	For $(a, x_1, x_2, x_3) \in \C^4$, let $\be=1-|x_1|^2-|x_2|^2+|x_3|^2$ and $w^2=x_1x_2-x_3$. Then the following are equivalent:
	\begin{enumerate}
		\item $(a, x_1, x_2, x_3) \in \CQN$;
		
		\vspace{0.15cm}
		
		\item $(x_1, x_2, x_3) \in \Ebar$ and one of the following holds: 
		\[
		(a, x_1, x_2, x_3)=(0, x_1, x_2, x_1x_2) \ \  \text{or} \ \ 		a \ne 0, 	\ \frac{\be-\sqrt{\be^2-4|w|^4}}{2} \leq |a|^2 \leq \frac{\be+\sqrt{\be^2-4|w|^4}}{2};
		\]
		
		\vspace{0.15cm}
		
		\item $(x_1, x_2, x_3) \in \Ebar$ and one of the following holds:
		\[
		(a, x_1, x_2, x_3)=(0, x_1, x_2, x_1x_2) \ \  \text{or} \ \  	a \ne 0, \ 	\det(I-A^*A) \geq 0 \ \ \text{for} \ \  A=\begin{pmatrix}
			x_1 & w^2\slash a\\
			a & x_2
		\end{pmatrix}.
		\]			
	\end{enumerate}
\end{thm}

\begin{proof}
The equivalence of $(2)$ \& $(3)$ follows from Corollary \ref{cor4.4}. So, it suffices to prove the equivalence of $(1)$ \&  $(3)$. If $(a, x_1, x_2, x_3) \in \overline{\Q}_N$ and $r \in (0,1)$, then by Proposition \ref{prop4.6}, there exists $A \in M_2(\C)$ with $\|A\| \leq 1$ and $\pi(A)=(a, x_1, x_2, x_3)$. Evidently, $\|rA\| <1$ and so, $\pi(rA)=(ra, rx_1, rx_2, r^2x_3) \in \QN$. Consequently, $(ra, rx_1, rx_2, r^2x_3) \in \QN$ for all $r \in (0, 1)$ if and only if $(a, x_1, x_2, x_3) \in \CQN$. Moreover, it is evident from Theorem \ref{thm4.4} that $(ra, rx_1, rx_2, r^2x_3) \in \QN$ if and only if $(rx_1, rx_2, r^2x_3) \in \E$, and one of the following holds:
\[
ra=0, \ (rx_1)(rx_2)=r^2x_3 \quad \text{or} \quad 	ra \ne 0, \ 	\det(I-A_r^*A_r)>0 \ \ \text{for} \ \  A_r=\begin{pmatrix}
	rx_1 & rw^2\slash a\\
	ra & rx_2
\end{pmatrix}.
\]
Note that $A_r=rA$ and so, $\det(I-A^*A)=\underset{r \to 1}{\lim}\det(I-A_r^*A_r) \geq 0$ whenever $\det(I-A_r^*A_r)> 0$ for all $r \in (0, 1)$. Combining these things together, we have that $(a, x_1, x_2, x_3) \in \CQN$ if and only if $(x_1, x_2, x_3) \in \Ebar$, and one of the following holds:
\[
(a, x_1, x_2, x_3)=(0, x_1, x_2, x_1x_2) \ \  \text{or} \ \  a \ne 0, \ 	\det(I-A^*A) \geq 0 \ \ \text{for} \ \  A=\begin{pmatrix}
	x_1 & w^2\slash a\\
	a & x_2
\end{pmatrix}.
\]		
Therefore, we have the equivalence of $(1)$ \& $(3)$, which completes the proof.
\end{proof}

We have by Proposition \ref{prop2.7}
that $\HB$ is not a domain in $\C^4$ and its interior consists of the points $(a, x_1, x_2, x_3) \in \HB$ for which $a\ne 0$. We show that similar results hold for $\QN$.

\begin{thm}\label{QN_not_open}
	$\QN$ is connected but is not an open set in $\C^4$.
\end{thm}

\begin{proof}
$\QN$ is connected as it is the continuous image of the domain $\mathbb{B}_{\|.\|}$ under the map $\pi$. Let $\al_1, \al_2 \in \D$. Then $
A=\begin{pmatrix}
	\al_1 & 0 \\
	0 & \al_2
\end{pmatrix}$ has norm less than $1$ and so, $\pi(A)=(0, \al_1, \al_2, \al_1\al_2) \in \QN$. Let if possible, $(0, \al_1, \al_2, \al_1\al_2)\in int(\QN)$. One can choose $\epsilon>0$ such that $B(0, \epsilon) \times B(\al_1, \epsilon) \times B(\al_2, \epsilon) \times B(\al_1\al_2, \epsilon) \subset \QN$. If this were to happen, then $
(0, \al_1, \al_2, \al_1\al_2+\epsilon \slash 2) \in \QN$, which is a contradiction to  Theorem \ref{thm4.4}. Therefore, $\QN$ is not open and the proof is complete.
\end{proof}

Since $\QN$ is not open, our next step is to determine its interior. Clearly, the proof of Theorem \ref{QN_not_open} shows that a point $(a, x_1, x_2, x_3) \in \QN$ with $a=0$ does not lie in $int(\QN)$. However, we show that all other points are in $int(\QN)$.

\begin{prop}\label{intQN}
	$int(\QN)=\{(a, x_1, x_2, x_3) \in \QN : a \ne 0 \}$. Moreover, $int(\QN)=int(\CQN)$.
\end{prop}	

\begin{proof}
Consider the function $	\phi: \C\setminus \{0\} \times \E \to \mathbb{R}$ defined as
\[
\phi(a, x_1, x_2, x_3)=|a|^{-2}\left(-|a|^4+|a|^2(1-|x_1|^2-|x_2|^2+|x_3|^2)-|x_1x_2-x_3|^2\right).
\]
By Lemma \ref{lem4.3} and Theorem \ref{thm4.4}, $\phi^{-1}(0, \infty)=\{(a, x_1, x_2, x_3) \in \QN  :  a \ne 0\}$. Since $\phi$ is continuous and $\C\setminus\{0\} \times \E$ is open in $\C^4$, we have that $\phi^{-1}(0, \infty)$ is an open set in $\C^4$ contained in $\QN$. Thus, $\phi^{-1}(0, \infty) \subseteq int(\QN) \subseteq \QN$. The proof of Theorem \ref{QN_not_open} shows that no point $(a, x_1, x_2, x_3) \in \QN$ with $a=0$ can be in its interior. Hence, $int(\QN)=\phi^{-1}(0, \infty)=\{(a, x_1, x_2, x_3) \in \QN  :  a \ne 0\}$. Since $\QN \subseteq \CQN \subseteq \C \times \Ebar$, we have that $int(\QN) \subseteq int(\CQN) \subseteq int(\C \times \Ebar)=\C \times \E$. Let $(a, x_1, x_2, x_3) \in int(\CQN)$. Let if possible, $a=0$. It follows from Theorem \ref{thm4.6} that $x_3=x_1x_2$. Following the proof of Theorem \ref{QN_not_open}, we see that there is an open neighbourhood of $(0, x_1, x_2, x_1x_2)$ in $\CQN$ that contains a point of the form $(0, y_1, y_2, y_3)$ with $y_1y_2 \ne y_3$. This is a contradiction to Theorem \ref{thm4.6} and so, $a \ne 0$. Let $x=(x_1, x_2, x_3)$. Then $x \in \E$. let us define $\be=1-|x_1|^2-|x_2|^2+|x_3|^2$ and $w^2=x_1x_2-x_3$. It follows from the proof of Lemma \ref{lem4.3} that $\be^2-4|w|^4>0$. Let 
\[
m=\frac{\be-\sqrt{\be^2-4|w|^4}}{2} \quad \text{and} \quad  M= \frac{\be+\sqrt{\be^2-4|w|^4}}{2}.
\]
By Theorem \ref{thm4.6}, $m \leq |a|^2 \leq M$. Let if possible, $|a|^2=M$. Since $(a, x_1, x_2, x_3) \in int(\CQN)$, there exists $\epsilon>0$ such that $ B(a, \epsilon) \times B(x_1, \epsilon) \times B(x_2, \epsilon) \times B(x_3, \epsilon) \subset \CQN$. Then one can choose $b \in B(a, \epsilon)$ with $|b|>|a|$ and so, $(b, x_1, x_2, x_3) \in \CQN$. Since $b \ne 0$, we have by Theorem \ref{thm4.6} that $|b|^2 \leq M=|a|^2$, a contradiction. Therefore, $|a|^2<M$. Similarly, one can prove that $|a|^2>m$. It follows from Theorem \ref{thm4.4} that $(a, x_1, x_2, x_3) \in \QN$.  Since $a \ne 0$, we have that $(a, x_1, x_2, x_3) \in \phi^{-1}(0, \infty)$ and so, $(a, x_1 ,x_2, x_3) \in int(\QN)$. The proof is now complete.
\end{proof}

It is evident from Proposition \ref{prop2.8} that $\HB$ is $(1,1,1,2)$-quasi-balanced but is neither starlike about the origin nor circled. Here we establish these geometric properties for $\QN$.

\begin{prop}
 The sets $\QN$ and $\CQN$ are $(1,1,1,2)$-quasi-balanced but are neither starlike about $(0,0,0,0)$ nor circled. Also, $\QN$ and $\QN \cap \mathbb{R}^4$ are not convex.		 
\end{prop}

\begin{proof}
The quasi-balanced property follows from the fact that $\pi(zA)=(za, zx_1, zx_2, z^2x_3)$ whenever $\pi(A)=(a, x_1, x_2, x_3)$ for all $A \in M_2(\C)$ and $z \in \C$. To show that $\QN$ and $\CQN$ are not starlike about $(0,0,0,0)$, it suffices to find $x \in \QN$ and $r \in (0, 1)$ such that $rx \notin \CQN$. Let $x_1, x_2 \in \D \setminus \{0\}$. It follows from Theorem \ref{thm4.4} that $x=(0, x_1, x_2, x_1x_2) \in \QN$. By Theorem \ref{thm4.6}, $rx=(0, rx_1, rx_2, rx_1x_2) \notin \CQN$ for any $r \in (0,1)$. Hence, $\QN$ and $\CQN$ are not starlike about $(0,0,0,0)$. The point $x=(0,1,1,1)$ is in $\CQN$ but $ix \notin \CQN$. Thus, neither $\QN$ nor $\CQN$ is circled. We have by  Theorem \ref{thm4.4} that $x=(0, \frac{1}{2}, \frac{1}{2}, \frac{1}{4}), y=(0, \frac{-1}{2}, \frac{-1}{2}, \frac{1}{4})$ are in $\QN \cap \mathbb{R}^4$ but $
(x+y)\slash 2=\left(0, 0, 0, \frac{1}{4}\right)
$ is not in $\QN$. Hence, both $\QN$ and $\QN \cap \mathbb{R}^4$ are non-convex sets.
\end{proof} 

We have proved in Theorem \ref{thm4.2} that $\QN \subsetneq \HB$. The same conclusion holds for their closures as the following result explains.

\begin{thm}\label{thm4.13}
	$\CQN$ is a proper subset of $\CHB$.
\end{thm}

\begin{proof}
Obviously it follows from Theorem \ref{thm4.2} that $\CQN \subset \CHB$. By Lemma \ref{lem4.1}, there is a matrix $A=(a_{ij})_{i, j=1}^2$ and $r>0$ such that $a_{21} \ne 0$ and $\mu_{\text{hexa}}(A)\leq r<\|A\|$. Then $\mu_{\text{hexa}}(r^{-1}A) \leq 1<\|r^{-1}A\|$ and by Theorem \ref{CHB}, $\pi(r^{-1}A) \in \CHB$. Let if possible, $\pi(r^{-1}A) \in \CQN$. It follows from Proposition \ref{prop4.6} that there exists $B=(b_{ij})_{i, j=1}^2 \in M_2(\C)$ such that $\|B\|\leq 1$ and $\pi(B)=\pi(r^{-1}A)$. By Lemma \ref{rem3.6}, $B=r^{-1}A$ which is not possible as $\|B\|\leq 1<\|r^{-1}A\|$. Therefore, $\pi(r^{-1}A) \in \CHB \setminus \CQN$.	
\end{proof} 	

We now explicitly produce some points in $\CHB \setminus \CQN$ and show that it contains an open set.

\begin{prop}\label{prop4.11}
	For $\alpha \in \overline{\D}, \ (\alpha, 0, 0, 1) \in \CHB$. Also, $(\alpha, 0, 0, 1) \in \CQN$ if and only if $\alpha \in \T$.
\end{prop}

\begin{proof}
We have by Theorem \ref{closedtetrablock} that $(0, 0, 1) \in \Ebar$. Let $\alpha \in \overline{\D}$ and let $(a, x_1, x_2, x_3)=(\alpha, 0, 0, 1)$. Then $\be=1-|x_1|^2-|x_2|^2+|x_3|^2=2, w^2=x_1x_2-x_3=-1$ and so, $\be^2-4|w|^4=0$. Since $x_1x_2 \ne x_3$, it follows from Theorem \ref{thm4.6} that $(a, x_1, x_2, x_3) \in \CQN$ if and only if $|a|^2=1$. Hence, $(\alpha, 0, 0, 1) \in \CQN$ if and only if $\alpha \in \T$. In particular, $(1, 0, 0, 1) \in \CQN$ and so, $(1, 0, 0, 1) \in \CHB$ by Theorem \ref{thm4.13}. It follows from Theorem \ref{CHB} that $\underset{z_1, z_2 \in \D}{\sup}|\Pz(1, 0, 0, 1)| \leq 1$. Let $\alpha \in \overline{\D}$. Then
\[
\sup_{z_1, z_2 \in \D}|\Pz(\al, 0, 0, 1)|= |\al|\sup_{z_1, z_2 \in \D}|\Pz(1, 0, 0, 1)| \leq 1.
\]
Thus, Theorem \ref{CHB} shows that $(\al, 0, 0, 1) \in \CHB$. This completes the proof.
\end{proof}

\begin{prop}\label{prop_open}
	$\HB \setminus \QN=int(\HB) \setminus \QN$ and it contains a non-empty open set in $\C^4$.
\end{prop}

\begin{proof} Evidently $int(\HB) \setminus \QN \subseteq \HB \setminus \QN$. Let $q_0=(a, x_1, x_2, x_3) \in \HB \setminus \QN$. For $a=0$, it follows from Corollary \ref{cor2.6} that $x_1x_2=x_3$ and $(x_1, x_2, x_3)\in \E$. By Theorem \ref{thm4.4}, we have $q_0 \in \QN$, which is a contradiction. Thus, $a \ne 0$ and so, by Proposition \ref{prop2.9}, $q_0\in int(\HB)$. Consequently, $\HB \setminus \QN \subseteq int(\HB) \setminus \QN$. We now prove the existence of a non-empty open set contained in $\HB \setminus \QN$.   By Lemma \ref{lem4.1}, one can choose $A=(a_{ij})_{i, j=1}^2 \in M_2(\C)$ with $\mu_{\text{hexa}}(A)<1 <\|A\|$ and $a_{21} \ne 0$. Let $q=\pi(A)$. Then $q \in \HB \setminus \QN$ and the first component of $q$ is non-zero. By Lemma \ref{rem3.6}, $A$ is the unique matrix in $M_2(\C)$ with $q=\pi(A)$. Since $\HB \setminus \QN= int(\HB) \setminus \QN$, it follows that $q \in int(\HB) \setminus \QN$.  We show that there is an open set in $\C^4$ that contains $q$ but does not intersect $\QN$. If not, then $q \in \CQN$. By Proposition \ref{prop4.6}, there exists $B \in M_2(\C)$ with $\|B\| \leq 1$ and $\pi(B)=q$. By uniqueness, $B=A$ which is not possible as $\|B\| \leq 1 < \|A\|$. Hence, there is an open set $U$ in $\C^4$ such that $q \in U$ and $U \cap \QN =\emptyset$. Clearly, $q \in U \cap int(\HB)$. Also,  $U \cap int(\HB)$ is an open set in $\C^4$, and it is contained in $\HB \setminus \QN$. The proof is now complete. 
\end{proof}

\section{The hexablock}\label{hexa}

\noindent As discussed in the previous sections, the $\mu$-hexablock $\HB$ and the normed hexablock $\QN$ have similar origin as that of $\Gg, \E$ and $\Pe$ in the context of the $\mu$-synthesis problem. However, it turns out that neither $\HB$ nor $\QN$ are domains in $\C^4$. We also mention that the map $(s, p) \mapsto (0, s, p)$ for $(s, p) \in \Gg$ is a natural embedding of $\Gg$ into $\Pe$. The authors of \cite{AglerIV} capitalized this embedding to establish geometric and function theoretic properties of $\Pe$ using the theory of $\Gg$. In fact, they also obtained an alternative description of the pentablock (see Theorem 5.2 in \cite{AglerIV}) given by
\[
\Pe=\left\{(a, s, p) \in \C \times \Gg : \sup_{z \in \D}|\psi_z(a, s, p)|<1 \right\},
\]
where the maps $\psi_z(a, s, p)$ on $\C \times \Gg$ were defined in \eqref{eqn_303}. As discussed in the beginning of Section \ref{psi_map}, the maps $\Pz$ for $(z_1, z_2) \in \D^2$ are generalization of the maps $\psi_z$ for $z \in \D$ in the sense that $\psi_{z, z}(a, s\slash 2, s\slash 2, p)$ coincides with $\psi_z(a,s,p)$. Motivated by the aforementioned description of $\Pe$, we define the set \textit{hexablock} $\Q$ in the following way:
\[
\Q:=\left\{(a, x_1, x_2, x_3) \in \C \times \E :  \sup_{z_1, z_2 \in \D}|\Pz(a, x_1, x_2, x_3)|<1 \right\}.
\]
In this section, our central object of study is the hexablock. First we prove that $\Q$ is a domain in $\C^4$. Clearly, $\{0\} \times \E \subseteq \Q$ and so, we have an embedding of $\E$ inside $\Q$ via the natural map $(x_1, x_2, x_3) \mapsto (0, x_1, x_2, x_3)$. Another motivation behind the hexablock comes from a simple fact that the interior of a set is always an open set. So, even though $\HB$ and $\QN$ are not open sets, the sets $int(\CHB), int(\HB), int(\QN)$ and $int(\CQN)$ are always open in $\C^4$. In fact, we shall prove that each of these sets is connected also and hence is a domain in $\C^4$, and that
\[
int(\CQN)=int(\QN) \subsetneq int(\HB) \subsetneq int(\CHB).
\]
We have by Propositions \ref{prop2.9} \& \ref{intQN} that the domains $int(\HB), int(\QN)$ and $int(\CQN)$ do not intersect with $\{0\} \times \E$, making it impossible to embed $\{0\} \times \E$ canonically into these three domains. However, we shall see later that $\{0\} \times \E$ is contained in $int(\CHB)$. Since we wish to have a domain in $\C^4$ closely related to $\HB$ and $\QN$, with $\{0\} \times \E$ inside it, $int(\CHB)$ seems to be another domain of our interest other than the hexablock. Interestingly, it turns out that \[
\Q= int(\CHB) \quad \text{and} \quad \QN \subsetneq \HB \subsetneq \Q. 
\] 
This establishes a connection of the hexablock $\Q$ with the $\mu$-hexablock $\HB$ and the norm hexablock $\QN$. Furthermore, we prove various complex geometric properties of $\Q$ and to begin with, let us have the following characterizations of the points in $\Q$.

\begin{thm}\label{thm5.1}
	Let $(a, x_1, x_2, x_3) \in \C^4$. Then the following are equivalent:
	\begin{enumerate}
		\item $(a, x_1, x_2, x_3) \in \Q$;
		\item $(x_1, x_2, x_3) \in \E$ and $|a|K_*(x_1, x_2, x_3)<1$;
		\item $(x_1, x_2, x_3) \in \E$ and $|\Pz(a, x_1, x_2, x_3)|<1$ for all $z_1, z_2 \in \D$;
		\item either $(a, x_1, x_2, x_3) \in \HB$ or $(a, x_1, x_2, x_3) \in \{0\} \times \E$,
	\end{enumerate}
	where $K_*$ is the map as in \eqref{eqn_K_*}.
	Consequently, $\Q=\HB \cup (\{0\} \times \E)=int(\HB) \cup (\{0\} \times \E)$ and $\QN \subset \HB \subset \Q$.
\end{thm}	

\begin{proof}
The equivalences of $(1)$, $(2$) \& $(3)$ follows from Corollary $\ref{cor2.3}$. 

\medskip 

\noindent $(1) \implies (4)$. Let $(a, x_1, x_2, x_3) \in \Q$. If $a=0$, then it is trivial that $(a, x_1, x_2, x_3) \in \{0\} \times \E$. Let $a \ne 0$. Then we have by Lemma \ref{rem3.6} that there exists $A \in M_2(\C)$ such that $\pi(A)=(a, x_1, x_2, x_3)$. It follows from Theorem \ref{thm2.4} that $\mu_{\text{hexa}}(A)<1$ and so, $(a, x_1, x_2, x_3) \in \HB$.  

\medskip 

\noindent $(4) \implies (1)$. Evidently, $\{0\} \times \E \subset \Q$ as $\underset{z_1, z_2 \in \D}{\sup}|\Pz(0, x_1, x_2, x_3)|=0$ for  $(x_1, x_2, x_3) \in \E$. Let $(a, x_1, x_2, x_3) \in \HB$. By Corollary \ref{cor2.6}, $(x_1, x_2, x_3) \in \E$ and $\underset{z_1, z_2 \in \D}{\sup}|\Pz(a, x_1, x_2, x_3)|<1$. Thus, $(a, x_1, x_2, x_3) \in \Q$ and so, $(1)$ holds. 

\medskip 

\noindent Hence, $\Q=\HB \cup (\{0\} \times \E)$ and by Theorem \ref{thm4.2}, we have that $\QN \subset \HB \subset \Q$. It follows from Proposition \ref{prop2.9} that $int(\HB)=\HB \setminus (\{0\} \times \E)$ and so, $\HB \cup (\{0\} \times \E)=int(\HB) \cup (\{0\} \times \E)$.
\end{proof} 

The next result provides one of the first domains in $\C^4$ related to $\HB$ and $\QN$.

\begin{prop}\label{prop5.2}
	The hexablock $\Q$ is a domain in $\C^4$.
\end{prop}	

\begin{proof}
We have by Proposition \ref{prop2.7} that  $\HB$ is connected. Moreover, $\{0\} \times \E$ is also connected in $\C^4$. It follows from Theorem \ref{thm5.1} that $\Q$ is the union of two connected sets $\HB$ and $\{0\} \times \E$ whose intersection is non-empty. Indeed, $\HB \cap (\{0\} \times \E)=\{(0, x_1, x_2, x_3) \in \C \times \E :  x_1x_2=x_3\}$ which follows from Corollary \ref{cor2.6}. Consequently, $\Q$ is connected. To show that $\Q$ is open in $\C^4$, we consider the following map:
\[
f: \C \times \E \to \mathbb{R}, \quad f(a, x_1, x_2, x_3)=|a|K_*(x_1, x_2, x_3).
\] 
From the discussion above Corollary \ref{cor2.3}, it follows that $K_*$ is a continuous map on $\E$. Therefore, $f$ is a continuous function on $\C \times \E$. It follows from Theorem \ref{thm5.1} that $f^{-1}(-\infty, 1)=\{(a, x_1, x_2, x_3) \in \C \times \E \ : \ |a|K_*(x_1, x_2, x_3)<1 \}=\Q$. Thus, $\Q$ is open in $\C \times \E$. Since $\C \times \E$ is open in $\C^4$, we have that $\Q$ is an open subset of $\C^4$ and hence, $\Q$ is a domain in $\C^4$.
\end{proof} 

The following result show that $\Q$ and $\HB$ are closely related.

\begin{thm}\label{thm5.3}
	$int(\CHB)=\Q$ and 
	\[
	\CHB=\overline{\Q}=\left\{(a, x_1, x_2, x_3) \in \C \times \overline{\E} \ : \ \bigg|\frac{a\sqrt{(1-|z_1|^2)(1-|z_2|^2)}}{1-x_1z_1-x_2z_2+x_3z_1z_2}\bigg| \leq 1 \ \text{for every} \ z_1, z_2 \in \D \right\}.
	\]
	Moreover, $\CQN \subseteq \CQ \subseteq \overline{\D} \times \overline{\D} \times \overline{\D} \times \overline{\D}$. 
\end{thm}	

\begin{proof}
We have by Theorems \ref{CHB}, \ref{thm4.2} \&  \ref{thm5.1} that $\QN \subseteq \HB \subset \Q \subset \CHB$. As a  consequence, $\CQN \subseteq \CHB=\CQ$ and $int(\overline{\Q})=int(\CHB)$. Therefore, $\Q \subseteq int(\overline{\Q})=int(\CHB)$.  Since $\E$ is a bounded domain, we have that $int(\overline{\E})=\E$. Thus, $\CHB \subseteq \C \times \overline{\E}$ and so, $int(\CHB) \subseteq \C \times int(\overline{\E})=\C \times \E$. Let $(a, x_1, x_2, x_3) \in int(\CHB)$. Then $(x_1, x_2, x_3) \in \E$ and it follows from the discussion above Corollary \ref{cor2.3} that $K_*(x_1, x_2, x_3)$ is a positive real number. It remains to show that $|a|K_*(x_1, x_2, x_3)<1$. Since $(a, x_1, x_2 , x_3) \in \CHB$, it follows from Theorem \ref{CHB} that $|a| K_*(x_1, x_2, x_3) \leq 1$. Let if possible, $|a|=K_*(x_1, x_2, x_3)^{-1}$. Choose an $\epsilon>0$ such that $B(a, \epsilon) \times B(x_1, \epsilon) \times B(x_2, \epsilon) \times B(x_3, \epsilon) \subset \CHB$. Then there exists $(b, x_1, x_2, x_3) \in \CHB$ with $|b|>|a|=K_*(x_1, x_2, x_3)^{-1}$. This is not possible because $|b| \leq K_*(x_1, x_2, x_3)^{-1}$ by Theorem \ref{CHB}. Hence, $|a|<K_*(x_1, x_2, x_3)^{-1}$ and so, we have by Theorem \ref{thm5.1} that $(a, x_1, x_2, x_3) \in \Q$. Therefore, $\Q=int(\CHB)$. For any $(a, x_1, x_2, x_3) \in \CQ$, it is clear from Theorem \ref{CHB} that $(x_1, x_2, x_3) \in \Ebar$. We have by Theorem \ref{closedtetrablock} that $\Ebar \subset \overline{\D}^3$ and so, $(x_1, x_2, x_3) \in \overline{\D}^3$. Moreover, 
\[
\bigg|\frac{a\sqrt{(1-|z_1|^2)(1-|z_2|^2)}}{1-x_1z_1-x_2z_2+x_3z_1z_2}\bigg| \leq 1 \quad \text{for all $z_1, z_2 \in \D$}.
\]
Choosing $z_1=z_2=0$ implies that $|a| \leq 1$ and so, $\overline{\Q} \subseteq \overline{\D}^4$. The proof is now complete.
\end{proof} 	

An immediate consequence to the above theorem is the following result. 

\begin{cor}\label{cor_CQ_r_Q}
	Let $(a, x_1, x_2, x_3) \in \CQ$. Then $(ra, rx_1, rx_2, r^2x_3) \in \Q$ for every $r \in (0,1)$. 
\end{cor}

\begin{proof} 
Recall that $\pi_{\E}(A)=(a_{11}, a_{22}, \det(A))$ for $A=(a_{ij})_{i, j=1}^2 \in M_2(\C)$. Let $r \in (0,1)$. Since $(a, x_1, x_2, x_3) \in \CQ$, we have by Theorem \ref{thm5.3} that $(x_1, x_2, x_3) \in \Ebar$. By part-$(5)$ of Theorem \ref{closedtetrablock}, there exists $A \in M_2(\C)$ with $\|A\| \leq 1$ such that $\pi_{\E}(A)=(x_1, x_2, x_3)$. Clearly, $\|rA\|\leq r<1$ and $\pi_{\E}(rA)=(rx_1, rx_2, r^2x_3)$. It follows from part-$(6)$ of Theorem \ref{tetrablock} that $(rx_1, rx_2, r^2x_3) \in \E$. For $z_1, z_2 \in \D$, we have that
\begin{equation*}
	\begin{split}
		\bigg|\frac{ra\sqrt{(1-|z_1|^2)(1-|z_2|^2)}}{1-rx_1z_1-rx_2z_2+r^2x_3z_1z_2}\bigg| & 
		\leq 		\bigg|\frac{ra\sqrt{(1-|rz_1|^2)(1-|rz_2|^2)}}{1-rx_1z_1-rx_2z_2+r^2x_3z_1z_2}\bigg|\\
		& <	\bigg|\frac{a\sqrt{(1-|rz_1|^2)(1-|rz_2|^2)}}{1-x_1(rz_1)-x_2(rz_2)+x_3(rz_1)(rz_2)}\bigg|\\
		& \leq 1,
	\end{split}
\end{equation*}
where the last inequality follows from Theorem \ref{thm5.3} as $(a, x_1, x_2, x_3) \in \CQ$. It follows from Part (3) of Theorem \ref{thm5.1} that $(ra, rx_1, rx_2, r^2x_3) \in \Q$. The proof is now complete. 
\end{proof}

Theorem \ref{thm5.3} tells us that $\CQ \subseteq \DC^4$. We further refine this bound on $\CQ$ and show that $\CQ \subseteq \overline{\mathbb{B}}_2 \times \DC \times \DC$, where $\mathbb B_2$ is the biball, i.e., $\mathbb{B}_2=\{(z_1, z_2) \in \C^2: |z_1|^2+|z_2|^2<1\}$. 

\begin{lem}\label{lem6.5}
	Let $(a, x_1, x_2, x_3) \in \CQ$. Then $|a|^2+|x_1|^2 \leq 1$ and $|a|^2+|x_2|^2 \leq 1$. Moreover, if $(a, x_1, x_2, x_3) \in \Q$, then $|a|^2+|x_1|^2<1$.
\end{lem}

\begin{proof}
It follows from Theorem \ref{thm5.3} that $(x_1, x_2, x_3) \in \overline{\E}$ and
\begin{equation}\label{eqn6.1}
	\bigg|\frac{a\sqrt{(1-|z_1|^2)(1-|z_2|^2)}}{1-x_1z_1-x_2z_2+x_3z_1z_2}\bigg| \leq 1\quad \text{for every} \quad z_1, z_2 \in \D.
\end{equation}
Moreover, $|x_1|, |x_2| \leq 1$ as $\overline{\E} \subseteq  \overline{\D}^3$. If $|x_1|<1$, then putting $(z_1, z_2)=(\overline{x}_1, 0)$ in  \eqref{eqn6.1}, we have that $|a| \leq \sqrt{1-|x_1|^2}$. If $|x_1| = 1$, then putting $z_1=r\overline{x}_1 \ (0<r<1)$ and $z_2=0$ in \eqref{eqn6.1} lead to
\[
\frac{|a|\sqrt{1-r^2}}{1-r} \leq 1 \quad  \text{and so,} \quad  |a| \leq \frac{\sqrt{1-r}}{\sqrt{1+r}}.
\]
Taking limit as $r \to 1$, we have $|a| =0=1-|x_1|^2$. In either case, $|a|^2+|x_1|^2 \leq 1$. The rest of the conclusion can be proved similarly.
\end{proof}

For points $(x_1, x_2, x_3) \in b\E$, we have a simple criteria for $a \in \C$ so that $(a, x_1, x_2, x_3) \in \overline{\Q}$.

\begin{cor}\label{cor_605}
	A point $(a, x_1, x_2, x_3)$ in $\C \times b\E$ belongs to $\CQ$ if and only if $|a|^2+|x_1|^2 \leq 1$. 
\end{cor}

\begin{proof}
The forward part follows from Lemma \ref{lem6.5}. Let $(x_1, x_2, x_3) \in b\E$ and $|a|^2+|x_1|^2 \leq 1$. By part-$(1)$ of Theorem \ref{thm6.1}, $x_1=\overline{x}_2x_3, |x_1| \leq 1$ and $|x_3|=1$. If $a=0$, then it follows from Theorem \ref{thm5.3} that $(a, x_1, x_2, x_3) \in \CQ$. Assume that $a \ne 0$. Then $|x_1| \leq \sqrt{1-|a|^2}<1$ and so, $|x_1|=|x_2|<1$. It follows from Theorem \ref{thm_sup_bE} that
\[
\sup_{z_1, z_2 \in \D} \bigg|\frac{a\sqrt{(1-|z_1|^2)(1-|z_2|^2)}}{1-x_1z_1-x_2z_2+x_3z_1z_2}\bigg|=\frac{|a|}{\sqrt{1-|x_1|^2}} \leq 1.
\]
By Theorem \ref{thm5.3}, $(a, x_1, x_2, x_3) \in \CQ$ which completes the proof.
\end{proof}

We now prove some useful geometric properties of the hexablock. An immediate one is the following whose proof follows directly from Theorems \ref{thm5.3} \& \ref{thm2.9}.

\begin{thm}\label{thm_Qpolyc}
	$\overline{\Q}$ is a polynomially convex set. 
\end{thm}	

\begin{defn} 
	A domain $\Omega$ in $\C^n$ is said to be \textit{polynomially convex} if for each compact subset $K$ of $\Omega$, the polynomial hull $\widehat{K}$ of $K$ is contained in $\Omega$.
\end{defn}

\begin{prop}\label{Qpolyconvex}
	$\Q$ is a polynomially convex domain.
\end{prop}

\begin{proof}
For $r \in (0,1)$, define the compact set
\[
\Q^{(r)}=\left\{(a, x_1, x_2, x_3) \in \C^4 : (r^{-1}a, r^{-1}x_1, r^{-1}x_2, r^{-2}x_3) \in \overline{\Q} \right\}.
\]
We show that $\Q^{(r)} \subset \Q$. Let $(a, x_1, x_2, x_3) \in \Q^{(r)}$. Then  $(r^{-1}a, r^{-1}x_1, r^{-1}x_2, r^{-2}x_3) \in \overline{\Q}$ and so, $(r^{-1}x_1, r^{-1}x_2, r^{-2}x_3) \in \overline{\E}$ by Theorem \ref{thm5.3}. It follows from Lemma 4.1.1 in \cite{AbouhajarI} that $(x_1, x_2, x_3) \in \E$ and so, $(a, x_1, x_2, x_3) \in \Q$ if $a=0$. Assume that $a \ne 0$. We have by Lemma \ref{rem3.6} that there is a unique $A_r \in M_2(\C)$ such that $\pi(A_r)=(r^{-1}a, r^{-1}x_1, r^{-1}x_2, r^{-2}x_3)$ and so, $\pi(A_r) \in \CQ$. It follows from Proposition \ref{prop2.5} and Theorem \ref{thm5.3} that $\mu_{\text{hexa}}(A_r) \leq 1$.  Let $A=rA_r$. Then $\mu_{\text{hexa}}(A) \leq r<1$ and thus, $\pi(A)=(a, x_1, x_2, x_3) \in \HB$. We have by Theorem \ref{thm5.1} that $(a, x_1, x_2, x_3) \in \Q$. Therefore, $\Q^{(r)} \subset \Q$ for every $r \in (0,1)$.

\medskip

Let $(a, x_1, x_2, x_3) \in \Q$. Then $(x_1, x_2, x_3) \in \E$. We show that $(a, x_1, x_2, x_3) \in \Q^{(r)}$ for some $r \in (0,1)$. If $a=0$, then Corollary 4.1.4 in \cite{AbouhajarI} ensures the existence of $r \in (0,1)$ such that $(r^{-1}x_1, r^{-1}x_2, r^{-2}x_3) \in \overline{\E}$. It follows from Theorem \ref{thm5.3} that $(0, r^{-1}x_1, r^{-1}x_2, r^{-2}x_3) \in \CQ$ and so, $(0, x_1, x_2, x_3) \in \Q^{(r)}$. Assume that $a \ne 0$. Then $(a, x_1, x_2, x_3) \in \HB$ by Theorem \ref{thm5.1}. By definition, there exists $A \in M_2(\C)$ such that $\pi(A)=(a, x_1, x_2, x_3)$ and $\mu_{\text{hexa}}(A) <1$. Choose $r$ such that $\mu_{\text{hexa}}(A) < r < 1$. For $B=r^{-1}A$, we have that $\mu_{\text{hexa}}(B) \leq 1$. By Theorems \ref{CHB} \& \ref{thm5.3}, $\pi(B)=(r^{-1}a, r^{-1}x_1, r^{-1}x_2, r^{-2}x_3) \in \CQ$ and thus, $(a, x_1, x_2, x_3) \in \Q^{(r)}$. Putting everything together, we have that
\[
\bigcup_{0<r<1} \Q^{(r)}=\Q.
\]
Evidently, $\Q^{(r)}=\left\{(ra, rx_1, rx_2, r^2x_3) \in \C^4 : (a, x_1, x_2, x_3) \in \overline{\Q} \right\}$ for $r \in (0,1)$ and so, the map $f_r: \overline{\Q} \to \Q^{(r)}$ given by $f_r(a, x_1, x_2, x_3)=\left(ra, rx_1, rx_2, r^2x_3\right)$ is a homeomorphism. Thus, $\Q^{(r)}$ is compact for all $r \in (0, 1)$. Let $\al \notin \Q^{(r)}$. Then $f_r^{-1}(\al) \notin \CQ$. By Theorem \ref{thm_Qpolyc}, $\CQ$ is polynomially convex and so, there is a polynomial $p$ in $4$-variables such that $\|p\|_{\infty, \CQ} \leq 1 <|p(f_r^{-1}(\al))|$. Then $g=p \circ f_r^{-1}$ is a polynomial that satisfies $\|g\|_{\infty, \Q^{(r)}} \leq 1<|g(\al)|$. Consequently, $\Q^{(r)}$ is a polynomially convex set for $r \in (0, 1)$.

\smallskip 

Let $K$ be a compact set of $\Q$. One can choose $r$ sufficiently close to $1$ such that $
K \subseteq \Q^{(r)}$. Assume on the contrary that no such $r$ exists. Then for every $r \in (0, 1)$, there exists $\alpha_r=(a_{1r}, x_{1r}, x_{2r}, x_{3r}) \in K$ such that $\alpha_r \notin \Q^{(r)}$. Let $\beta_r=(r^{-1}a_{1r}, r^{-1}x_{1r}, r^{-1}x_{2r}, r^{-2}x_{3r})$. It follows from the definition of $\Q^{(r)}$ that $\beta_r \in  \C^4 \setminus \CQ$. Note that $\{\al_r : 0<r<1\}$ is a bounded subset of a compact set $K$ and so, it has a convergent subsequence, say $\{\al_{r_n}\}$. Let $\underset{r_n \to 1}{\lim}\al_{r_n}=\al=(a, x_1, x_2, x_3)$. Consequently, $\al \in K$ and $\underset{r_n \to 1}{\lim}\beta_{r_n}=\al$. Since each $\beta_{r_n} \in \C^4 \setminus \CQ$, we must have $\al \in \overline{\C^4 \setminus \CQ}=\C^4 \setminus int(\CQ)=\C^4 \setminus \Q$, where the last equality follows from Theorem \ref{thm5.3}. This gives a contradiction as $\al \in K \subset \Q$. Hence, one can find $\Q^{(r)}$ that contains $K$ for some $r \in (0, 1)$. Since $\Q^{(r)}$ is polynomially convex, we have that $\widehat{K} \subset \widehat{\Q}^{(r)} =\Q^{(r)} \subset \Q$. Thus, $\Q$ is polynomially convex and the proof is complete.
\end{proof}

\begin{prop}\label{prop_qb}
	$\Q$ is $(1,1,1,2)$-quasi-balanced.
\end{prop}

\begin{proof}It follows from part-$(6)$ of Theorem \ref{thm5.1} that $\Q=\HB \cup (\{0\} \times \E)$. We have by Proposition \ref{prop2.8} and Theorem  \ref{tetrablock} that both $\HB$ and $\{0\} \times \E$ are $(1,1,1,2)$-quasi-balanced and so, the desired conclusion follows. 
\end{proof}

We now discuss about linear convexity of the hexablock. Note that a domain $\Omega \subset \C^n$ is said to be \textit{linearly convex} if for any $z \notin \Omega$, there is a complex hyperplane passing through $z$ that does not intersect $\Omega$. It was established separately in \cite{NikolovII}, \cite{Zwonek} and \cite{KosinskiII} the linear convexity of the domains $\Gg, \E$ and $\Pe$, respectively.

\begin{prop}\label{prop_Qlc}
	$\Q$ is linearly convex.
\end{prop}	

\begin{proof}
Let $(a, x_1, x_2, x_3) \notin \Q$. We must find a complex hyperplane passing through $(a, x_1, x_2, x_3)$ that does not intersect $\Q$. By definition, either $(x_1, x_2, x_3) \notin \E$ or there is some $(z_1, z_2) \in \D^2$ such that $|\Pz(a, x_1, x_2, x_3)|\geq 1$. Let $(x_1, x_2, x_3) \notin \E$. Since $\E$ is linearly convex, there is a complex plane $P$ passing through $(x_1, x_2, x_3)$ and omitting $\E$. Then $\C \times  P$ is a complex hyperplane passing through $(a, x_1, x_2, x_3)$ without any intersection with $\Q$. Suppose that $(x_1, x_2, x_3) \in \E$. We have by Theorem \ref{thm5.1} that there exist $(z_1, z_2) \in \D^2$ and $\omega \in \C \setminus \D$ such that $\Pz(a, x_1, x_2, x_3)=\omega$. Note that $|\Pz(b, y_1, y_2, y_3)|<1$ for any $(b, y_1, y_2, y_3) \in \Q$. Thus, $\{(b, y_1, y_2, y_3) \in \C^4 : \Pz(b, y_1, y_2, y_3) =\omega \}$
is a complex hyperplane that passes through $(a, x_1, x_2, x_3)$ without intersecting $\Q$.
\end{proof}

\subsection*{Hexablock as a Hartogs domain} Following the description of the hexablock $\Q$, we learn that it is closely related to the tetrablock $\E$. We show that $\Q$ is a Hartogs domain over $\E$.

\begin{defn}For a domain $\Omega \subseteq \C^n$ and a continuous function $\ell$ on $\Omega$, the domain
	\[
	\Omega_0=\{(a, z) \in \C^d \times \Omega : \|a\|^2 < \ell(z)\}
	\]
	is called a \textit{Hartogs domain} in $\C^{n+d}$ over $\Omega$.
\end{defn}

It was proved in \cite{Guicong} that $\Pe$ is a Hartogs domain over $\Gg$. An analogous result holds for $\Q$, as shown below.

\begin{thm}\label{thm_Hartogs}
	$\Q$ is a Hartogs domain in $\C^4$ over $\E$.
\end{thm}

\begin{proof} 
	Let $(x_1, x_2, x_3) \in \E$. By Proposition  \ref{prop2.1}, there is a unique point $(z_1(x), z_2(x))\in \D^2$ such that 
\[
(z_1, z_2) \mapsto \bigg| \frac{1-x_1z_1-x_2z_2+x_3z_1z_2}{\sqrt{(1-|z_1|^2)(1-|z_2|^2)}} \bigg|^{-1}
\]
attains its supremum over $\D^2$. It follows from  Theorem \ref{thm5.1} that 
\begin{equation}\label{eqn_Hartogs}
	\Q=\left\{(a, x_1, x_2, x_3) \in \C \times \E : |a|^2< e^{-u(x_1, x_2, x_3)}\right\},
\end{equation}
where 	
\begin{equation*}
	u(x_1, x_2, x_3)=-2\log K_*(x_1, x_2, x_3)^{-1}=2\log\bigg| \frac{1-x_1z_1(x)-x_2z_2(x)+x_3z_1(x)z_2(x)}{\sqrt{(1-|z_1(x)|^2)(1-|z_2(x)|^2)}} \bigg|
\end{equation*}
for $(x_1, x_2, x_3) \in \E$. Let 
\[
g(x)+ih(x)=\frac{1-x_1z_1(x)-x_2z_2(x)+x_3z_1(x)z_2(x)}{\sqrt{(1-|z_1(x)|^2)(1-|z_2(x)|^2)}} \qquad (x=(x_1, x_2, x_3) \in \E)
\]
and so, $u(x)=\log(g^2(x)+h^2(x))$ is an infinitely many times real-differentiable function. 
\end{proof} 

\subsection*{Subdomains of the hexablock} We have proved in Proposition \ref{prop5.2} that $\Q$ is a domain in $\C^4$. Moreover, by Theorems \ref{thm5.1} \& \ref{thm5.3}, 
$
\Q=int(\CHB)=\HB \cup (\{0\} \times \E)=int(\HB) \cup (\{0\} \times \E).
$
So, $int(\CHB)$ and $\HB \cup (\{0\} \times \E)$ are domains in $\C^4$. In the same spirit, we study the sets 
\[
int(\CQN), \quad int(\QN), \quad int(\HB) \quad \text{and} \quad \QN \cup (\{0\} \times \E),
\]
and obtain out of them a few domains in $\C^4$ that are contained in $\mathbb{H}$. To begin with, we have the following proposition. 

\begin{prop}
	$int(\CQN)=int(\QN) \subsetneq int(\HB) \subsetneq int(\CHB)=\Q$.
\end{prop}

\begin{proof}  It follows directly from  Theorem \ref{thm4.2}, Proposition \ref{intQN}, and Theorem \ref{thm5.3} that $int(\CQN)= int(\QN) \subseteq int(\HB) \subseteq int(\CHB)=\Q$. By Theorems \ref{thm5.1} \& \ref{thm5.3}, 
$
\Q=int(\CHB)=int(\HB) \cup (\{0\} \times \E).
$ So, $int(\HB) \ne int(\CHB)$. Let if possible, $int(\QN)=int(\HB)$. Now, it follows from Theorem \ref{thm4.4} and Proposition \ref{intQN} that 
\begin{align*}
	\QN =int(\QN) \cup \{(0, x_1, x_2, x_3) \in \E : x_1x_2=x_3\}=int(\HB) \cup \{(0, x_1, x_2, x_3) \in \E : x_1x_2=x_3\}.
\end{align*} 
Also, by Corollary \ref{cor2.6} and Proposition \ref{prop2.9}, $\HB=int(\HB) \cup \{(0, x_1, x_2, x_3) \in \E : x_1x_2=x_3\}$ and so, $\HB=\QN$. This is a contradiction to Theorem \ref{thm4.2}. The proof is now complete.
\end{proof}

Our next few results show that $int(\QN)$ and $int(\HB)$ are connected and hence domains in $\C^4$.

\begin{prop}
	$int(\HB)$ is a domain in $\C^4$. Moreover, 
	\[
	int(\HB)=\Q \setminus (\{0\} \times \E)=\left\{(a, x_1, x_2, x_3) \in \Q : a \ne 0\right\}.
	\]
\end{prop}

\begin{proof} By Theorems \ref{thm5.1} \& \ref{thm5.3}, $\Q=\HB \cup (\{0\} \times \E)=int(\HB) \cup (\{0\} \times \E)$. It now follows from Proposition \ref{prop2.9} that
$int(\HB)=\Q \setminus (\{0\} \times \E)=\left\{(a, x_1, x_2, x_3) \in \Q : a \ne 0\right\}$. Since $\Q$ is a domain, we have by Proposition 1.1 in \cite{FG} that $int(\HB)$ is a domain in $\C^4$. 
\end{proof}

\begin{prop}\label{prop_cintCQN}
	$int(\QN)$ is a domain in $\C^4$ and $\overline{int(\QN)}=\CQN$.
\end{prop}

\begin{proof} By Proposition \ref{intQN}, $int(\QN)=int(\CQN)$. Evidently, $int(\QN)$ is open in $\C^4$. Define 
\[
B_*=\{A=(a_{ij})_{i,j=1}^2 \in \mathbb{B}_{\|.\|} : a_{21} \ne 0\} \quad \text{and so,} \quad B_*=\mathbb{B}_{\|.\|} \setminus \{(a_{ij})_{i, j=1}^2 \in M_2(\C) : a_{21}=0\}.
\]
Since $\mathbb{B}_{\|.\|}$ is a domain, it follows from Proposition 1.1 in \cite{FG} that $B_*$ is a domain in $M_2(\C)$. Note that $\pi(B_*)=\{\pi(A) : A \in B_*\}=\{(a, x_1, x_2, x_3) \in \QN: a \ne 0\}$. We have by Proposition \ref{intQN} that $int(\QN)=\pi(B_*)$ and so, $int(\QN)$ is connected. Therefore, $int(\QN)$ is a subdomain of $\Q$. 

\smallskip

It remains to show that $\overline{int(\QN)}=\CQN$. Again by Proposition \ref{intQN}, $int(\CQN)=int(\QN)=\{(a, x_1, x_2, x_3) \in \QN : a \ne 0\}$ and so, $\overline{int(\QN)} \subseteq \CQN$. Let $\alpha=(a, x_1, x_2, x_3) \in \CQN$. Then there is a sequence of elements $\al_n=(a^{(n)}, x_1^{(n)}, x_{2}^{(n)}, x_3^{(n)})$ in $\QN$ that converges to $\al$. Suppose that there are only finitely many $n \in \mathbb{N}$ with $a^{(n)}=0$. Then there exists some $N \in \mathbb{N}$ such that $a^{(n)} \neq 0 $ for $n \geq N$. Then $\al_n \in int(\QN)$ for $n \geq N$ and so, $\alpha \in \overline{int(\QN)}$. We now assume that there are infinitely many $n \in \mathbb{N}$ with $a^{(n)}= 0$. One can choose a subsequence $\{a^{(n_k)}\}$ with $a^{(n_k)}=0$ and so, $a=\underset{n_k \to \infty}{\lim}a^{(n_k)}=0$. By Theorem \ref{thm4.6}, $(x_1, x_2, x_3) \in \Ebar$ and $x_1x_2=x_3$. Then $\al=(0, x_1, x_2, x_1x_2)$. For $\epsilon_n =1\slash n$, define
$
\beta_n=(a_n, x_{1n}, x_{2n}, x_{3n})=(\epsilon_n, x_1(1-\epsilon_n), x_2(1-\epsilon_n), x_1x_2(1-\epsilon_n)^2)
$
for $n \in \mathbb{N}$. It is evident that $\underset{n \to \infty}{\lim} \beta_n=\alpha$ and $|x_{1n}|, |x_{2n}| <1$. Note that 
\begin{equation*}
	\begin{split}
		& \ |a_n|^{-2}\left(-|a_n|^4+|a_n|^2(1-|x_{1n}|^2-|x_{2n}|^2+|x_{3n}|^2)-|x_{1n}x_{2n}-x_{3n}|^2 \right)\\
		&=-\epsilon_n^2+(1-|x_{1}|^2(1-\epsilon_n)^2)(1-|x_{2}|^2(1-\epsilon_n)^2)\\
		& \geq -\epsilon_n^2+(1-(1-\epsilon_n)^2)^2 \quad [\text{as } |x_{1}|, |x_{2}| <1 ]\\
		&=\epsilon_n(1-\epsilon_n)(1-(1-\epsilon_n)^2+\epsilon_n)\\
		&>0.
	\end{split}
\end{equation*}
We have by Lemma \ref{lem4.3} and Theorem \ref{thm4.4} that $\be_n \in \QN$ for $n \in \mathbb{N}$. As the first component of each $\beta_n$ is non-zero, we have that $\beta_n \in int(\QN)$ and so, $\alpha \in \overline{int(\QN)}$. The proof is now complete. 
\end{proof} 

It has been discussed prior to Proposition \ref{prop_cintCQN} that $\HB \cup (\{0\} \times \E)$ is a domain in $\C^4$. However, the same result does not hold for its norm analogue $\QN \cup (\{0\} \times \E)$, as the next proposition shows. 

\begin{prop}
	The set $\QN \cup (\{0\} \times \E)$ is connected but not open.
\end{prop}

\begin{proof} It is evident from Theorem \ref{QN_not_open} that $\QN$ is connected. Moreover, $\{0\} \times \E$ is also connected in $\C^4$. It follows from Theorem \ref{thm4.4} that $\QN \cap (\{0\} \times \E)=\{(0, x_1, x_2, x_3) \in \C \times \E :  x_1x_2=x_3\}$. Consequently, $\QN \cup (\{0\} \times \E)$ is connected. We show that $(0, 0, 0, 0) \in \QN \cup (\{0\} \times \E)$ is not an interior point. On the contrary, assume that there exists $\epsilon \in (0,1)$ such that 
\[
B(0, \epsilon) \times B(0, \epsilon) \times B(0, \epsilon) \times B(0, \epsilon) \subseteq \QN \cup (\{0\} \times \E) \ \ \text{and so,} \ \ \left(\frac{\epsilon}{2}, 0, \frac{\epsilon}{2}, \frac{\epsilon}{2}\right) \in \QN \cup (\{0\} \times \E).
\]
Since $\epsilon \ne 0$, we have that $(\epsilon \slash 2, 0, \epsilon \slash 2, \epsilon \slash 2) \in \QN$. Let $(a, x_1, x_2, x_3)=(\epsilon \slash 2, 0, \epsilon \slash 2, \epsilon \slash 2)$. Then $\be=1-|x_1|^2-|x_2|^2+|x_3|^2=1$ and $w^2=x_1x_2-x_3=-\epsilon \slash 2$. We have by Theorem \ref{thm4.4} that 
\[
\be-\sqrt{\be^2-4|w|^4}<2|a|^2 \quad \text{and so}, \quad 1 -\sqrt{1-\epsilon^2}<\epsilon^2\slash 2 \, ,
\]
which is not possible for any $\epsilon \in (0,1)$. Therefore, $\QN \cup (\{0\} \times \E)$ is not open. 
\end{proof}

\section{The intersection of hexablock with $\mathbb{R}^4$}\label{realslice} 

\noindent In this section, we study the real hexablock $\Q \cap \mathbb{R}^4$ and explain the importance of the term ``hexablock" in describing $\Q$. Recall that 
\begin{align}\label{eqn_Q_describe}
	\Q \cap \mathbb{R}^4=\left\{(a, x_1, x_2, x_3) \in \mathbb{R} \times (\E \cap \mathbb{R}^3) :  \ \left|\frac{a\sqrt{(1-|z_1|^2)(1-|z_2|^2)}}{1-x_1z_1-x_2z_2+x_3z_1z_2}\right|<1 \ \ \text{for all} \ z_1, z_2 \in \D   \right\}.
\end{align}
The motivation to study $\Q \cap \mathbb{R}^4$ comes from the geometric description of $\Pe \cap \mathbb{R}^3$ in \cite{AglerIV}. The authors of \cite{AglerIV} referred to $\Pe \cap \mathbb{R}^3$ as the real pentablock, which is given by 
\[
\Pe \cap \mathbb{R}^3=\left\{(a, s, p) \in \mathbb{R}^3: (s, p) \in  \Gg \cap \mathbb{R}^2 \ \  \text{and} \ \ |a|< \left|1-\frac{\frac{1}{2}s^2\slash (1+p)}{1+\sqrt{1-(s\slash (1+p))^2}}\right|\right\}.
\]
In fact, it was proved in \cite{AglerIV} that $\mathbb{P} \cap \mathbb{R}^3$ is a convex body bounded by five faces, three of them flat and two curved (see Theorem 9.3 in \cite{AglerIV}). For this reason, the domain $\Pe$ is called the pentablock. We focus on the study of $\Q \cap \mathbb{R}^4$ in the same direction. To do so, we first refine the description of $\Q \cap \mathbb{R}^4$ as in \eqref{eqn_Q_describe} by capitalizing the results from Section \ref{psi_map}. To begin with, let $(x_1, x_2, x_3) \in \E \cap \mathbb{R}^3$. It follows from Corollary \ref{cor_sup} that 
\[
\sup_{(z_1, z_2)\in \DC^2}\left|\frac{\sqrt{(1-|z_1|^2)(1-|z_2|^2)}}{1-x_1z_1-x_2z_2+x_3z_1z_2}\right| 
=\left|\frac{\sqrt{(1-|z_1(x)|^2)(1-|z_2(x)|^2)}}{1-x_1z_1(x)-x_2z_2(x)+x_3z_1(x)z_2(x)}\right|,
\]
where $z_1(x), z_2(x)$ are points in $(-1, 1)$ given as in Corollary \ref{cor2.2}.
Consequently, we have that 
\[
\Q \cap \mathbb{R}^4=\left\{(a, x_1, x_2, x_3) \in \mathbb{R} \times (\E \cap \mathbb{R}^3) : \left|\displaystyle\frac{a\sqrt{(1-|z_1(x)|^2)(1-|z_2(x)|^2)}}{1-x_1z_1(x)-x_2z_2(x)+x_3z_1(x)z_2(x)}\right|<1  \right\}.
\]
It follows from Corollary \ref{cor2.2} that
\[
\sup_{z_1, z_2\in [-1, 1]}\frac{\sqrt{(1-|z_1|^2)(1-|z_2|^2)}}{1-x_1z_1-x_2z_2+x_3z_1z_2} 
=\frac{\sqrt{(1-|z_1(x)|^2)(1-|z_2(x)|^2)}}{1-x_1z_1(x)-x_2z_2(x)+x_3z_1(x)z_2(x)}.
\]
Since the left-hand side is at least $1$ (by taking $z_1=z_2=0$), it follows that
\begin{align}\label{eqn_K} 
	1 \leq \frac{\sqrt{(1-|z_1(x)|^2)(1-|z_2(x)|^2)}}{1-x_1z_1(x)-x_2z_2(x)+x_3z_1(x)z_2(x)}=\left|\frac{\sqrt{(1-|z_1(x)|^2)(1-|z_2(x)|^2)}}{1-x_1z_1(x)-x_2z_2(x)+x_3z_1(x)z_2(x)}\right|.
\end{align}
Consequently, we  arrive at the following description of \(\Q \cap \mathbb{R}^4\).
\begin{prop}\label{prop_real_hexa}
	The real hexablock is given by
	\[
	\Q \cap \mathbb{R}^4=\left\{(a, x_1, x_2, x_3) \in \mathbb{R} \times (\E \cap \mathbb{R}^3) : |a|<\frac{1-x_1z_1(x)-x_2z_2(x)+x_3z_1(x)z_2(x)}{\sqrt{(1-|z_1(x)|^2)(1-|z_2(x)|^2)}} \right\},
	\]
	where $z_1(x), z_2(x)$ are points in $(-1, 1)$ corresponding to $x=(x_1, x_2, x_3) \in \E$ as in Corollary \ref{cor2.2}.	
\end{prop}

Agler et al. \cite{AglerIV} proved that $\Gg \cap \mathbb{R}^2$ and $\Pe \cap \mathbb{R}^3$ are convex sets. An interested reader is referred to Section 9 in \cite{AglerIV} for further details. A similar result holds for the tetrablock. In fact, Abouhajar et al. \cite{Abouhajar} proved that $\mathbb{E} \cap \mathbb{R}^3$ is an open tetrahedron and so, a convex set. We conjecture that the intersection of the hexablock with $\mathbb{R}^4$ is also convex. However, verifying this appears to be challenging at this point, mainly because of elaborated and complicated computations. We present an equivalent criterion for determining if $\Q \cap \mathbb{R}^4$ is convex.  To see this, we need a few terminologies from the literature. 

\begin{defn}
	Let $C \subseteq \mathbb{R}^n$ be a convex set. A function $f: C \to \mathbb{R}$ is said to be \textit{concave} if 
	\[
	tf(x)+(1-t)f(y) \leq  f(tx+(1-t)y)
	\]
	for every $x, y \in C$ and $t \in (0,1)$. Also, a function $f: C \to \mathbb{R}$ is said to be \textit{convex} if 
	\[
	tf(x)+(1-t)f(y) \geq  f(tx+(1-t)y)
	\]
	for every 	$x, y \in C$ and $t \in (0,1)$.
\end{defn}

\begin{prop}\label{prop_H_K}
	The set $\Q \cap \mathbb{R}^4$ is convex if and only if the map $K: \E \cap \mathbb{R}^3 \to \mathbb{R}$ given by
	\[
	K(x)=\frac{1-x_1z_1(x)-x_2z_2(x)+x_3z_1(x)z_2(x)}{\sqrt{(1-(z_1(x))^2)(1-(z_2(x))^2)}} 
	\] 
	is concave, where $z_1(x), z_2(x)$ are the points in the open interval $(-1, 1)$ corresponding to $x=(x_1, x_2, x_3) \in \E$ as in Corollary \ref{cor2.2}.
\end{prop}

\begin{proof}
Assume that $\Q \cap \mathbb{R}^4$ is convex. Let $x=(x_1, x_2, x_3), y=(y_1, y_2, y_3) \in \E \cap \mathbb{R}^3$ and let $t \in (0, 1)$. For $\epsilon \in (0, 1)$, we define $a_\epsilon(x)=\epsilon K(x)$ and $b_\epsilon(y)=\epsilon K(y)$. By Proposition \ref{prop_real_hexa}, we have that $(a_\epsilon(x), x_1, x_2, x_3), (b_\epsilon(y), y_1, y_2, y_3) \in \Q \cap \mathbb{R}^4$, because $|a_\epsilon(x)K(x)^{-1}|=|b_\epsilon(y)K(y)^{-1}|=\epsilon$. Since $\Q \cap \mathbb{R}^4$ is convex, it follows that 
\[
\left(ta_\epsilon(x)+(1-t)b_\epsilon(y), \ tx_1+(1-t)y_1, tx_2+(1-t)y_2, \ tx_3+(1-t)y_3\right) \in \Q \cap \mathbb{R}^4.
\]
For the sake of brevity, let $tx+(1-t)y=(tx_1+(1-t)y_1, \ tx_2+(1-t)y_2, \ tx_3+(1-t)y_3)$. Since $\E \cap \mathbb{R}^3$ is convex, we have that $tx+(1-t)y \in \E \cap \mathbb{R}^3$. It follows from Proposition \ref{prop_real_hexa} that
\[
\epsilon(tK(x)+(1-t)K(y))=ta_\epsilon(x)+(1-t)b_\epsilon(y) < K(tx+(1-t)y)
\]
for all $\epsilon \in (0,1)$. Letting $\epsilon \to 1$, it follows that $tK(x)+(1-t)K(y) \leq  K(tx+(1-t)y)$ and so, $K$ is a concave function on $\E \cap \mathbb{R}^3$. Conversely, let $K$ be a concave function and choose any points $(a, x_1, x_2, x_3), (b, y_1, y_2, y_3) \in \Q \cap \mathbb{R}^4$. By Proposition \ref{prop_real_hexa}, $x=(x_1, x_2, x_3), y=(y_1, y_2, y_3) \in  \E \cap \mathbb{R}^3, |a|<K(x)$ and $|b|<K(y)$. By Theorem 2.8 in \cite{Abouhajar}, $\E \cap \mathbb{R}^3$ is convex and so, $tx+(1-t)y \in \mathbb E \cap \mathbb R^3$. Thus,
$
|ta+(1-t)b| \leq t|a|+(1-t)|b| <tK(x)+(1-t)K(y) \leq K(tx+(1-t)y)
$
for all $t \in (0,1)$ and so,  $\mathbb{H} \cap \mathbb{R}^4$ is convex. The proof is now complete. 
\end{proof} 

We recall from the literature a useful characterization of concave functions, which provides another equivalent criterion for the convexity of the real hexablock.

\begin{thm}[\cite{Mas}, Mathematical Appendix, Theorem M.C.2]\label{thm_concave}
	Let $f$ be a twice continuously differentiable real-valued function on an open convex set $C$ in $\mathbb{R}^n$. Then $f$ is concave if and only if its Hessian matrix
	\[
	\begin{pmatrix}
		\displaystyle \frac{\partial ^2 f}{\partial x_i \partial x_j}(x_1, \dotsc, x_n)
	\end{pmatrix}_{i, j=1}^n	
	\] 
	is negative semi-definite for every $(x_1, \dotsc, x_n) \in C$. 	
\end{thm}

Evidently, the function $K(x)$ as in Proposition \ref{prop_H_K} is a twice continuously differentiable function on the open convex set $\E \cap \mathbb{R}^3$. Therefore, Theorem \ref{thm_concave} tells us that the convexity of $\Q \cap \mathbb{R}^4$ can be determined by examining the Hessian of $K(x)$. A similar approach was employed in \cite{AglerIV} to establish the convexity of $\mathbb{P} \cap \mathbb{R}^3$. However, computing the partial derivatives involved in the Hessian matrix	
\[
\begin{pmatrix}
	\displaystyle \frac{\partial ^2 K}{\partial x_i \partial x_j}(x_1, x_2, x_3)
\end{pmatrix}_{i, j=1}^3
\]
is a long and laborious process due to the complexity of the expressions. For this reason, establishing the convexity of $\Q \cap \mathbb{R}^4$ remains challenging. Naturally, we consider the smallest convex set in $\mathbb{R}^4$ that contains $\Q \cap \mathbb{R}^4$, i.e., the intersection of all convex sets containing $\Q \cap \mathbb{R}^4$. Recall that the intersection of all the convex sets containing a given subset $S$ of $\mathbb{R}^n$ is called the \textit{convex hull} of $S$ and is denoted by $\text{conv}(S)$. Since the intersection of an arbitrary collection of convex sets is convex, we have that $\text{conv}(S)$ is the smallest convex set containing $S$. In fact, we have the following description of the convex hull.

\begin{thm}[\cite{Rock}, Theorem 2.3 \& \cite{Robinson}, Lemma 21.9]\label{thm_conv_hull}
	For any $S \subseteq \mathbb{R}^n$, the convex hull of $S$ is given by
	\[
	\text{conv}(S)=\left\{\overset{m}{\underset{i=1}{\sum}}\lambda_ix_i \ \bigg|\lambda_1 \geq 0, \dotsc, \lambda_m \geq 0, \ \lambda_1+\dotsc+\lambda_m=1, \ x_1, \dotsc, x_m \in S \ \ \text{and} \ \ m \in \mathbb{N} \right\}
	\]
\end{thm}
We discuss a few fundamental theorems about convex sets in $\mathbb{R}^n$. 

\begin{thm}[Carath\'{e}odory, \cite{Cara} \& \cite{Rock}, Theorem 17.1]\label{thm_cara}
	For any $S \subseteq \mathbb{R}^n$ and $x \in \text{conv}(S)$, there exists $(n+1)$  points $x_1, \dotsc, x_{n+1} \in S$ such that $x \in \text{conv}(\{x_1, \dotsc, x_{n+1}\})$.
\end{thm} 

It clearly follows from Theorem \ref{thm_conv_hull} that
\[
\text{conv}(\Q \cap \mathbb{R}^4)=\left\{\overset{m}{\underset{i=1}{\sum}}\lambda_i\alpha_i \ \bigg|\lambda_1 \geq 0, \dotsc, \lambda_m \geq 0, \ \ \overset{m}{\underset{i=1}{\sum}}\lambda_i=1, \ \ \alpha_1, \dotsc, \alpha_m \in \Q\cap \mathbb{R}^4, \ \ m \in \mathbb{N}\right\}.
\]
The next result gives a description of the convex hull of $\CQ \cap \mathbb{R}^4$.

\begin{prop}
	$\overline{\text{conv}(\Q \cap \mathbb{R}^4)}=\text{conv}(\CQ \cap \mathbb{R}^4)$.
\end{prop}

\begin{proof} Let $\alpha=(a, x_1, x_2, x_3) \in \overline{\text{conv}(\Q \cap \mathbb{R}^4)}$. Then there exists a sequence of elements 
\[
\alpha^{(m)}=(a^{(m)}, x_1^{(m)}, x_2^{(m)}, x_3^{(m)}) \in \text{conv}(\Q \cap \mathbb{R}^4)
\]
such that $\alpha^{(m)} \to \alpha$ as $m \to \infty$. It follows from Theorem \ref{thm_cara} that each $\alpha^{(m)}$ can be written as convex combination of at most $5$ points in $\Q\cap \mathbb{R}^4$, i.e.,
\[
\left(a^{(m)}, x_1^{(m)}, x_2^{(m)}, x_3^{(m)}\right)=\overset{5}{\underset{i=1}{\sum}}\lambda_{im}\left(a_i^{(m)}, \ x_{1i}^{(m)}, \ x_{2i}^{(m)}, \ x_{3i}^{(m)}\right),
\]
where $(a_i^{(m)}, \ x_{1i}^{(m)}, \ x_{2i}^{(m)}, \ x_{3i}^{(m)}) \in \Q \cap \mathbb{R}^4, \lambda_{im} \geq 0$ and $\overset{5}{\underset{i=1}{\sum}}\lambda_{im}=1$ for $m \in \mathbb{N}$. It follows from Theorem \ref{thm5.3} that $\Q \subseteq \DC^4$. Thus, $\{\lambda_{im}\}_{m=1}^\infty$ and $\left\{(a_i^{(m)}, \ x_{1i}^{(m)}, \ x_{2i}^{(m)}, \ x_{3i}^{(m)})\right\}_{m=1}^\infty$ are finitely many bounded sequences for $1 \leq i \leq 5$. So, there exist convergent subsequences
\[
\{\lambda_{im_k}\}_{k=1}^\infty \quad \text{and} \quad \left\{(a_i^{(m_k)}, \ x_{1i}^{(m_k)}, \ x_{2i}^{(m_k)}, \ x_{3i}^{(m_k)})\right\}_{k=1}^\infty
\]
with the limits $\lambda_i$ and $(a_1, x_{1i}, x_{2i}, x_{3i})$, respectively, for $1 \leq i \leq 5$. It is clear that  $(a_i, x_{1i}, x_{2i}, x_{3i}) \in \overline{\Q \cap \mathbb{R}^4}=\CQ \cap \mathbb{R}^4$ for every $i=1, \dotsc, 5$. Moreover, $\lambda_i \geq 0$ for $1 \leq i \leq 5$ and $\overset{5}{\underset{i=1}{\sum}}\lambda_{i}=1$. By uniqueness of the limits, we have that 
\begin{align*}
	\alpha=(a, x_1, x_2, x_3)
	=\underset{k \to \infty}{\lim}\left(a^{(m_k)}, \ x_{1}^{(m_k)}, \ x_{2}^{(m_k)}, \ x_{3}^{(m_k)}\right)
	&=\lim_{k \to \infty}\overset{5}{\underset{i=1}{\sum}}\lambda_{im_k}\left(a_i^{(m_k)}, \ x_{1i}^{(m_k)}, \ x_{2i}^{(m_k)}, \ x_{3i}^{(m_k)}\right)\\
	&=\overset{5}{\underset{i=1}{\sum}}\lambda_{i}\left(a_i, \ x_{1i}, \ x_{2i}, \ x_{3i}\right).
\end{align*}
Hence, $\alpha \in \text{conv}(\CQ \cap \mathbb{R}^4)$ and so, $\overline{\text{conv}(\Q \cap \mathbb{R}^4)} \subseteq \text{conv}(\CQ \cap \mathbb{R}^4)$. We now prove that these two sets are equal. To do so, let $\alpha=(a, x_1, x_2, x_3) \in \text{conv}(\CQ \cap \mathbb{R}^4)$. Then
\begin{align*}
	(a, x_1, x_2, x_3)=\overset{m}{\underset{i=1}{\sum}}\lambda_{i}\left(a_i, \ x_{1i}, \ x_{2i}, \ x_{3i}\right),
\end{align*}
where $m \in \mathbb{N}, (a_i, \ x_{1i}, \ x_{2i}, \ x_{3i}) \in \CQ \cap \mathbb{R}^4, \lambda_i \geq 0$ and $\lambda_1+\dotsc+\lambda_m=1$ for $1 \leq i \leq m$. Let $r \in (0,1)$. By Corollary \ref{cor_CQ_r_Q}, $(ra_i, \ rx_{1i}, \ rx_{2i}, \ r^2x_{3i}) \in \Q \cap \mathbb{R}^4$ for $1 \leq i \leq m$. Consequently, we have that
\begin{align*}
	(ra, rx_1, rx_2, r^2x_3)=\overset{m}{\underset{i=1}{\sum}}\lambda_{i}\left(ra_i, \ rx_{1i}, \ rx_{2i}, \ r^2x_{3i}\right) \in \text{conv}(\Q \cap \mathbb{R}^4)
\end{align*}
for every $r \in (0,1)$ and so, $(a, x_1, x_2, x_3) \in \overline{\text{conv}(\Q \cap \mathbb{R}^4)}$. The proof is now complete. 
\end{proof}  

We now recall from \cite{AglerIV} the geometric description of the real pentablock. Before doing this, we also recall from the literature a few terminologies that shall be used frequently.

\begin{defn}[\cite{Rudin_FA}, Chapter 3 \& \cite{Voigt}, Chapter 17]\label{defn_extreme}
	Let $K$ be a non-empty subset of $\mathbb{R}^n$. A non-empty set $S \subseteq K$ is said to be an \textit{extreme set} (or an \textit{extreme subset}) of $K$ if $x \in K, y \in K, 0<t<1$ are such that $(1-t)x+ty \in S$, then $x 
	\in S$ and $y \in S$. 
\end{defn}

\begin{defn}[\cite{Rock}, Section 18 \& \cite{Voigt}, Chapter 17]\label{defn_face}
	Let $C \subseteq \mathbb{R}^n$ be a non-empty convex set. A non-empty convex subset $F$ of $C$ is said to be a \textit{face} of $C$ if, for any $x, y \in C$ and any real number $t \in (0,1)$ such that $(1-t)x + t y \in F$, it follows that $x \in F$ and $y \in F$. Equivalently, a face of $C$ is a convex extreme subset of $C$.
\end{defn}

We mention a very simple but interesting property of a face of a convex set, which follows directly from its definition. 

\begin{lem}\label{lem_face}
	Let $F$ be a face of a convex set $C \subseteq \mathbb{R}^n$. If $x_1, \dotsc, x_n \in C$, $\lambda_1, \dotsc, \lambda_n \in (0,1)$ with 
	$\overset{n}{\underset{i=1}{\sum}}\lambda_i=1$ and $\overset{n}{\underset{i=1}{\sum}}\lambda_ix_i \in F$, then $x_i \in F$ for all $i=1, \dotsc, n$. 
\end{lem}

\subsection*{A note on the faces of the real pentablock} Agler et al. \cite{AglerIV} proved that $\Pe \cap \mathbb{R}^3$ is a convex set bounded by the following five sets: 
\begin{enumerate}
	\item  the triangle $\mathfrak{T}_1$ with vertices $(0, 2, 1), (1, 0, -1)$ and $(-1, 0, -1)$ together with its interior;
	\item  the triangle  $\mathfrak{T}_2$ with vertices $(0, -2, 1), (1, 0, -1)$ and $(-1, 0, -1)$ together with its interior;
	\item the ellipse $\mathfrak{E}= \{(a, s, 1): a^2+s^2\slash 4 \leq 1, -2 \leq s \leq 2 \}$ together with its interior;
	\item the surface 
	\begin{align}\label{eqn_S1}
		\mathfrak{S}_1=\left\{(a, s, p) \in \mathbb{R} \times (\Gg \cap \mathbb{R}^2) :  0 \leq a \leq 1, \ a=\left|1-\frac{\frac{1}{2}s^2\slash (1+p)}{1+\sqrt{1-(s\slash (1+p))^2}}\right|\right\};
	\end{align}
	\item the surface 
	\begin{align}\label{eqn_S2}
		\mathfrak{S}_2=\left\{(a, s, p) \in \mathbb{R} \times (\Gg \cap \mathbb{R}^2) :  -1 \leq a \leq 0, \ a=-\left|1-\frac{\frac{1}{2}s^2\slash (1+p)}{1+\sqrt{1-(s\slash (1+p))^2}}\right|\right\}.
	\end{align}
\end{enumerate}
Since the closure of a convex set is again convex and so, $\overline{\Pe \cap \mathbb{R}^3}=\Pbar \cap \mathbb{R}^3$ is convex. It is evident that the triangles $\mathfrak{T}_1, \mathfrak{T}_2$ and the ellipse $\mathfrak{E}$ are convex sets contained in $\Pbar \cap \mathbb{R}^3$. Moreover, it is not difficult to see that $\mathfrak{T}_1, \mathfrak{T}_2$ and $\mathfrak{E}$ are extreme subsets of $\Pbar \cap \mathbb{R}^3$. Hence, each of these three sets are faces of $\Pbar \cap \mathbb{R}^3$. However, our next results show that the surfaces $\mathfrak{S}_1$ and $\mathfrak{S}_2$ do not satisfy the criteria to be faces of $\Pbar \cap \mathbb{R}^3$ in the sense of Definition \ref{defn_face}.

\begin{lem}\label{lem_S1}
	The set $\mathfrak{S}_1$ as in \eqref{eqn_S1} is not a face of $\Pbar \cap \mathbb{R}^3$. Moreover, $\mathfrak{S}_1$ is neither convex nor an extreme subset of $\Pbar \cap \mathbb{R}^3$.
\end{lem}

\begin{proof}
Let $\displaystyle (a, s, p)=\left(\frac{2+\sqrt{5}}{3+\sqrt{5}}, 1, \frac{1}{2}\right)$. Since $|s-\overline{s}p|<1-|p|^2$, it follows from Theorem \ref{thmG_2} that $(s, p) \in \Gg \cap \mathbb{R}^2$. Moreover, a routine calculation shows that
\begin{align*}
	\left|1-\frac{\frac{1}{2}s^2\slash (1+p)}{1+\sqrt{1-(s\slash (1+p))^2}}\right|
	=\left|1-\frac{s^2}{2\left((1+p)+\sqrt{(1+p)^2-s^2}\right)}\right|
	=\frac{2+\sqrt{5}}{3+\sqrt{5}}=a.
\end{align*}
Consequently, we have by \eqref{eqn_S1} that $(a, s, p)$ and $(a, -s, p)$ belong to $\mathfrak{S}_1$. Let us consider the mid point of $(a, s, p)$ and $(a, -s, p)$ which is given by
\[
(a_0, s_0, p_0)=\frac{1}{2}(a, s, p)+\frac{1}{2}(a, -s, p)=(a, 0, p).
\] 
Then
\begin{align*}
	\left|1-\frac{\frac{1}{2}s_0^2\slash (1+p_0)}{1+\sqrt{1-(s_0\slash (1+p_0))^2}}\right|
	=1 \ne \frac{2+\sqrt{5}}{3+\sqrt{5}}=a_0 \quad \text{and so}, \quad (a_0, s_0, p_0) \notin \mathfrak{S}_1.
\end{align*}
Therefore, $\mathfrak{S}_1$ is not convex. We also show that $\mathfrak{S}_1$ is not an extreme set of $\Pbar \cap \mathbb{R}^3$. Consider 
\[
(a_1, s_1, p_1)=(1, 0, 0), \quad (a_2, s_2, p_2)=(1, 0, 1) \quad \text{and} \quad (a_3, s_3, p_3)=(1, 0, 1\slash 2).
\]
Recall from \eqref{eqn_pi_Pe} that $\pi_\Pe(A)=(a_{21}, \text{tr}(A), \det(A))$ for $A=(a_{ij})_{i,j=1}^2 \in M_2(\C)$. Let us define
\[
A_1=\begin{pmatrix}
	0 & 0\\
	1 & 0
\end{pmatrix}, \quad A_2=\begin{pmatrix}
	0 & -1\\
	1 & 0
\end{pmatrix} \quad \text{and} \quad A_3=\begin{pmatrix}
	0 & -1\slash 2\\
	1 & 0
\end{pmatrix}.
\]
We have that $\pi_{\Pe}(A_1)=(1, 0, 0), \pi_{\Pe}(A_2)=(1, 0, 1)$ and $\pi_{\Pe}(A_3)=(1, 0, 1\slash 2)$. Since $\|A_j\| \leq 1$, we have by part (2) of Theorem \ref{pentablock_c} that $(a_j, s_j, p_j)$ belongs to $\Pbar \cap \mathbb{R}^3$ for $j=1, 2, 3$. By Theorem \ref{thmG_2}, $(s_2, p_2)=(0, 1) \notin \mathbb{G}_2$ and so, $(a_2, s_2, p_2) \notin \mathfrak{S}_1$. The point $(a_3, s_3, p_3)=(1, 0, 1\slash 2) \in \mathfrak{S}_1$, because 
\[
(0, 1\slash 2) \in \mathbb{G}_2 \cap \mathbb{R}^2 \quad \text{and} \quad 
\left|1-\frac{\frac{1}{2}s_3^2\slash (1+p_3)}{1+\sqrt{1-(s_3\slash (1+p_3))^2}}\right|
=a_3=1. 
\]
Consequently, we have that 
\[
(a_1, s_1, p_1), (a_2, s_2, p_2) \in \Pbar \cap \mathbb{R}^3, \quad (a_3, s_3, p_3)=\frac{1}{2}(a_1, s_1, p_1)+\frac{1}{2}(a_2, s_2, p_2)  \in \mathfrak{S}_1.
\]
However, $(a_2, s_2, p_2) \notin \mathfrak{S}_1$ and so, $\mathfrak{S}_1$ is not an extreme set of $\Pbar \cap \mathbb{R}^3$.  Therefore, $\mathfrak{S}_1$ is not a face of $\Pbar \cap \mathbb{R}^3$. 
\end{proof}  

\begin{lem}\label{lem_S2}
	The set $\mathfrak{S}_2$ as in \eqref{eqn_S2} is not a face of $\Pbar \cap \mathbb{R}^3$. Moreover, $\mathfrak{S}_2$ is neither convex nor an extreme subset of $\Pbar \cap \mathbb{R}^3$.
\end{lem}

\begin{proof} It follows from the definition of $\mathfrak{S}_1$ and $\mathfrak{S}_2$ as in \eqref{eqn_S1} and \eqref{eqn_S2} that $(a, s, p) \in \mathfrak{S}_1$ if and only if $(-a, s, p) \in \mathfrak{S}_2$. Since the map $\varphi: \C^3 \to \C^3$ defined as $\varphi(a, s, p)=(-a, s, p)$ preserves the convexity and $\varphi(\mathfrak{S}_2)=\mathfrak{S}_1$, we have by Lemma \ref{lem_S1} that $\mathfrak{S}_2$ is not convex.  We now show that $\mathfrak{S}_2$ is not an extreme set of $\Pbar \cap \mathbb{R}^3$. Consider 
\[
(a_1, s_1, p_1)=(-1, 0, 0), \quad (a_2, s_2, p_2)=(-1, 0, 1) \quad \text{and} \quad (a_3, s_3, p_3)=(-1, 0, 1\slash 2).
\]
Recall from \eqref{eqn_pi_Pe} that $\pi_\Pe(A)=(a_{21}, \text{tr}(A), \det(A))$ for $A=(a_{ij})_{i,j=1}^2 \in M_2(\C)$. Let us define
\[
A_1=\begin{pmatrix}
	0 & 0\\
	-1 & 0
\end{pmatrix}, \quad A_2=\begin{pmatrix}
	0 & 1\\
	-1 & 0
\end{pmatrix} \quad \text{and} \quad A_3=\begin{pmatrix}
	0 & 1\slash 2\\
	-1 & 0
\end{pmatrix}.
\]
We have that $\pi_{\Pe}(A_1)=(-1, 0, 0), \pi_{\Pe}(A_2)=(-1, 0, 1)$ and $\pi_{\Pe}(A_3)=(-1, 0, 1\slash 2)$. Since $\|A_j\| \leq 1$, we have by part (2) of Theorem \ref{pentablock_c} that $(a_j, s_j, p_j)$ belongs to $\Pbar \cap \mathbb{R}^3$ for $j=1, 2, 3$. By Theorem \ref{thmG_2}, $(s_2, p_2)=(0, 1) \notin  \mathbb{G}_2$ and so, $(a_2, s_2, p_2) \notin \mathfrak{S}_2$. Also, $(a_3, s_3, p_3)=(-1, 0, 1\slash 2) \in \mathfrak{S}_2$, because 
\[
(0, 1\slash 2) \in \mathbb{G}_2 \cap \mathbb{R}^2 \quad \text{and} \quad 
\left|1-\frac{\frac{1}{2}s_3^2\slash (1+p_3)}{1+\sqrt{1-(s_3\slash (1+p_3))^2}}\right|
=-a_3=1. 
\]
Thus, $(a_1, s_1, p_1), (a_2, s_2, p_2) \in \Pbar \cap \mathbb{R}^3, (a_3, s_3, p_3)=\frac{1}{2}(a_1, s_1, p_1)+\frac{1}{2}(a_2, s_2, p_2)  \in \mathfrak{S}_2$. However, $(a_2, s_2, p_2) \notin \mathfrak{S}_2$ and so, $\mathfrak{S}_2$ is not an extreme set of $\Pbar \cap \mathbb{R}^3$. Thus, $\mathfrak{S}_2$ is not a face of $\Pbar \cap \mathbb{R}^3$. 
\end{proof}  

\subsection*{The topological boundary of $\Q \cap \mathbb{R}^4$} It is evident from the discussion preceding Lemma \ref{lem_S1} regarding the description of $\mathbb{P} \cap \mathbb{R}^3$ (see Theorem 9.3 in \cite{AglerIV}) that the topological boundary of $\Pe \cap \mathbb{R}^3$ can be written as 
\[
\partial(\Pe \cap \mathbb{R}^3)=\mathfrak{T}_1 \cup \mathfrak{T}_2 \cup \mathfrak{E} \cup \mathfrak{S}_1 \cup \mathfrak{S}_2.
\]
For such a decomposition of $\partial(\Pe \cap \mathbb{R}^3)$ into the above five sets, the domain $\Pe$ is referred to as the pentablock. Employing similar techniques as in Theorem 9.3 of \cite{AglerIV}, we obtain a similar decomposition of $\partial(\Q \cap \mathbb{R}^4)$ into six sets, which gives rise to the name hexablock for the set $\Q$. To begin with, we have that 
\[
\partial (\Q \cap \mathbb{R}^4)=\overline{\Q \cap \mathbb{R}^4}\setminus (\Q \cap \mathbb{R}^4)=(\CQ \cap \mathbb{R}^4)\setminus (\Q \cap \mathbb{R}^4)=(\CQ \setminus \Q)\cap \mathbb{R}^4=\partial \Q \cap \mathbb{R}^4. 
\]
By Theorem \ref{thm5.3}, $(a, x_1, x_2, x_3) \in \partial \Q \cap \mathbb{R}^4$ if and only if at least one of the following holds.
\begin{enumerate}
	\item[(i)] $(a, x_1, x_2, x_3) \in \mathbb{R} \times (\partial\E \cap \mathbb{R}^3)$ and $|\Pz(a, x_1, x_2, x_3)| \leq 1$ for every $z_1, z_2 \in \D$;
	
	\item[(ii)] $(a, x_1, x_2, x_3) \in \mathbb{R} \times (\E \cap \mathbb{R}^3)$ and $ |\Pz(a, x_1, x_2, x_3)|=1$ for some $z_1, z_2 \in \D$.			
\end{enumerate}
Also, we have by Theorem \ref{thm5.3} that $\partial \CQ=\partial \Q$ and so, $\partial\CQ \cap\mathbb{R}^4=\partial \Q \cap\mathbb{R}^4$. Before proceeding further, we also discuss the geometric description of the real tetrablock $\E \cap \mathbb{R}^3$. We have by Theorem 2.8 in \cite{Abouhajar}, $\E \cap \mathbb{R}^3$ is an open tetrahedron with the vertices given by 
\[
A(1, 1, 1), \quad B(1, -1, -1), \quad C(-1, 1, -1) \quad  \text{and} \quad D(-1, -1, 1).
\] 
We show that the four faces of the open tetrahedron $\E \cap \mathbb{R}^3$ induce four sets $C_1, \dotsc, C_4$ contained in $\partial \CQ \cap \mathbb{R}^4$. In addition, the two faces $\mathfrak{S}_1, \mathfrak{S}_2$ of $\mathbb{P} \cap \mathbb{R}^3$ induce two more sets $C_5, C_6$ contained in $\partial \CQ \cap \mathbb{R}^4$ so that 
\[
\partial \CQ \cap \mathbb{R}^4=C_1 \cup C_2 \cup C_3 \cup C_4 \cup C_5 \cup C_6.
\]
Evidently, $\overline{\E \cap \mathbb{R}^3}=\Ebar \cap \mathbb{R}^3$ and so, $\partial(\E \cap \mathbb{R}^3)=\partial \E \cap \mathbb{R}^3$. Therefore, $\partial \E \cap \mathbb{R}^3$ consists of the following faces of tetrahedron determined by the vertices $A, B, C$ and $D$:
\begin{equation*}
	\begin{split}
		& (i) \ \Delta_1=\{(x_1, x_2, x_3) : -x_1+x_2-x_3+1 =0  \} \quad (iii) \ \Delta_3=\{(x_1, x_2, x_3) :  x_1+x_2+x_3+1 =0 \} \\
		& (ii) \ \Delta_2=\{(x_1, x_2, x_3) :  -x_1-x_2+x_3+1 =0  \} \quad (iv) \ \Delta_4=\{(x_1, x_2, x_3)  :  x_1-x_2-x_3+1 =0  \}.
	\end{split}
\end{equation*}
Naturally, these four faces of the convex set $\E \cap \mathbb R^3$ induce four closed sets in $\partial \CQ \cap \mathbb{R}^4$. Indeed, the sets given by
\begin{align}\label{eqn_Cj}
	C_j=\left\{(a, x_1, x_2, x_3) \in \mathbb{R} \times  \Delta_j :  |\Pz(a, x_1, x_2, x_3)| \leq 1 \ \text{for all} \ z_1, z_2 \in \D\right\} \ \ (1 \leq j \leq 4),
\end{align} 
are contained in $\partial \CQ \cap \mathbb{R}^4$. We show that each $C_j$ is closed. To do so, let $\{q_n\}$ be a sequence of elements in $C_j$ for $n \in \mathbb{N}$ and let $q=(a, x_1, x_2, x_3)$ be its limit. Since $C_j$ is contained in the closed set $\mathbb{R} \times \Delta_j$, we have that $q \in \mathbb{R} \times \Delta_j$. Also, by Theorem \ref{thm5.3}, $C_j \subseteq \CQ$ and so, $q \in \CQ$. Consequently, by Theorem \ref{thm5.3}, we have that $|\Pz(a, x_1, x_2, x_3)| \leq 1$ for every $z_1, z_2 \in \D$. Therefore, $q \in C_j$ and hence, we arrive at the following result. 

\begin{prop}
	The sets $C_1, \dotsc, C_4$ as in \eqref{eqn_Cj} are closed in $\mathbb{R}^4$. 	
\end{prop}

Since verifying the convexity of $C_1, \dotsc, C_4$ seems too difficult, we establish that these sets are extreme subsets of $\CQ \cap \mathbb{R}^4$. In addition, we prove that the convex hull of each of these sets is a face of the convex hull of $\CQ \cap \mathbb{R}^4$. 

\begin{prop}
	The sets $C_1, \dotsc, C_4$ as in \eqref{eqn_Cj} are extreme subsets of $\CQ \cap \mathbb{R}^4$. Moreover, the sets $\text{conv}(C_1), \dotsc, \text{conv}(C_4)$ are faces of $\text{conv}(\CQ \cap \mathbb{R}^4)$.
\end{prop}

\begin{proof}
Let $\alpha=(a, x_1, x_2, x_3), \beta=(b, y_1, y_2, y_3) \in \CQ \cap \mathbb{R}^4$ and let $t \in (0,1)$ be such that 
\[
t\alpha+(1-t)\beta=(ta+(1-t)b, tx_1+(1-t)y_1, tx_2+(1-t)y_2, tx_3+(1-t)y_3) \in C_1.
\]
Since $C_1 \subseteq \CQ$, we have by Theorem \ref{thm5.3} that $(x_1, x_2, x_3), (y_1, y_2, y_3) \in \Ebar \cap \mathbb{R}^3$, and 
\[
\left|\Pz\left(ta+(1-t)b, tx_1+(1-t)y_1, tx_2+(1-t)y_2, tx_3+(1-t)y_3\right)\right| \leq 1 \quad \text{for all} \ z_1, z_2 \in \D.
\] 
By definition of $C_1$, $(tx_1+(1-t)y_1, tx_2+(1-t)y_2, tx_3+(1-t)y_3) \in \Delta_1$. Since $\Delta_1$ is a face of the convex set $\Ebar \cap \mathbb{R}^3$ (see Theorem 2.8 in \cite{Abouhajar}), we have that $(x_1, x_2, x_3), (y_1, y_2, y_3) \in \Delta_1$. Hence, $C_1$ is an extreme subset of $\CQ\cap \mathbb{R}^4$. We now show that $\text{conv}(C_1)$ is a face of $\text{conv}(\CQ \cap \mathbb{R}^4)$. Let 
\[
\alpha=\overset{n}{\underset{i=1}{\sum}}\lambda_i \alpha_i \in \text{conv}(\CQ \cap \mathbb{R}^4) \quad \text{and} \quad \beta=\overset{n}{\underset{i=1}{\sum}}\mu_i \beta_i \in \text{conv}(\CQ \cap \mathbb{R}^4),
\]
where 
\[
\alpha_i=(a_i, x_{1i}, x_{2i}, x_{3i}), \ \beta_i=(b_i, y_{1i}, y_{2i}, y_{3i}) \in \CQ \cap \mathbb{R}^4
\] 
and $\lambda_i, \mu_i > 0$ for $1 \leq i \leq n$ with $\overset{n}{\underset{i=1}{\sum}}\lambda_i=\overset{n}{\underset{i=1}{\sum}}\mu_i=1$. Let  \[
t\alpha+(1-t)\beta=\overset{n}{\underset{i=1}{\sum}}\bigg[t\lambda_i\left(a, x_{1i}, x_{2i}, x_{3i}\right)+(1-t)\mu_i\left(b, y_{1i}, y_{2i}, y_{3i}\right)\bigg] \in \text{conv}(C_1)
\] 
for some $t \in (0,1)$. Then we can write
\[
t\alpha+(1-t)\beta=\overset{n}{\underset{i=1}{\sum}}\bigg[t\lambda_i\left(a, x_{1i}, x_{2i}, x_{3i}\right)+(1-t)\mu_i\left(b, y_{1i}, y_{2i}, y_{3i}\right)\bigg]=\overset{m}{\underset{j=1}{\sum}}\theta_j(c_j, w_{1j}, w_{2j}, w_{3j}),
\]
where $\theta_j > 0$ with $\overset{n}{\underset{j=1}{\sum}}\theta_j=1$ and $(c_j, w_{1j}, w_{2j}, w_{3j}) \in C_1$ for $1 \leq j \leq m$. Thus, we have that
\[
\overset{n}{\underset{i=1}{\sum}}\bigg[t\lambda_i\left(x_{1i}, x_{2i}, x_{3i}\right)+(1-t)\mu_i\left(y_{1i}, y_{2i}, y_{3i}\right)\bigg]=\overset{m}{\underset{j=1}{\sum}}\theta_j(w_{1j}, w_{2j}, w_{3j}) \in \text{conv}(\Delta_1)=\Delta_1,
\]
where the last equality follows because $\Delta_1$ is a face of $\Ebar \cap \mathbb{R}^3$. Since $(x_{1i}, x_{2i}, x_{3i}), (y_{1i}, y_{2i}, y_{3i}) \in \Ebar$, we have by Lemma \ref{lem_face} that  $(x_{1i}, x_{2i}, x_{3i}), (y_{1i}, y_{2i}, y_{3i}) \in \Delta_1$ for $1 \leq i \leq n$. By definition of $\Delta_1$, we have that $\alpha_i, \beta_i \in C_1$ for $1 \leq i \leq n$ and so, the points $\alpha, \beta \in \text{conv}(C_1)$. Consequently, $\text{conv}(C_1)$ is a face of $\text{conv}(\CQ \cap \mathbb{R}^4)$. Similarly, one can prove the desired conclusions for $C_2, C_3$ and $C_4$.
\end{proof} 

It is evident that $(\partial \CQ \cap \mathbb{R}^4) \setminus (C_1 \cup C_2 \cup C_3 \cup C_4)$ equals the set
\[
C_0=\left\{(a, x_1, x_2, x_3) \in \mathbb{R} \times  (\E \cap \mathbb{R}^3)  :  |\psi_{\al_1, \al_2}(a, x_1, x_2, x_3)|=1 \ \text{for some} \ (\al_1, \al_2) \in \D^2 \right\}.
\] 
We show that it is further possible to write $C_0$ as a union of two disjoint sets in a canonical way. Let $(a, x_1, x_2, x_3) \in C_0$ and let $x=(x_1, x_2, x_3)$. Then  $|\psi_{\al_1, \al_2}(a, x_1, x_2, x_3)|=1$ for some $\al_1, \al_2 \in \D$. By Theorem \ref{thm5.3},  $|\Pz(a, x_1, x_2, x_3)| \leq 1$ for all $(z_1, z_2) \in \D^2$ and so, 
$\underset{z_1, z_2 \in \D}{\sup}|\Pz(a, x_1, x_2, x_3)|=|\psi_{\al_1, \al_2}(a, x_1, x_2, x_3)|=1$. Since $x \in \E$, we have by Proposition \ref{prop2.1} that $(\al_1, \al_2)=(z_1(x), z_2(x)) \in \mathbb{R}^2$, where $z_1(x), z_2(x)$ the points in $\D$ associated with $x=(x_1, x_2, x_3)$ as in Proposition \ref{prop2.1}. Hence, $\al_1, \al_2 \in (-1, 1)$. So, we have that
\[
C_0=\left\{(a, x_1, x_2, x_3) \in \mathbb{R} \times  (\E \cap \mathbb{R}^3)  :  |\psi_{z_1(x), z_2(x)}(a, x_1, x_2, x_3)|=1 \right\}.
\] 
Since for $\al_1, \al_2 \in (-1,1)$ and $(a, x_1, x_2, x_3) \in \CQ \cap \mathbb{R}^4$, we have that $\psi_{\al_1, \al_2}(a, x_1, x_2, x_3) \in \mathbb{R}$. Consequently, the set $C_0$ can be written as the union of disjoint sets given by
\begin{align}\label{eqn_C5C_6}
	C_{5}=\left\{(a, x_1, x_2, x_3) \in \mathbb{R} \times (\E \cap \mathbb{R}^3) : 0 \leq a \leq 1, \   \frac{a\sqrt{(1-|z_1(x)|^2)(1-|z_2(x)|^2)}}{|1-x_1z_1(x)-x_2z_2(x)+x_3z_1(x)z_2(x)|}=1 \right\}, \ \ \notag \\
	C_{6}=\left\{(a, x_1, x_2, x_3) \in \mathbb{R} \times (\E \cap \mathbb{R}^3) : -1 \leq a \leq 0,  \frac{-a\sqrt{(1-|z_1(x)|^2)(1-|z_2(x)|^2)}}{|1-x_1z_1(x)-x_2z_2(x)+x_3z_1(x)z_2(x)|}=1 \right\}.
\end{align}
For the rest of this section, $C_5$ and $C_6$ always denote the above two sets as in \eqref{eqn_C5C_6}. The following result establishes a connection between $C_5$ and $\mathfrak{S}_1$, as well as between $ C_6$ and $\mathfrak{S}_2$.

\begin{lem}\label{lem_CS}
	An element $(a, s, p) \in \mathfrak{S}_1$ (respectively, $\mathfrak{S}_2$) if and only if $(a, s\slash2, s\slash2, p) \in C_5$ (respectively, $C_6$), where $\mathfrak{S}_1$ and $\mathfrak{S}_2$ are as in \eqref{eqn_S1} and \eqref{eqn_S2} respectively. 
\end{lem}

\begin{proof}  Let $x=(s\slash 2, s \slash 2, p) \in \E \cap \mathbb{R}^3$. We have by Proposition \ref{prop2.1} and Corollary \ref{cor_302} that
\[
z_1(x)=z_2(x)=\frac{s\slash (1+p)}{1+\sqrt{1-(s\slash (1+p))^2}}.
\]
A laborious but routine calculation gives
\begin{align*}
	\bigg|\displaystyle\frac{\sqrt{(1-|z_1(x)|^2)(1-|z_2(x)|^2)}}{1-(s\slash 2)z_1(x)-(s\slash 2)z_2(x)+pz_1(x)z_2(x)}\bigg|
	&=\left|1-\frac{\frac{1}{2}s^2\slash (1+p)}{1+\sqrt{1-(s\slash (1+p))^2}}\right|^{-1}.
\end{align*}
Consequently,  
\begin{equation}\label{eqn_CSII}
	\begin{split} 
		\bigg|\displaystyle\frac{a\sqrt{(1-|z_1(x)|^2)(1-|z_2(x)|^2)}}{1-(s\slash 2)z_1(x)-(s\slash 2)z_2(x)+pz_1(x)z_2(x)}\bigg|
		&=|a|\left|1-\frac{\frac{1}{2}s^2\slash (1+p)}{1+\sqrt{1-(s\slash (1+p))^2}}\right|^{-1}.
	\end{split}
\end{equation}
By part (2) of Theorem \ref{thmG_2} and part (5) of \ref{tetrablock}, $(s, p) \in \mathbb{G}_2$ if and only if $(s\slash2, s\slash 2, p) \in \E$. From this fact and \eqref{eqn_CSII}, the desired conclusion follows.
\end{proof} 

As discussed below \eqref{eqn_Cj}, the sets $C_1, \dotsc, C_4$ are closed sets in $\mathbb{R}^4$. However, $C_5$ and $C_6$ are not closed in $\mathbb{R}^4$. To see this, let us define $r_n=1-1\slash n$ and $q_n=(1, 0, 0, r_n)$ for $n \in \mathbb{N}$. By Theorem \ref{tetrablock}, we have that $x^{(n)}=(0, 0, r_n)$ is in $\E$ for every $n \in \mathbb{N}$. It follows from Proposition \ref{prop2.1} that $\left(z_1(x^{(n)}), z_2(x^{(n)})\right)=(0, 0)$ and so, $\psi_{z_1(x^{(n)}), z_2(x^{(n)})}(1, 0, 0, r_n)=1$. Hence, $\{q_n\}$ is a sequence of elements in $C_5$ with limit $(1, 0, 0, 1)$, which does not belong to $C_5$ since $(0, 0, 1) \notin \E$. Similarly, one can show that $p_n=(-1, 0, 0, r_n) \in C_6$ for $n \in \mathbb{N}$, and $\underset{n \to \infty}{\lim} p_n=(-1, 0, 0, 1) \notin C_6$. 

\begin{lem}
	The sets $C_5$ and $C_6$ are closed in $\mathbb{R} \times (\E \cap \mathbb{R}^3)$.
\end{lem}

\begin{proof}  Let $q_n=(a_{1n}, x_{1n}, x_{2n}, x_{3n}) \in C_5$ for $n \in \mathbb{N}$ and  let $q=(a, x_1, x_2, x_3) \in \mathbb{R} \times  (\E \cap \mathbb{R}^3)$ be the limit point of the sequence $\{q_n\}$. Assume that $x^{(n)}=(x_{1n}, x_{2n}, x_{3n})$ and $x=(x_1, x_2, x_3)$ for $n \in \mathbb{N}$. Since $q_n \in C_5$, we have that $|\psi_{\al_{1n}, \al_{2n}}(a_n, x_{1n}, x_{2n}, x_{3n})|=1$, where $(\alpha_{1n}, \alpha_{2n})=(z_1(x^{(n)}), z_2(x^{(n)}))$ (as in Proposition \ref{prop2.1}) for $n \in \mathbb{N}$. Note that 
\[
x^{(n)}=(x_{1n}, x_{2n}, x_{3n}) \in \E \quad \& \quad \sup_{z_1, z_2 \in \D}|\Pz(a_n, x_{1n}, x_{2n}, x_{3n})|=|\psi_{\al_{1n}, \al_{2n}}(a_n, x_{1n}, x_{2n}, x_{3n})|=1.
\] 
It is evident from Proposition \ref{prop2.1} that $\underset{n\to \infty}{\lim} \alpha_{1n}=z_1(x)$ and $\underset{n\to \infty}{\lim} \alpha_{2n}=z_2(x)$, where $(z_1(x), z_2(x))$ is the unique point in $\D^2$ at which $|\Pz(a, x_1, x_2, x_3)|$ attains its supremum over $\D^2$. Since $q_n \in C_5$, we have that $a_{1n} \in [0,1]$ and so, $a \in [0,1]$. By continuity arguments, we have that
\begin{align*}
	|\psi_{z_1(x), z_2(x)}(a, x_1, x_2, x_3)|
	=\lim_{n \to \infty} |\psi_{\al_{1n}, \al_{2n}}(a_n, x_{1n}, x_{2n}, x_{3n})|=1
\end{align*}
and so, $(a, x_1, x_2, x_3) \in C_5$. Therefore, $C_5$ is closed in $\mathbb{R} \times  (\E \cap \mathbb{R}^3)$. Similarly, one can show that $C_6$ is closed in $\mathbb{R} \times  (\E \cap \mathbb{R}^3)$. 
\end{proof}  

An immediate consequence of Lemma \ref{lem_S1} is the following result regarding the convexity of the subsets $C_5$ and $C_6$ of $\partial \Q \cap \mathbb{R}^4$.

\begin{cor}
	The sets $C_5$ and $C_6$ are not convex. Also, $C_5$ and $C_6$ are not extreme subsets of $\CQ \cap \mathbb{R}^4$. 
\end{cor}

\begin{proof} Let $\displaystyle (a, s, p)=\left(\frac{2+\sqrt{5}}{3+\sqrt{5}}, 1, \frac{1}{2}\right)$. Following the proof of Lemma \ref{lem_S1}, we have that 
\[
(a, s, p), (a, -s, p) \in \mathfrak{S}_1 \quad \text{and} \quad (a_0, s_0, p_0)=\frac{1}{2}(a, s, p)+\frac{1}{2}(a, -s, p)=(a, 0, p) \notin \mathfrak{S}_1. 
\]
Consequently, we have by Lemma \ref{lem_CS} that 
\[
\left(a, \frac{s}{2}, \frac{s}{2}, p\right), \left(a, \frac{-s}{2}, \frac{-s}{2}, p\right)  \in C_5 \ \ \text{and} \ \   \frac{1}{2}\left[\left(a, \frac{s}{2}, \frac{s}{2}, p\right)+ \left(a, \frac{-s}{2}, \frac{-s}{2}, p\right)\right]=  \left(a, 0, 0, p\right) \notin C_5.
\]
Hence, $C_5$ is not convex. Now, we consider the points 
\[
(a_1, s_1, p_1)=(1, 0, 0), \quad (a_2, s_2, p_2)=(1, 0, 1) \quad \text{and} \quad (a_3, s_3, p_3)=(1, 0, 1\slash 2).
\]
Following the proof of Lemma \ref{lem_S1}, we have that each of these three points is in $\Pbar \cap \mathbb{R}^3$. Moreover, $(a_3, s_3, p_3) \in \mathfrak{S}_1$ and $(a_2, s_2, p_2) \notin \mathfrak{S}_1$. By Theorem \ref{CHB}, it follows that
\[
\left(a_1, \frac{s_1}{2}, \frac{s_1}{2}, p_1\right), \left(a_2, \frac{s_2}{2}, \frac{s_2}{2}, p_2\right), \left(a_3, \frac{s_3}{2}, \frac{s_3}{2}, p_3\right) \in \CQ \cap \mathbb{R}^4.
\]
We have by Lemma \ref{lem_CS} that 
\[
\frac{1}{2}\left[\left(a_1, \frac{s_1}{2}, \frac{s_1}{2}, p_1\right)+ \left(a_2, \frac{s_2}{2}, \frac{s_2}{2}, p_2\right)\right]=\left(a_3, \frac{s_3}{2}, \frac{s_3}{2}, p_3\right) \in C_5 \ \ \text{and} \ \
\left(a_2, \frac{s_2}{2}, \frac{s_2}{2}, p_2\right) \notin C_5.
\]
Consequently, $C_5$ is not an extreme subset of $\CQ \cap \mathbb{R}^4$. The desired conclusions for $C_6$ follow in a similar manner to $C_5$, using the same techniques and applying Lemma \ref{lem_S2}. The proof is now complete. 
\end{proof}

Recall that for an open subset $\Omega$ in $\mathbb{R}^n$ and $k \in \mathbb{N}$, the class $C^k(\Omega)$ denotes the space of functions with continuous derivatives up to order $k$ on $\Omega$ and $C^\infty(\Omega)=\cap_{k=1}^\infty C^k(\Omega)$. A subset $S$ of $\mathbb{R}^n$ is called a \textit{hypersurface of class $C^k (1 \leq k \leq \infty)$} if for every $x_0 \in S$, there is an open set $U \subseteq \mathbb{R}^n$ containing $x_0$ and a real-valued function $\rho \in C^k(U)$ such that grad $\rho$ is non-vanishing on $S \cap U$ and $S\cap U=\{x \in U : \rho(x)=0\}$. An interested reader is referred to Chapter 0 in \cite{Folland} for further details. 

\begin{lem}
	The sets $C_5$ and $C_6$ are hypersurfaces of class $C^\infty$ in $\mathbb{R}^4$. 
\end{lem}

\begin{proof} Consider the map $K: \E \cap \mathbb{R}^3 \to \mathbb{R}$ given by
\[
K(x_1, x_2, x_3)=\frac{1-x_1z_1(x)-x_2z_2(x)+x_3z_1(x)z_2(x)}{\sqrt{(1-(z_1(x))^2)(1-(z_2(x))^2)}}, 
\] 
where $z_1(x), z_2(x)$ are the points in $(-1, 1)$ corresponding to $x=(x_1, x_2, x_3) \in \E$ as in Corollary \ref{cor2.2}. Again by Corollary \ref{cor2.2}, the map $K$ is in $C^\infty(\E \cap \mathbb{R}^3)$. For such $z_1(x), z_2(x)$, we have by \eqref{eqn_K} that $1-x_1z_1(x)-x_2z_2(x)+x_3z_1(x)z_2(x) >0$ and $0 < K(x) \leq 1$. Thus, it follows from \eqref{eqn_C5C_6} that
\begin{align*}
	& C_{5}=\left\{(a, x_1, x_2, x_3) \in \mathbb{R} \times (\E \cap \mathbb{R}^3) :   a=K(x_1, x_2, x_3) \right\} \quad \text{and} \\
	& C_{6}=\left\{(a, x_1, x_2, x_3) \in \mathbb{R} \times (\E \cap \mathbb{R}^3) :   a=-K(x_1, x_2, x_3) \right\}.
\end{align*}
Let $U=\mathbb{R} \times (\E \cap \mathbb{R}^3)$. Consider the maps $\rho-, \rho_+: U \to \mathbb{R}$ given by  $\rho_-(a, x_1, x_2, x_3)=a-K(x_1, x_2, x_3)$ and $\rho_+(a, x_1, x_2, x_3)=a+K(x_1, x_2, x_3)$. Evidently, the gradients of $\rho_-$ and $\rho_+$ are non-vanishing on the open set $U$. Also, we have that
\begin{align*}
	& C_5 \cap U=\{(a, x_1, x_2, x_3) \in U : \rho_-(a, x_1, x_2, x_3)=0\} \quad \text{and}  \\
	&C_6 \cap U=\{(a, x_1, x_2, x_3) \in U : \rho_+(a, x_1, x_2, x_3)=0\}.
\end{align*}
Since $\rho_-$ and $\rho_+$ are maps in  $C^\infty(U)$, we have the desired conclusion. 
\end{proof}

Evidently, combining all facts together, we conclude that 
\[
\partial \CQ \cap \mathbb{R}^4=C_1 \cup C_2 \cup C_3 \cup C_4 \cup C_{5} \cup C_{6},
\]
where $C_1, \dotsc, C_6$ are as in \eqref{eqn_Cj} and \eqref{eqn_C5C_6}. Since the topological boundary of $\mathbb{H} \cap \mathbb{R}^4$ decomposes into six such sets $C_1, \dotsc, C_6$, the domain $\Q$ is referred to as the hexablock. The four sets $C_1, \dotsc, C_4$ that are obtained from the four faces of the convex set $\E \cap \mathbb{R}^3$, are closed in $\mathbb{C}^4$. Moreover, these four sets are extreme subsets of $\CQ \cap \mathbb{R}^4$ and their convex hulls are faces of the convex hull of $\CQ \cap \mathbb{R}^4$. The remaining two mutually disjoint hypersurfaces $C_{5}$ and $C_{6}$ are closed in $\mathbb{R} \times  (\E \cap \mathbb{R}^3)$, and are induced by the two surfaces $\mathfrak{S}_1, \mathfrak{S}_2$ contained in $\partial(\mathbb{P} \cap \mathbb{R}^3)$ in the sense of Lemma \ref{lem_CS}.

\medskip

We conclude this article here. As a continuation of investigation in this direction, we shall discuss in a couple of upcoming papers about the automorphisms of $\Q$ and determine its distinguished boundary. Moreover, we shall characterize the rational inner functions into $\Q$ and prove a Schwarz type lemma for $\Q$. Deeper connection of the geometry and function theory of $\Q$ with that of the domains $\Gg$, $\E$, $\Pe$ will also be explored.

\subsection*{Funding} The first named author is partially supported by a J. C. Bose Fellowship (JBR\slash 2023\slash 000003). The second named author is supported by the `Core Research Grant (CRG)' with Award No. CRG/2023/005223 of Anusandhan National Research Foundation (ANRF), Govt. of India. The third named author was supported by the Prime Minister's Research Fellowship (PMRF) with PMRF ID No. 1300140 of the Govt. of India during the course of this work. At present the third named author is supported via the IIT Bombay RDF Grant of the second named author with Project Code RI/0115-10001427.

\end{document}